

\magnification 1200
\hsize 13.2cm
\vsize 20cm
\parskip 3pt plus 1pt
\parindent 5mm

\def\\{\hfil\break}


\font\seventeenbf=cmbx10 at 17.28pt

\font\twelvebf=cmbx10 at 12pt
\font\eightbf=cmbx8
\font\sixbf=cmbx6

\font\eighti=cmmi8
\font\sixi=cmmi6

\font\eightrm=cmr8
\font\sixrm=cmr6

\font\eightsy=cmsy8
\font\sixsy=cmsy6

\font\eightit=cmti8
\font\eighttt=cmtt8
\font\eightsl=cmsl8

\font\seventeenbsy=cmbsy10 at 17.28pt

\font\twelvebsy=cmbsy10 at 12pt
\font\tenbsy=cmbsy10
\font\eightbsy=cmbsy8
\font\sevenbsy=cmbsy7
\font\sixbsy=cmbsy6
\font\fivebsy=cmbsy5

\font\tenmsa=msam10

\font\sevenmsa=msam7
\font\fivemsa=msam5
\newfam\msafam
  \textfont\msafam=\tenmsa
  \scriptfont\msafam=\sevenmsa
  \scriptscriptfont\msafam=\fivemsa

\font\tenmsb=msbm10
\font\eightmsb=msbm8
\font\sevenmsb=msbm7
\font\fivemsb=msbm5
\newfam\msbfam
  \textfont\msbfam=\tenmsb
  \scriptfont\msbfam=\sevenmsb
  \scriptscriptfont\msbfam=\fivemsb
\def\Bbb{\fam\msbfam\tenmsb}

\font\tenCal=eusm10
\font\sevenCal=eusm7
\font\fiveCal=eusm5
\newfam\Calfam
  \textfont\Calfam=\tenCal
  \scriptfont\Calfam=\sevenCal
  \scriptscriptfont\Calfam=\fiveCal
\def\Cal{\fam\Calfam\tenCal}

\font\teneuf=eusm10
\font\teneuf=eufm10
\font\seveneuf=eufm7
\font\fiveeuf=eufm5
\newfam\euffam
  \textfont\euffam=\teneuf
  \scriptfont\euffam=\seveneuf
  \scriptscriptfont\euffam=\fiveeuf

\font\seventeenbfit=cmmib10 at 17.28pt

\font\twelvebfit=cmmib10 at 12pt
\font\tenbfit=cmmib10
\font\eightbfit=cmmib8
\font\sevenbfit=cmmib7
\font\sixbfit=cmmib6
\font\fivebfit=cmmib5
\newfam\bfitfam
  \textfont\bfitfam=\tenbfit
  \scriptfont\bfitfam=\sevenbfit
  \scriptscriptfont\bfitfam=\fivebfit


\catcode`\@=11
\def\eightpoint{%
  \textfont0=\eightrm \scriptfont0=\sixrm \scriptscriptfont0=\fiverm
  \def\rm{\fam\z@\eightrm}%
  \textfont1=\eighti \scriptfont1=\sixi \scriptscriptfont1=\fivei
  \def\oldstyle{\fam\@ne\eighti}%
  \textfont2=\eightsy \scriptfont2=\sixsy \scriptscriptfont2=\fivesy
  \textfont\itfam=\eightit
  \def\it{\fam\itfam\eightit}%
  \textfont\slfam=\eightsl
  \def\sl{\fam\slfam\eightsl}%
  \textfont\bffam=\eightbf \scriptfont\bffam=\sixbf
  \scriptscriptfont\bffam=\fivebf
  \def\bf{\fam\bffam\eightbf}%
  \textfont\ttfam=\eighttt
  \def\tt{\fam\ttfam\eighttt}%
  \textfont\msbfam=\eightmsb
  \def\Bbb{\fam\msbfam\eightmsb}%
  \abovedisplayskip=9pt plus 2pt minus 6pt
  \abovedisplayshortskip=0pt plus 2pt
  \belowdisplayskip=9pt plus 2pt minus 6pt
  \belowdisplayshortskip=5pt plus 2pt minus 3pt
  \smallskipamount=2pt plus 1pt minus 1pt
  \medskipamount=4pt plus 2pt minus 1pt
  \bigskipamount=9pt plus 3pt minus 3pt
  \normalbaselineskip=9pt
  \setbox\strutbox=\hbox{\vrule height7pt depth2pt width0pt}%
  \let\bigf@ntpc=\eightrm \let\smallf@ntpc=\sixrm
  \normalbaselines\rm}
\catcode`\@=12

\def\eightpointbf{%
 \textfont0=\eightbf   \scriptfont0=\sixbf   \scriptscriptfont0=\fivebf
 \textfont1=\eightbfit \scriptfont1=\sixbfit \scriptscriptfont1=\fivebfit
 \textfont2=\eightbsy  \scriptfont2=\sixbsy  \scriptscriptfont2=\fivebsy
 \eightbf
 \baselineskip=10pt}

\def\tenpointbf{%
 \textfont0=\tenbf   \scriptfont0=\sevenbf   \scriptscriptfont0=\fivebf
 \textfont1=\tenbfit \scriptfont1=\sevenbfit \scriptscriptfont1=\fivebfit
 \textfont2=\tenbsy  \scriptfont2=\sevenbsy  \scriptscriptfont2=\fivebsy
 \tenbf}

\def\twelvepointbf{%
 \textfont0=\twelvebf   \scriptfont0=\eightbf   \scriptscriptfont0=\sixbf
 \textfont1=\twelvebfit \scriptfont1=\eightbfit \scriptscriptfont1=\sixbfit
 \textfont2=\twelvebsy  \scriptfont2=\eightbsy  \scriptscriptfont2=\sixbsy
 \twelvebf
 \baselineskip=14.4pt}

\def\seventeenpointbf{%
 \textfont0=\seventeenbf  \scriptfont0=\twelvebf  \scriptscriptfont0=\eightbf
 \textfont1=\seventeenbfit\scriptfont1=\twelvebfit\scriptscriptfont1=\eightbfit
 \textfont2=\seventeenbsy \scriptfont2=\twelvebsy \scriptscriptfont2=\eightbsy
 \seventeenbf
 \baselineskip=20.736pt}


\newdimen\srdim \srdim=\hsize
\newdimen\irdim \irdim=\hsize
\def\NOSECTREF#1{\noindent\hbox to \srdim{\null\dotfill ???(#1)}}
\def\SECTREF#1{\noindent\hbox to \srdim{\csname REF\romannumeral#1\endcsname}}
\def\INDREF#1{\noindent\hbox to \irdim{\csname IND\romannumeral#1\endcsname}}
\newlinechar=`\^^J
\def\openauxfile{
  \immediate\openin1\jobname.aux
  \ifeof1
  \message{^^JCAUTION\string: you MUST run TeX a second time^^J}
  \let\sectref=\NOSECTREF \let\indref=\NOSECTREF
  \else
  \input \jobname.aux
  \message{^^JCAUTION\string: if the file has just been modified you may
    have to run TeX twice^^J}
  \let\sectref=\SECTREF \let\indref=\INDREF
  \fi
  \message{to get correct page numbers displayed in Contents or Index
    Tables^^J}
  \immediate\openout1=\jobname.aux
  \let\END=\end \def\end{\immediate\closeout1\END}}

\newbox\titlebox   \setbox\titlebox\hbox{\hfil}
\newbox\sectionbox \setbox\sectionbox\hbox{\hfil}
\def\folio{\ifnum\pageno=1 \hfil \else \ifodd\pageno
           \hfil {\eightpoint\copy\sectionbox\kern8mm\number\pageno}\else
           {\eightpoint\number\pageno\kern8mm\copy\titlebox}\hfil \fi\fi}
\footline={\hfil}
\headline={\folio}

\def\titlerunning#1{\setbox\titlebox\hbox{\eightpoint #1}}
\def\title#1{\noindent\hfil$\smash{\hbox{\seventeenpointbf #1}}$\hfil
             \titlerunning{#1}\medskip}

\newcount\numbersection \numbersection=-1
\def\sectionrunning#1{\setbox\sectionbox\hbox{\eightpoint #1}
  \immediate\write1{\string\def \string\REF
      \romannumeral\numbersection \string{%
      \noexpand#1 \string\dotfill \space \number\pageno \string}}}
\def\section#1{%
  \par\vskip0.666cm\penalty -100
  \vbox{\baselineskip=14.4pt\noindent{{\twelvepointbf #1}}}
  \vskip2pt
  \penalty 500
  \advance\numbersection by 1
  \sectionrunning{#1}}

\def\subsection#1{%
  \par\vskip0.5cm\penalty -100
  \vbox{\noindent{{\tenpointbf #1}}}
  \vskip1pt
  \penalty 500}

\newcount\numberindex \numberindex=0
\def\index#1#2{%
  \advance\numberindex by 1
  \immediate\write1{\string\def \string\IND #1%
     \romannumeral\numberindex \string{%
     \noexpand#2 \string\dotfill \space \string\S \number\numbersection,
     p.\string\ \space\number\pageno \string}}}

\newdimen\itemindent \itemindent=\parindent

\def\item#1{\par\noindent\hangindent\itemindent%
            \rlap{#1}\kern\itemindent\ignorespaces}
\def\itemitem#1{\par\noindent\hangindent2\itemindent%
            \kern\itemindent\rlap{#1}\kern\itemindent\ignorespaces}
\def\itemitemitem#1{\par\noindent\hangindent3\itemindent%
            \kern2\itemindent\rlap{#1}\kern\itemindent\ignorespaces}

\long\def\claim#1|#2\endclaim{\par\vskip 5pt\noindent
{\tenpointbf #1.}\ {\it #2}\par\vskip 5pt}

\def\proof{\noindent{\it Proof}}

\def\today{\ifcase\month\or
January\or February\or March\or April\or May\or June\or July\or August\or
September\or October\or November\or December\fi \space\number\day,
\number\year}

\catcode`\@=11
\newcount\@tempcnta \newcount\@tempcntb
\def\timeofday{{%
\@tempcnta=\time \divide\@tempcnta by 60 \@tempcntb=\@tempcnta
\multiply\@tempcntb by -60 \advance\@tempcntb by \time
\ifnum\@tempcntb > 9 \number\@tempcnta:\number\@tempcntb
  \else\number\@tempcnta:0\number\@tempcntb\fi}}
\catcode`\@=12

\def\bibitem#1&#2&#3&#4&%
{\hangindent=0.8cm\hangafter=1
\noindent\rlap{\hbox{\eightpointbf #1}}\kern0.8cm{\rm #2}{\it #3}{\rm #4.}}


\def\bQ{{\Bbb Q}}
\def\bR{{\Bbb R}}

\def\bZ{{\Bbb Z}}

\def\cA{{\Cal A}}
\def\cC{{\Cal C}}

\def\cF{{\Cal F}}

\def\cI{{\Cal I}}

\def\cO{{\Cal O}}
\def\cR{{\Cal R}}
\def\cS{{\Cal S}}


\def\square{{\hfill \hbox{
\vrule height 1.453ex  width 0.093ex  depth 0ex
\vrule height 1.5ex  width 1.3ex  depth -1.407ex\kern-0.1ex
\vrule height 1.453ex  width 0.093ex  depth 0ex\kern-1.35ex
\vrule height 0.093ex  width 1.3ex  depth 0ex}}}
\def\qed{\kern10pt$\square$}
\def\hexnbr#1{\ifnum#1<10 \number#1\else
 \ifnum#1=10 A\else\ifnum#1=11 B\else\ifnum#1=12 C\else
 \ifnum#1=13 D\else\ifnum#1=14 E\else\ifnum#1=15 F\fi\fi\fi\fi\fi\fi\fi}
\def\msatype{\hexnbr\msafam}
\def\msbtype{\hexnbr\msbfam}
\mathchardef\restriction="3\msatype16   
\mathchardef\boxsquare="3\msatype03
\mathchardef\preccurlyeq="3\msatype34
\mathchardef\compact="3\msatype62
\mathchardef\smallsetminus="2\msbtype72   
\mathchardef\subsetneq="3\msbtype28
\mathchardef\supsetneq="3\msbtype29
\mathchardef\leqslant="3\msatype36   
\mathchardef\geqslant="3\msatype3E   
\mathchardef\stimes="2\msatype02
\mathchardef\ltimes="2\msbtype6E
\mathchardef\rtimes="2\msbtype6F

\def\ddbar{\partial\overline\partial}

\let\ol=\overline

\let\wt=\widetilde
\let\wh=\widehat
\let\text=\hbox
\def\buildo#1^#2{\mathop{#1}\limits^{#2}}
\def\buildu#1_#2{\mathop{#1}\limits_{#2}}
\def\ort{\mathop{\hbox{\kern1pt\vrule width4.0pt height0.4pt depth0pt
    \vrule width0.4pt height6.0pt depth0pt\kern3.5pt}}}

\def\vlra{\mathrel{\smash-}\joinrel\mathrel{\smash-}\joinrel%
    \kern-2pt\longrightarrow}
\def\srelbar{\vrule width0.6ex height0.65ex depth-0.55ex}
\def\merto{\mathrel{\srelbar\kern1.3pt\srelbar\kern1.3pt\srelbar
    \kern1.3pt\srelbar\kern-1ex\raise0.28ex\hbox{${\scriptscriptstyle>}$}}}


\def\card{\mathop{\rm card}\nolimits}

\def\Supp{\mathop{\rm Supp}\nolimits}

\def\NS{\mathop{\rm NS}\nolimits}

\def\NS{{\rm NS}}

\def\nd{{\rm nd}}

\long\def\InsertFig#1 #2 #3 #4\EndFig{\par
\hbox{\hskip #1mm$\vbox to#2mm{\vfil\special{"
(/home/demailly/psinputs/grlib.ps) run
#3}}#4$}}
\long\def\LabelTeX#1 #2 #3\ELTX{\rlap{\kern#1mm\raise#2mm\hbox{#3}}}


\itemindent = 7mm

\title{Relative critical exponents,}
\title{non-vanishing and}
\title{metrics with minimal singularities}
\titlerunning{ }

\vskip10pt

\centerline {\tenrm Mihai P\u aun} 
\smallskip
\centerline {Institut \'Elie Cartan, Nancy}

\vskip20pt

\noindent{\bf Abstract. \it {In this article we prove a non-vanishing statement, as well as several properties of metrics with minimal singularities of adjoint bundles. Our arguments involve many ideas from Y.-T. Siu's analytic proof of the finite generation of the canonical ring. An important technical tool is 
the notion of relative critical exponent of two closed positive currents with respect to a measure.
}}

%
\section{\S0 Introduction}
\noindent The main theme of this article is the notion of {\sl relative critical exponent of
two closed positive currents with respect to a singular measure}. Before presenting the results we obtain in connection with this definition, we recall the general context in [4], which will be relevant for us.

Let $X$ be a compact complex manifold and
let $\Theta_1, \Theta_2$ be closed positive currents of (1,1)--type on $X$. We also consider a finite open covering $(U_\alpha)$ of $X$ and a family of functions 
$\psi:= (\psi_\alpha:= \psi_\alpha^1-\psi_\alpha^2)$ which are defined in terms of auxiliary  functions
$\psi_\alpha^j: U_\alpha\to [-\infty, \infty[$ and which have the properties :
\smallskip
{\itemindent 2mm
\item {\noindent $\bullet$} The difference $\displaystyle \psi_\alpha- \psi_\beta$ is non-singular 
on $U_\alpha\cap U_\beta$, and $\psi_\alpha^j$ is plurisubharmonic, for each $\alpha, \beta, j$ ;
\smallskip

\item {\noindent $\bullet$} We have $\int_{U_\alpha} \exp(-\psi_\alpha)d\lambda< \infty$, for each $\alpha$, where $d\lambda$ the Lebesgue measure.}

\medskip
\noindent Then we define
$$C_{\Theta_1, \psi}(X, \Theta_2):= \sup \{t\geq 0: \exp \big(t(\varphi_\alpha^{1}-
\varphi_\alpha^{2})-\psi_\alpha \big) \hbox { is in } L^1(U_\alpha),  \forall \hskip 0.3mm \alpha \}$$ 
where for each $j= 1,2$, 
the function $\varphi_\alpha^{j}$ in the preceding expression is a local potential of $\Theta_j$.

The normalization we use in the definition of the potentials above is such that if
$\Theta_2= [D], \Theta_1= 0$ and $\psi= 0$ (where $D$ is an effective $\bQ$-divisor on $X$), we recover the usual notion of
{\sl log-canonical threshold} of $D$ in  algebraic geometry. Thus, the quantity 
$\displaystyle C_{\Theta_1, \psi}(X, \Theta_2)$ is some kind of generalisation of that  notion in algebraic geometry. 

\smallskip

\noindent If $(\Theta_j, \psi)_{j= 1, 2}$ above have arbitrary singularities,  it seems very difficult to say something meaningful about $\displaystyle C_{\Theta_1, \psi}(X, \Theta_2)$, i.e.,
this notion is far too general to work with. For example, it is not clear whether $\displaystyle C_{\Theta_1, \psi}(X, \Theta_2)$ is  in general nonzero.
Another problem would be to define a notion of {\sl multiplier sheaf} in this context.

For these reasons, in the article [4] we are forced to confine ourselves most of the time to the 
{\sl logarithmic singularities} case. An important idea in the proof of the next result
is the consideration of a version of the threshold above which only takes into account a well chosen ``logarithmic part" of the preceding 
currents/functions.

\claim 0.1 Theorem|Let $X$ be a projective manifold, and let $\theta_L\in \NS_{\bR}(X)$ be a 
cohomology class in the real Neron-Severi space of $X$, such that: 
{\itemindent 5mm

\item {(a)} The adjoint class $c_1(K_X)+ \theta_L$ is pseudoeffective, i.e. there exist a 
closed positive current $$ \Theta_{K_X+ L}\in c_1(K_X)+ \theta_L ;$$

\item {(b)} The class $\theta_L$ contains a K\"ahler 
current $\Theta_L$ such that for any index $\alpha$ we have
$$e^{(1+\varepsilon_0)(\varphi_{K_X+ L}-\varphi_L)}\in L^1(U_\alpha)
\leqno (\star)$$
where $\varepsilon_0$ is a positive real number, and $\varphi_{K_X+ L}$ 
(resp. $\varphi_{L}$) is a local potential of the current $\displaystyle \Theta_{K_X+ L}$
(resp. $\displaystyle \Theta_{L}$).
}

\noindent Then the adjoint class $c_1(K_X)+ \theta_L$ contains an effective $\bR$--divisor, i.e. there exist a finite family of positive reals $\mu^{j}$ and hypersurfaces 
$W_j\subset X$ such that 
$$\sum_{j= 1}^N\mu^{j}[W_j]\in c_1(K_X)+ \theta_L.$$

\endclaim

\noindent We use the subscript ``$L$" in the statement 0.1 in order to suggest
that in some cases, $\theta_L$ is the Chern class of a line bundle $L$. \hfill \qed
\medskip

Firstly,  we would like to mention that the above result generalizes
the classical ``non-vanishing" theorems of  V. Shokurov and Y. Kawamata, cf. [23], [33], [34]. 

Secondly, it may seem that the integral hypothesis $(\star)$ in 0.1 is more 
general than the assumption that the critical exponent of $\Theta_L$ is greater than 1
(compare with [5]), but after  suitable modification of $X$ one can see that it is enough to consider this case. This is a consequence of a theorem due to H. Skoda, see e.g. the paragraph 1.B.1 and 1.C ; however, we prefer to state our result in this form, since the hypothesis $(b)$ is {\sl almost canonical}, in the sense that the quantity $(\star)$ with $\varepsilon_0= 0$ is a global measure on $X$. We would also like to mention that the result in [5]
is stronger : the authors obtain a statement within the {\sl numerical equivalence class} rather than in cohomology. 
 Of course, if 
$L$ is a $\bQ$--bundle and $\theta_L$ is its Chern class, then 0.1 imply the existence of a section of some multiple of $K_X+ L$ ; therefore, in the rational case we obtain the same non-vanishing theorem as in [5]. 

One important aspect of our proof is that is Char $p$-free, we avoid the {\sl explicit} use of 
the minimal model program algorithm.

A theorem similar to 0.1 was proved by Y.-T. Siu in [39], pages 31-46. Even if the hypothesis in his statement are much more restrictive than in the theorem 0.1, a substantial part of the arguments from his work will be used here. \hfill\qed 

It seems to us that 0.1 is not optimal : it is quite likely that one could prove a 
similar result with $\varepsilon_0= 0$. Also, we expect our result to hold under weaker positivity assumptions on $\Theta_L$ : the best result one could hope would be 
to replace the assumption ``$\Theta_L$ is a K\"ahler current" in $(b)$ above with 
the requirement $\Theta_L\geq 0$.
However, the difficulties one has to deal with in this case appear to be rather severe. 

Other possible generalization of 0.1 would be to work in {\sl purely transcendental} setting
i.e. without assuming that $\theta_L\in \NS_{\bR}(X)$, but so far it is not clear what  
the statement
should be. \hfill\qed

\bigskip
\noindent Along the next lines we will give a very rough overview of our arguments ; 
as a starting point, we recall the classical non-vanishing result due to V. Shokurov (see [34]).

\claim Theorem ([34])|Let $X$ be a projective manifold and let $D$, respectively $G= \sum_j\rho^jZ_j$
be a nef line bundle, resp. a $\bQ$-divisor. We assume that the following relations hold : 
{
\itemindent 5mm

\item {(i)} The $\bQ$--divisor $D+ G-K_X$ is nef and big ;
\smallskip
\item {(ii)} The critical exponent of the current associated to $-G$ is greater than 1.

}
\noindent Then for all large enough integers $m\in \bZ_+$, the bundle $mD+ \wh G$
is effective, where $\wh G$ is the ``round-up" of $G$.  \hfill\qed

\endclaim

\medskip
\noindent 
In order to explain the link between 0.1 and the above statement in a simple way, we assume for a moment the existence of a finite family of hypersurfaces 
$(Y_j)$ of $X$ with normal crossing intersections, such that
$$\Theta_L= \sum_{j\in J}a^{j}_L[Y_{j}]+ \Lambda_L\in  \theta_L,$$
where $\Lambda_L$ is a K\"ahler metric, $J$ is a finite set  and $a^{j}_L> 0$.
By the decomposition of closed positive currents theorem due to Y.-T. Siu
(see [36]) we have

$$\Theta_{K_X+ L}= \sum_{j\in J}a^{j}_{K_X+ L}[Y_{j}]+ \sum_{i\in I}x^i[W_i]+ \Lambda_{K_X+ L}\in c_1(K_X)+ \theta_L.$$
The numbers 
$(x^i, a^j_{K_X+L})$ in the formula above are positive reals and $(W_i)$ is a set of hypersurfaces of $X$ disjoint from $(Y_j)$. The cardinality of the index set $I$ could be infinite, and $\Lambda_{K_X+ L}$ is a closed positive current whose Lelong level sets have codimension at least 2.

Next, one can show that 
under the assumption $(\star)$ we have
$$a^{j}_{K_X+ L}- a^{j}_{L}> -1$$
(see the section 1.C of this article). On the other hand, by the equality above we infer the following relation :
$$\sum_{i\in I}x^i[W_i]+ \Lambda_{K_X+ L}+ \sum_{j\in J}(a^{j}_{K_X+ L}- a^j_L)[Y_{j}]
\in c_1(K_X+ \Lambda_L).$$
We introduce the following notations :

\noindent $\bullet$ $D^\prime:= \sum_ix^i[W_i]+ \Lambda_{K_X+ L}$ ;
\smallskip
\noindent $\bullet$ $G^\prime:= \sum_{j\in J}(a^{j}_{K_X+ L}- a^j_L)[Y_{j}]$

\noindent and the relation above shows that the Chern class of the $\bR$--bundle $D^{\prime}+
G^{\prime}-K_X$ contains a 
positive and non-singular representative, namely $\Lambda_L$. Moreover, since the coefficients of 
$G^{\prime}$ are greater than $-1$, the hypothesis $(ii)$ of the above theorem is verified.

We assume further that the current $D^\prime$ above correspond to a {\sl nef line bundle}, and that the coefficients $\displaystyle a^j_{K_X+ L}, a^j_L$ are {\sl rational}.
Under these additional hypothesis, the theorem 0.1 is a direct consequence of the aforementioned result of V. Shokurov, as follows.

For each  $j\in J$, we denote by $m^j$ the smallest integer greater than 
$\displaystyle a^j_{K_X+ L}- a^j_L$ (i.e., the {\sl round-up} of this difference). 
For any large enough integer $m\gg 0$, the theorem [34] show the existence of an effective $\bQ$-section
$$T_m\in c_1\Big(D^\prime+ {1\over m}\wh G^\prime\Big)$$
and we remark that 
$$T_m+ \sum_{j\in J}\Big(a^j_{K_X+ L}-{{m^j}\over {m}}\Big)Y_j\in c_1(K_X)+ \theta_L$$
is effective, by the definition of $(m^j)$.

\medskip
In this perspective, there are two important differences between 0.1 and the theorem of Shokurov.
To start with, the currents $D^\prime$ and $G^\prime$ do not necessarily correspond to a line bundle, 
respectively to a $\bQ$-line bundle. More seriously, the cohomology class of the current
$D^{\prime}$ defined by the first bullet above may not be numerically effective. Actually, $D^{\prime}$ has two components : a  
``divisor-like" part --although the sum in question could be infinite-- and a part 
(corresponding to $ \Lambda_{K_X+ L}$) which is
{\sl nef in codimension 1} in the terminology [7]. The first component admits a 
well-defined restriction to any of the hypersurfaces $(Y_j)$, and the same thing is true
for the second component, modulo a standard regularization process (see [7], [11]).

These excellent restriction properties of $D^\prime$ indicate that the proof of the theorem of Shokurov could eventually be adapted to the present setting. 
The principal use of the numerical effectivity of $D$ in [34] is to apply the vanishing theorems
during an inductive process, so that twisted pluricanonical sections defined on some well-chosen hypersurface of $X$ extend to the whole manifold.
In view of the very powerful extension theorems for the pluricanonical sections which were established since Y.-T. Siu's invariance of plurigenera breakthrough (see [37], [38]), we show that indeed, despite a few serious technical difficulties, the general outline of the proof of the classical non-vanishing result still works. \hfill\qed

 \medskip

Our article is organized as follows. 
We first review some standard facts  concerning the notion of {\sl numerical dimension} of a real $(1,1)$--class. If the dimension of $c_1(K_X)+ \theta_L$ is equal to zero, then the theorem 0.1 is a consequence of a result due to S. Boucksom, see [7]. If this is not the case, we use the numerical positivity of 
$c_1(K_X)+ \theta_L$, together with (a version of) the critical exponent defined above, in order to identify a hypersurface $S$ (the {\sl minimal center}) of some modification of $X$ such that by restriction to $S$ we reproduce the same context as in 0.1, except that the dimension drops
(see 1.C and 1.D). This part of our proof could be seen as a generalization of the classical arguments used in the Fujita conjecture literature (see [39] and the references therein).

During the restriction to the minimal center process, we will use in an essential 
manner the regularization techniques of Demailly ;
a diophantine approximation argument is also involved, to reduce to the case where the
geometric objects we are dealing with are rational (see 1.F, and also [5], [39]). Finally, we use the extension techniques of Siu and Hacon-McKernan adapted to the present situation (see 1.G and 1.H). The main technical point in 1.H is an {\sl ad hoc} version of the invariance of plurigenera.

The important steps in our proof of the main theorem have their origin in the notes [39], [40] by Siu ; at the beginning of each concerned paragraph we will make this more explicit.   
Most of the subtle points in our arguments are equally observable in the algebraic geometry proof [5], as it was kindly explained to us by J. McKernan and S. Druel ; it would be very interesting to have a precise comparison between the two approaches. 
\hfill\qed

\medskip

The second part of this article is a corollary of the first. 
Suppose given a $\bQ$-bundle $L$ which is big and endowed with a metric with
positive curvature current whose critical exponent is greater than 1. Assume further that
some multiple of the adjoint bundle $K_X+ L$  is effective.
In this context (in fact, in a slightly more general context, see  section 2), we want to compare the {\sl metric with minimal singularities} 
$\varphi_{\rm min}$ on the bundle $K_X+ L$ with 
its algebraic approximations, i.e. induced by {\sl finite} families of sections.

Very roughly, the main result we obtain is as follows. We assume that there exist an algebraic metric $\varphi_\alpha$ on $K_X+ L$ which is strictly more singular than $\varphi_{\rm min}$ ; then we obtain a modification
$\wh X\to X$ and a new algebraic metric $\varphi_{\alpha^\prime}$
such that 
$$\nu_S\big(\varphi_\alpha\circ \mu \big)> \nu_S(\varphi_{\rm min}\circ \mu)$$
and
$$\nu_S(\varphi_{\alpha^\prime}\circ \mu )= \nu_S(\varphi_{\rm min}\circ \mu)$$
for {\sl some} hypersurface $S\subset \wh X$.
The new metric is produced using  the non-vanishing theorem 0.1.\hfill\qed

\medskip


\medskip 

\claim Acknowledgments|It is a pleasure to thank S. Boucksom, B. Claudon and D. Varolin
for their very 
interesting and pertinent comments about the present text, as well as for suggesting 
infinitely many improvements/short-cuts. Also, we would like to mention that this paper
was completed during our visit to the  Mittag-Leffler Institute ; we are very grateful to the organizers for the invitation and support. We owe a debt of gratitude to J.-P. Demailly, S. Druel, L. Ein, J. McKernan and 
Y.-T. Siu, who shared with generosity and good humor their feelings about the topics in this paper. A substantial part of the techniques in this article emerged from our
collaboration with B. Berndtsson~; qu'il  en soit  chaleureusement remerci\'e!
\endclaim

\medskip
 
\vskip 10pt

\bigskip

\section {\S 1. A non-vanishing result}
\vskip 10pt

\noindent In this section we are going to prove the theorem 0.1, which is a version of the theorem  obtained in [5] as a by-product of their fundamental result on the finiteness of the canonical ring (see equally [16], [25] for interesting presentations of [5]). 
To start with, we give some precisions about the notions which were involved in the statement 0.1.

\noindent As it is well-known, an integral cohomology class in $H^2(X, \bZ)$
is the Chern class of a holomorphic line bundle if and only if it is of (1,1) type. The Neron-Severi group
$$\NS(X):= H^2(X, \bZ)\cap H^{1,1}(X, \bR)$$
is the set of cohomology classes of line bundles.
We denote by 
$$\NS_{\bR}(X):= \NS(X)\otimes_{\bZ}\bR\subset H^2(X, \bR)$$
the {\sl real Neron-Severi group}.
\smallskip
We also recall the following notions.
\claim Definition|A current $\Theta$ of type (1,1) is called a K\"ahler current if 
there exist a K\"ahler metric $\omega$ on $X$ such that
$\Theta\geq \omega.$ \hfill\qed
\endclaim

\medskip

\claim Definition|A function $\phi: X\to [-\infty, \infty[$ has logarithmic poles (or analytic singularities) if locally at each point $x\in X$ we have
$$\phi= \log(\sum_j|f_j|^2)$$
modulo $\cC^\infty$ functions, where $f_j\in \cO_{X, x}$ are local holomorphic functions.
 \hfill\qed
\endclaim

\medskip

\noindent Along the next lines we will  
use the {\sl proof of the  classical case} of the theorem 0.1 as a ``guide", together with the theory of closed positive currents and the invariance of plurigenera. We borrow a few techniques from both analytic as well as algebraic works on the subject~; however, we stress again that the characteristic $p$ methods, which seem to be essential in the later, are not used here. \hfill\qed

\bigskip

\noindent 
\subsection {\S 1.A Numerical dimension of pseudoeffective line bundles}
\smallskip

\noindent Let $X$ be a compact K\"ahler manifold endowed with a metric $\omega$ and let 
$\alpha$ be a non-singular (1,1)--form on $X$, which is assumed to be real and closed.
We denote its cohomology class by $\{\alpha\}\in H^{1,1}(X, \bR)$, and we assume this class to be pseudoeffective. 

Following [7], we denote by $\alpha[-\varepsilon\omega]$ the set of closed currents
$T\in \{\alpha\}$ with logarithmic poles and such that
$$T\geq -\varepsilon\omega.$$

The next fundamental result of J.-P. Demailly [11] is a quantitative version of the fact that if $\{\alpha\}$ is pseudoeffective, then for any $\varepsilon>0$ the set 
$\alpha[-\varepsilon\omega]$ is non empty.

\claim 1.A.1 Theorem([11])|Let $T= \alpha+ \sqrt{-1}\ddbar \varphi_T$ be a closed positive $(1, 1)$--current on a compact complex manifold $X$
(here the function $\varphi_T$ is globally defined on $X$). Then for any real number $\varepsilon> 0$
there exist a closed current 
$$T_\varepsilon:= \alpha+ \sqrt{-1}\ddbar \varphi_\varepsilon\in \alpha[-\varepsilon\omega]$$ 
such that we have the pointwise inequality
 $\varphi_\varepsilon\geq \varphi_T+ \cO (1)$.\hfill\qed

\endclaim
\medskip

\noindent In connection with his definition of mobile intersection of pseudoeffective classes, S. Boucksom proposed the next transcendental version of the classical notion of {\sl numerical dimension} of a nef line bundle, as follows :

$$\nd(\{\alpha\}):= \max \big\{k\in \bZ_+ : 
{\lim\sup}_{\varepsilon> 0, T_\varepsilon\in \alpha[-\varepsilon\omega]}\int_{X\setminus Z_\varepsilon} T_\varepsilon^k\wedge \omega^{n-k}> 0\big\}$$
where $Z_\varepsilon$ above is the singular set of the current $T_\varepsilon$.

If $\{\alpha\}= c_1(L)$ for some nef line bundle $L\to X$, then 
$\nd(\{\alpha\})$ above become the usual numerical dimension of $L$~; we refer to
[7] for a more detailed account about this notion and its properties.
\smallskip 
\noindent The statements which will follow assert the existence of geometric objects in the class
$\{\alpha\}$ and its approximations, according to the size of its numerical dimension.
The first one is due to S. Boucksom (see also the work of N. Nakayama, [31]).

\claim 1.A.2 Theorem ([7])|Let $\{\alpha\}$ be a $(1, 1)$--cohomology class which is 
pseudoeffective and such that $\nd(\{\alpha\})= 0$. Then there exist a closed positive current 
$$\Theta:= \sum_{j= 1}^\rho\nu^{j}[Y_j]\in \{\alpha\}.$$

\endclaim

\noindent For a more complete discussion about the properties of the current 
$\Theta$ above we refer to the article [7].\hfill\qed

\smallskip
\noindent Concerning the pseudoeffective classes $\{\alpha\}\in \NS_{\bR}(X)$ whose numerical dimension is strictly greater than 0, we have the following well-known statement.

\claim 1.A.3 Theorem|Let $X$ be a projective manifold, let $\{\alpha\}\in \NS_{\bR}(X)$ be a pseudoeffective class, such that 
$\nd (\{\alpha\})\geq 1$, and let $\beta$ be a K\"ahler current. Then for any $x\in X$ and $m\in \bZ_+$ there exist 
an integer $k_m$ and a closed positive current
$$T_{k, x}\in \{m\alpha+ \beta\}$$
with logarithmic poles, and such that $\nu (T_{m, x}, x)\geq k_m$ and 
$k_m\to\infty$ as $m\to\infty$. \hfill\qed
\endclaim

\proof. We fix an ample bundle $A\to X$, endowed with a metric $h$ with positive curvature. Let $N_0\in \bZ_+$ such that 
$$\beta \geq {{2}\over {N_0}}\Theta_{h}(A).\leqno (1)$$
Since $\nd (\{\alpha\})\geq 1$, there exist a positive constant $C> 0$ and a family of currents 
$\displaystyle T_\varepsilon \in \alpha[-\varepsilon\Theta_h(A)]$ such that
$$\int_{X\setminus Z_\varepsilon} T_\varepsilon\wedge \Theta_h(A)^{n-1}
\geq C> 0\leqno (2)$$
for any positive $\varepsilon$.

We will use now the hypothesis $\{\alpha\}\in \NS_{\bR}(X)$ : there exist 
a sequence of bundles $(L_m)_{m\in \bZ_+}$ such that 
$$\Vert c_1(L_m)-m\{\alpha\}\Vert\to 0\leqno (3)$$
as $m\to \infty$, where $\Vert \cdot\Vert$ denotes any norm on $\NS_\bR(X)$.
We are not going to explain the details of this claim, since it is a simple 
diophantine approximation argument, see e.g. [21], [5], [39], and the paragraph 1.F of this article, but rather indicate how to use the family of currents $T_\varepsilon$ above in order to obtain a lower bound of the quantity
$${{1}\over {p^n}}h^0\big(X, p(N_0L_m+ A)\big)$$
as $p\to \infty$. 

We recall that a very precise lower bound for the asymptotic behavior 
of the above dimension is provided by 
the holomorphic Morse inequalities, in the version obtained by L. Bonavero in [6]. 
In order to apply this result, we have to endow the bundle 
$N_0L_m+ A$ with a suitable metric. To this end, we remark that we have
$$N_0c_1(L_m)+ c_1(A)= N_0(c_1(L_m)- m\{\alpha\})+ mN_0\{\alpha\}+ c_1(A).\leqno (4)$$
The class $c_1(L_m)- m\{\alpha\}$ contains a non-singular representative $\rho_m$ which tend to zero as $m\to \infty$, by the relation (3). We take 
$$\varepsilon_m:= {{1}\over {2mN_0}}$$
and then we have 
$$\rho_m+ mN_0T_{\varepsilon_m}+ \Theta_h(A)\geq 0\leqno (5)$$
if $m\gg 0$. The relation (4) imply the existence of a metric $h_m$ on the bundle
$N_0L_m+ A$, whose associated curvature current is (5).

By the holomorphic Morse inequalities [6] we obtain
$${{1}\over {p^n}}h^0\big(X, p(N_0L_m+ A)\big)\geq C_0\int_{X\setminus Z_m}
\big(\rho_m+ mN_0T_{\varepsilon_m}+ \Theta_h(A)\big)^n
$$
where $Z_m$ is the set of singularities of $T_{\varepsilon_m}$ and $C_0$ is a positive constant, independent on $p$ and $m$. The inequality (2) show that the 
growth of the integral in the right hand side of the above inequality is at least linear with respect to $m$.

Now we invoke the usual linear algebra arguments (see [28]) and infer the existence of a $\bQ$--divisor 
$D_{m, x}\in  c_1(N_0L_m+ A)$ such that $\nu(D_{m, x}, x)\geq Cm^{1\over n}$
as $m\to\infty$. Then we define
$$T_{m, x}:= {{1}\over {N_0}}\big([D_{m, x}]- \Theta_h(A)\big)+ \beta-\rho_m ;$$
it is a closed, positive (1,1) current in the class $\{m\alpha+ \beta\}$, and its Lelong number at $x$ tend to infinity with $m$. Thus the statement 1.A.3 is completely proved. \hfill\qed

\medskip
\claim 1.A.4 Remark|{\rm The theorem 1.A.3 holds in a more general context, as follows.
\claim Theorem|Let $X$ be a compact K\"ahler manifold, 
endowed with a K\"ahler curent $\omega$, and let $\{\alpha\}$ be a pseudoeffective class, such that 
$\nd (\{\alpha\})\geq 1$. Then for any $x\in X$ and $m\in \bZ_+$ there exist 
an integer $k_m$ and a closed positive current
$$T_{m, x}\in \{m\alpha+ \omega\}$$
with logarithmic poles, and such that $\nu (T_{m, x}, x)\geq k_m$ and 
$k_m\to\infty$ as $m\to\infty$.

\endclaim

\smallskip The proof of this result will not be discussed here, since we do not need it. Let us just mention that the ``ancestor" of the above result
can be found in the beautiful article [12] ; see also [15] for an overview of the techniques involved in the proof (the Yau theorem [46], and of the maximum principle of Bedford-Taylor [1]).\hfill\qed

}\endclaim


\subsection {\S 1.B Dichotomy}

\smallskip After the preliminary discussion in the previous paragraph concerning the numerical dimension of the pseudoeffective classes and some of its properties, we start now the actual proof of the non-vanishing theorem. We denote by $\nu$ the numerical dimension of the 
class $c_1(K_X)+ \theta_L$, and we proceed as in [5], [23], [27], [39].
\smallskip
\noindent $\bullet$ If $\nu= 0$, then the existence of the $\bR$-section in the class 
$c_1(K_X)+ \theta_L$ is given by the theorem 1.A.2 above, therefore this first case is completely settled.\hfill\qed
\smallskip
\noindent $\bullet$ The second case $\nu\geq 1$ is much more involved~; we are going to use induction on the dimension of the manifold. 
In order to ease the comprehension of our approach, 
we will first prove the theorem 0.1 under some additional rationality and finiteness assumptions in the next subsection.


\subsection {\S 1.B.1 A particular case of 0.1}
\medskip

\noindent The aim of the present subsection is to give a detailed proof of the next result,
in order to illustrate the approach/difficulties for the general case.
\claim 1.B.1 Theorem|Let $X$ be a projective manifold and let $L\to X$ be a 
$\bQ$-line bundle ; we denote by 
$\theta_L\in \NS_\bR(X)$ its Chern class.
 We assume the existence of a closed positive 
current $\displaystyle \Theta_{K_X+ L}\in c_1(K_X)+ \theta_L$ with logarithmic poles and rational Lelong numbers. Moreover,
we assume that the class 
$\theta_L$ contains a K\"ahler current $\Theta_L$ such that
$$\int_X\exp(\varphi_{K_X+L}-\varphi_L)d\lambda< \infty.$$
Then $H^0\big(X, p(K_X+L)\big)\neq 0$ for all $p$ large and divisible enough. 
\endclaim

\smallskip
\noindent We remark that the integral condition above is less restrictive than
the hypothesis in the statement 0.1~; therefore, the (heavy) additional assumption 
is the existence of a current {\sl with log poles and rational singularities } in the class $\{K_X+ L\}$. \hfill\qed

\medskip
\proof. 
\noindent In the first place, we would like to mention that this version of the non-vanishing is {\sl almost} due to Shokurov, but we are going to prove it in a slightly different manner, which is better adapted for the illustration of the general case.

\smallskip
To start with, we remark that we can assume that the current 
$T_m$ in the statement 1.A.3 is given by an effective $\bQ$-section $D$
of the bundle $m(K_X+ L)+ L$. Indeed, as a consequence of the holomorphic Morse inequalities we have
$${{1}\over {p^n}}h^0\big(X, mp(K_X+ L)+ pL\big)\geq Cm^{\nu}$$
see e.g. [6]. Thus, the existence of the section $D$ above is provided by the usual
linear algebra arguments.
We fix $m$ large enough, so the the singularity of $\nu (D, x_0)\geq n+1$, where the point $x_0$ is chosen such that $\Theta_{K_X+L}$ and
$\Theta_L$ are non-singular at $x_0$.

\medskip

 Now let us consider the following relative critical exponent,
which is adapted to the current situation. We set
$$\tau:= C_{m\Theta_{K_X+L}+ \Theta_L, e^{\varphi_{K_X+L}-\varphi_L}}(X, D) ;$$
in other words we have
$$\tau = \sup\{t\in \bR_+:  \int_Xe^{t(\varphi_L+ m\varphi_{K_X+ L}- \varphi_D)}
e^{\varphi_{K_X+ L}- \varphi_L}d\lambda<\infty \}.\leqno (6)$$
\noindent Thus we consider the relative critical exponent {\sl with respect to the singular measure} 
of finite mass
$$e^{\varphi_{K_X+ L}- \varphi_L}d\lambda.$$

We note that we have the relations 
$$0< \tau< 1~;\leqno (7)$$
indeed, the first inequality is due to the fact that by hypothesis the relative critical exponent of $\Theta_L$ with respect to $\displaystyle \Theta_{K_X+ L}$
is greater than 1. As for the second one, it can be seen as a consequence of the fact that the 
singularity of $D$ et $x_0$ is large enough.

As in the proof of the non-vanishing result in [34],
we are going to use an appropriate modification of the manifold $X$ in order to have an interpretation of the quantity $\tau$.
\smallskip
\noindent  Let $\mu: \wh X\to X$ be a modification with such that the singular part of the inverse images of the currents below have normal crossing.

{\itemindent 7mm

\smallskip
$$\mu^\star\big(\Theta_{K_X+ L}\big)= \sum_{j\in J}a^{j}_{K_X+ L}[Y_{j}]+ \wh \Lambda_{K_X+ L}\leqno (8)$$

$$\mu^\star(\Theta_L)= \sum_{j\in J}a^{j}_L[Y_{j}]+ \wh \Lambda_L\leqno (9)$$

$$\mu^\star(\Theta_{D})= \sum_{j\in J}a^{j}_D[Y_{j}]
\leqno (10)$$

$$K_{\wh X/X}=  \sum_{j\in J}a^{j}_{\wh X/X}[Y_j]
\leqno (11)$$
where $J$ is a finite set, $(a^j)$ are non-negative real numbers (some of them may be zero, since we use the same family of indexes), $(Y_j)$ are divisors in $\wh X$, either proper transforms of divisors in $X$ or $\mu$-exceptional~; we will assume that they contain all the exceptional divisors of $\mu$ (as we can take the corresponding coefficients zero if necessary).
Finally $\wh \Lambda$
are non-singular, semi-positive (1, 1)--forms on $\wh X$. Moreover, we remark that
$$ \wh \Lambda_L\geq \mu^\star \omega.\leqno (12)$$
}

\smallskip

\noindent Let $\omega$ (respectively $\wh \omega$) be a K\"ahler metric on 
$X$ (respectively $\wh X$). The definition of our relative critical exponent (6) show that the 
quantity we have to evaluate is
$$\mu^\star \big(\Theta_\omega(K_X)+ t\big([D]- m\Theta_{K_X+ L}-\Theta_L\big)+ \Theta_L-\Theta_{K_X+ L}\big)$$
that is to say 
$$\mu^\star\big(\Theta_\omega(K_X)+ t[D]- \big(1+ tm\big)\Theta_{K_X+ L}+ (1-t)\Theta_L\big)\leqno (13)$$
where we denote by $[D]$ the current of integration associated to the 
$\bQ$-divisor $D$.

We first remark that {\sl the cohomology class of the current above is equal to zero}. This is an important point 
in our proof.

\medskip On the other hand, the equalities (8)--(11) show that for any positive real $t$ we have
$$ \leqno (14)
\eqalign {
\mu^\star \big(& \Theta_\omega(K_X) + t[D]+ (1-t)\Theta_{L}- (1+ mt)\Theta_{K_X+ L}\big)\equiv \cr
\equiv & \Theta_{\wh \omega}(K_{\wh X})+ \sum_{j\in J}(ta^{j}_D+ (1-t)a^{j}_L- (1+ mt)a^{j}_{K_X+L}-
a^j_{\wh X/X})[Y_{j}]+ \cr
+ & (1-t)\wh \Lambda_L- (1+ mt)\wh\Lambda_{K_X+L}\cr
}
$$
where the symbol $\equiv$ means that the two currents above have the same cohomology class.
Therefore, we get the next cohomological identity
$$(1+ mt)\wh\Lambda_{K_X+L}\equiv \Theta_{\wh \omega}(K_{\wh X})+ \sum_{j\in J}(ta^{j}_D+ (1-t)a^{j}_L- (1+ mt)a^{j}_{K_X+L}-
a^j_{\wh X/X})[Y_{j}]+ (1-t)\wh \Lambda_L.$$

\noindent By the definition (6), for any real $t<\tau$ and $j\in J$ we have the  
following inequalities :

\noindent $(\bullet)_j$ $t \big(a^{j}_D- a^{j}_L- ma^{j}_{K_X+L}\big)< 1+
a^{j}_{K_X+L}- a^{j}_{L}+ a^j_{\wh X/X}$.

\noindent In addition at least one of the previous inequality is 
in fact an equality for $t= \tau$. We will show next that we can further modify the inverse image of the current $\Theta_L$ in order to have equality 
for precisely on single index~; we stress on the fact that we can achieve this without changing 
the cohomology classes above, since $L$ is big.

To this end, we remark that given any family of positive rational numbers $\eta^{j}$, 
the $\bQ$--divisor
 $$\sum_{j\in J}(a^{j}_L+ \eta^{j})[Y_{j}] 
+ \wh\Lambda_{L,1}$$
belongs to the Chern class of the $\bQ$--line bundle 
$\mu^\star (L)$, provided that
$$\wh \Lambda_{L, 1}\equiv \wh \Lambda_L- \sum_{j\in J}\eta^{j}[Y_{j}].$$
Now we recall that for any K\"ahler metric $\omega$ on $X$, the cohomology class of the next current
$$\mu^\star (\omega)- \sum_{j\in J}\eta^{j}[Y_j]$$
contains a positive representative, for some family of positive rational numbers 
$(\eta^{j})$, which can be 
chosen as small as we want. Here we use the fact that all the $\mu$-exceptional divisors are among the $(Y_j)$ above ; as for the non-exceptional ones, they are absorbed by $\omega$.
 We remark that once such a family is fixed, any small enough perturbation of it
will have the same properties. 
Therefore we can assume that all the rational numbers
$${{1+a^{j}_{K_X+L}+ a^j_{\wh X/X}- a^{j}_L}\over {a^{j}_D- a^{j}_L- ma^{j}_{K_X+L}}}$$
are distinct, where the $a^j_L$ in the above quotient is the coefficient of 
$\nu^\star \Theta_L$ {\sl after perturbation}. In what will follow, we still use the same notations for 
the inverse image of $L$, but we keep in mind that $\wh \Lambda_L$ is now 
a K\"ahler metric on $\wh X$.

Our $L^2$ condition in the theorem 1.B.1 is translated via the blow up 
as follows
$$1+a^{j}_{K_X+L}+ a^j_{\wh X/X}- a^{j}_L> 0$$
for any index $j\in J$, therefore we have
the inequalities $(\bullet)_j$ above 
are automatically satisfied for the indexes $j$ such that 
$$a^{j}_D- a^{j}_L- ma^{j}_{K_X+L}\leq 0.$$
Therefore we obtain
$$\tau= {{1+a^{j}_{K_X+L}+ a^j_{\wh X/X}- a^{j}_L}\over {a^{j}_D- a^{j}_L- ma^{j}_{K_X+L}}}
$$
for a unique $j= j_0$ ; we remark that $\tau\in \bQ$, by our rationality hypothesis. We will denote by $\displaystyle S:= Y_{j_0}$.

\medskip
In conclusion, we have the numerical identity
$$
K_{\wh X}+ S+ \wh L\equiv \wh \Theta\leqno (15)$$
where $\wh L$ is any $\bQ$-bundle on $\wh X$
whose associated Chern class contains the current
$$\sum_{j\in J_p}(\tau a^{j}_D+ (1-\tau)a^{j}_{L}- 
(1+ m\tau )a^{j}_{K_X+L}- a^j_{\wh X/X})[Y_{j}]+ (1-\tau)\wh \Lambda_{L}.
\leqno (16)$$
We denote by
$$\wh \Theta:= \sum_{j\in J_n}((1+ m\tau )a^{j}_{K_X+L}+ a^j_{\wh X/X}- \tau a^{j}_D- (1-\tau)a^{j}_L)[Y_{j}]+ (1+ m\tau )\Lambda_{K_X+L}.
$$
The $J_p$ (respectively $J_n$) are the sets of indexes $j\in J$
 for which the corresponding
coefficient of $Y_j$ in the expression (16) and in the definition of 
$\wh \Theta$ is positive ; we remark that we have 
$$J= J_p\cup J_n\cup \{j_0\}.$$
\smallskip

\noindent It is at this point that one can see the utility of the 
relative critical exponent defined above :
 the $\bQ$-bundle $\wh L$ {\sl and its restriction to $S$} are big, and that the coefficients of their respective singular part 
are strictly smaller than 1. Also, the numerical identity $(15)$
restricted to $S$ show that $K_S+ \wh L_{|S}$ is pseudoeffective, and it has a metric with 
analytic singularities, since we have 
$$K_{S}+ \wh L_{|S}\equiv  \sum_{j\in J_n}((1+ m\tau )a^{j}_{K_X+L}+ a^j_{\wh X/X}- \tau a^{j}_D- (1-\tau)a^{j}_L)[Y_{j|S}]+ (1+ m\tau )\Lambda_{K_X+L|S}\leqno (17)$$

\noindent Therefore by induction we infer that 

$$H^0\big(S, p(K_S+ \wh L_{|S})\big)\neq 0\leqno (18)$$
for all large and divisible $p$.
\smallskip 
\noindent We claim next that any section 
$$u\in H^0\big(S, p(K_S+ \wh L_{|S})\big)$$
admits an extension $U$ to $\wh X$. This is a immediate application of the invariance of plurigenera of Siu, in the version due to
Hacon-McKernan [19]. Indeed, we have :
{\itemindent 4mm

\item {1.}  The bundle $\wh L$ is decomposed as sum of an ample line bundle 
and an effective one with critical exponent greater than 1~;

\item {2.} The bundle  $K_{\wh X}+ S+ \wh L$ admit a metric with positive curvature,
whose singular part is transversal to the support of the effective part of $\wh L$.

}

\noindent Thus by the extension theorem in [19] we obtain
$U\in H^0\big(\wh X, p(K_{\wh X}+ S+ \wh L)\big)$ such that $U_{|S}= u$.
\smallskip
A last observation is that since $0< \tau< 1$ we have 

$$ \mu^\star \big(\tau D+ (1-\tau)L\big)= S+ \wh L+ F$$
where some multiple of $F$ has non-zero sections. Indeed, by the formulas $(16)$ and (8)-(11) we have

$$F= \sum_{j\in J_p\cup \{j_0\}}\big((1+ m\tau )a^{j}_{K_X+L}+ a^j_{\wh X/X}\big)[Y_{j}]+
\sum_{j\in J_n}\big(\tau a^{j}_D+ (1-\tau)a^{j}_{L}\big)[Y_j].$$

\medskip
\noindent Now we twist an appropriate power of $U$ with the section of the corresponding multiple of $F$ and we obtain a section in a multiple of the bundle
$$K_X+ \tau \big(L+ m(K_X+ L)\big)+ (1-\tau)L$$
which is nothing but $(1+ \tau m)(K_X+ L)$.

\noindent Therefore, the particular case of the theorem 0.1 is
completely proved. \hfill\qed

\medskip
\claim Remark|{\rm The only difference between the theorem 1.B.1 and the result of Shokurov in [34] is the existence of the singular part $(Y_j)$ in the expression of $\wh L$ and $\wh \Theta$ above. However, the singular part encoded by $\wh L$
is small (the coefficients are smaller than 1) and the one inside $\wh \Theta$
is transversal to the singularities of $\wh L$ : it is {\sl precisely} for this reason that the invariance of plurigenera [17] still holds, and it is used to replace the Kawamata-Viehweg vanishing theorem (cf. [22], [45], [30]) to give the desired result. \hfill\qed

}
\endclaim

\bigskip

\subsection {\S 1.C Relative threshold of the logarithmic part of $\Theta_{K_X+ L}$}

\smallskip
\noindent We continue our proof of the general case of 0.1 by introducing a version of the relative threshold which is adapted to the case where
$\displaystyle \Theta_{K_X+L}$ does not necessarily have log poles. The main motivation
is that we still want to use the same approach as in the previous section and settle the general case by induction.

We recall that we have $\nd(\{K_X+L\})\geq 1$. Let $x_0\in X$ be a very general point,
such that $\displaystyle \nu(\Theta_{K_X+L}, x_0)= 0$.
By the theorem 1.A.3, there exist a closed positive current 
$$T\in m(c_1(K_X)+ \theta_L)+ \theta_L$$
with logarithmic poles, such that $\nu (T, x_0)> n+1$ (the positive integer $m$
above is large and fixed during the rest of the proof). We remark that $T$ is the substitute for the $\bQ$-section $D$ in the previous paragraph.

Let $\mu_0: \wt X\to X$ be a common log resolution of the currents $T$ and 
$\Theta_L$. By this we mean that $\mu_0$ is the composition of a sequence of 
blow-up maps with non-singular centers, such that we have

$$\mu_0^\star(\Theta_L)= \sum_{j\in J}a^{j}_L[Y_{j}]+ \wt \Lambda_L\leqno (19)$$

$$\mu_0^\star(T)= \sum_{j\in J}a^{j}_T[Y_{j}]+ \wt \Lambda_T
\leqno (20)$$

$$K_{\wt X/X}=  \sum_{j\in J}a^{j}_{\wt X/X}[Y_j]
\leqno (21)$$
where the divisors above are assumed to be non-singular and to have normal crossings.

We remark that the existence of the manifold $\wt X$ respectively of the map 
$\mu_0$ is a consequence of the fact that the currents $T$ and $\Theta_L$ 
have log poles. 

\noindent Now the current $\displaystyle \Theta_{K_X+ L}$ enter into the picture.
Let us consider its inverse image {\sl via} the map $\mu_0$ :

$$\mu_0^\star\big(\Theta_{K_X+ L}\big)= \sum_{j\in J}a^{j}_{K_X+ L}[Y_{j}]+ 
 \wt \Lambda_{K_X+ L}~;\leqno (22)$$
where $\wt \Lambda_{K_X+ L}$ in the relation above is a closed positive current, such that 
$$\nu_{Y_j}(\wt \Lambda_{K_X+ L})= 0\leqno (23)$$
for all $j\in J$ (the decomposition (22) is a direct consequence of a result due
to Siu in [36]). 

In other words, even if $\wt \Lambda_{K_X+ L}$ is not smooth anymore (as it was the case in the previous paragraph), {\sl its generic Lelong number along all the possible candidates for the hypersurface $S$ in the previous section is zero}. 
\medskip
\noindent Let us denote by $\wt D$ the first (divisor-type) part 
of the current $\displaystyle \mu^\star\big(\Theta_{K_X+ L}\big)$~; we consider the 
next quantity  
$$\tau:= \sup\{t\in \bR_+:  \int_{\wt X}\exp\big(t(\varphi_L\circ \mu_0+ 
{m\varphi_{\wt D}}- \varphi_{T}\circ \mu_0)+ \varphi_{\wt D}+ \varphi_{\wt X/X}- \varphi_L\circ \mu_0\big)d\lambda<\infty \}\leqno (24)$$

\noindent We will prove next that we have the relations
$$0<\tau< 1.$$
The  latter inequality is due to the fact that the singularity of $T$ at $x_0$ is large, and moreover 
$\mu_0^{-1} (x_0)$ is disjoint from the support of $\wt D$ (thus, it is here that the 
choice of a very generic point $x_0$ is important). As for the former inequality, 
we remark that by hypothesis we have
$$\int_{\wt X}\exp\big((1+ \varepsilon_0)(\varphi_{K_X+ L}\circ \mu_0- \varphi_L\circ \mu_0)+ \varphi_{\wt X/X}\big)d\lambda< \infty$$
so in particular
$$I:= \int_{\wt X}\exp\big((1+ \varepsilon_0)(\varphi_{\wt{D}}+ \varphi_{\wt \Lambda}- \varphi_L\circ \mu_0+ \varphi_{\wt X/X})\big)d\lambda< \infty$$
We consider a point $y\in \wt X$ such that 
$$\displaystyle\nu (\widetilde{\Lambda}_{K_{{X}}+ L}, y)= 0$$ 
and we have
$$
\int_{(\wt X, y)}\exp\big(\varphi_{\wt D}- \varphi_L\circ \mu_0+ \varphi_{\wt X/X}\big)d\lambda\leq  I^{1/{1+\varepsilon_0}} 
\Big(\int_{(\wt X, y)}\exp\big(-(1+1/\varepsilon_0)\varphi_{\wt \Lambda}\big)d\lambda
\Big)^{\varepsilon_0/1+\varepsilon_0} $$
and the last integral is convergent since the Lelong number of $\wt \Lambda_L$
at $y$ is equal to zero, therefore we can apply the theorem of Skoda, see [35].

In conclusion, the function $\exp\big(\varphi_{\wt D}- \varphi_L\circ \mu_0+ \varphi_{\wt X/X}\big)$ is in $L^1$ at the generic point of each $Y_j$, by the relation (18). Together with the fact that the hypersurfaces $(Y_j)$ have normal crossing, this prove that 
$\tau> 0$.\hfill\qed

\medskip

\noindent As before, we can obtain the explicit expression of the threshold 
$\tau$ (modulo perturbation) by using the modification $\mu_0$ as follows. Given any real number $t$, we have

$$ \eqalign {
\mu^\star \big(\Theta_\omega(K_X) + t(T - \Theta_{L})+ & \Theta_L\big)
\equiv  \hskip 1mm \Theta_{\wh \omega}(K_{\wt X})+ (1+ mt)\wt D + (1-t)\wt \Lambda_L+ t\wt \Lambda_T\cr
&+  \sum_{j\in J}(ta^{j}_T+ (1-t)a^{j}_L- (1+ mt)a^{j}_{K_X+L}- 
a^j_{\wt X/X})[Y_{j}].\cr
}
$$
We have $\mu^\star\big(\Theta_{K_X+ L}\big)= \wt D+ \wt \Lambda_{K_X+ L}$, 
and on the other hand the cohomology class of the current
$$t(T - \Theta_{L})+ \Theta_L- (1+ mt)\Theta_{K_X+ L}$$
is equal to the first Chern class of $X$, so
the previous relation can be written as

$$ \eqalign {
(1+ mt)\wt \Lambda_{K_X+ L} \equiv &
\Theta_{\wh \omega}(K_{\wt X})+ \sum_{j\in J}(ta^{j}_T+ (1-t)a^{j}_L- (1+ mt)a^{j}_{K_X+L}-a^j_{\wt X/X})[Y_{j}]+ \cr
& (1-t)\wt \Lambda_L+ t\wt \Lambda_T.\cr
}
$$
The perturbation argument used in the previous paragraph 
is still valid~; in conclusion, for $t:= \tau$ (or better say, a slight modification of this quantity, since we change the inverse image of the current $\Theta_L$ {\sl within the same 
cohomology class}) we can assume that we have

$$(1+ m\tau)\wt \Lambda_{K_X+ L}+ \wt G\equiv  c_1(K_{\wt X}+ \wt S)+ \theta_{\wt L} \leqno (25)$$
where the notations we use are as follows :
$$\leqno (26) \eqalign {
\theta_{\wt L}:\equiv 
& \sum_{j\in J_p} (\tau a^{j}_T+ (1-\tau)a^{j}_L- (1+ m\tau)a^{j}_{K_X+L}
-a^j_{\wt X/X})[Y_{j}] + \cr
+ & (1-\tau)\wt \Lambda_L+ \tau\wt \Lambda_T\cr
 }
 $$
 as well as
 $$\wt G:= \sum_{j\in J_n} \big((1+ m\tau)a^{j}_{K_X+L}+ 
 a^j_{\wt X/X}- \tau a^{j}_T- (1-\tau)a^{j}_L\big)[Y_{j}] ;
\leqno (27)$$
after perturbation, we can assume that $\wt \Lambda_L$ is a K\"ahler metric.
In the above relations, we have used the same conventions as in 1.B for the definition of $J_p, J_n$.

\medskip
\subsection {\S 1.D Properties of $\theta_{\wt L}$}
\smallskip

\noindent We collect in this paragraph the main features of the class 
$\theta_{\wt L}$ which will be needed later.
\smallskip

\noindent {$\bullet$} In the first place, by the definition of the threshold $\tau$, we see that the coefficients of the singular part in the expression (26) are strictly smaller than 1. Thus, $\theta_{\wt L}$ contains a closed positive current whose critical exponent is strictly greater than 1~; it is equally a K\"ahler current, as it dominates $(1-\tau)\wt \Lambda_L$, and this form is positive definite on $\wt X$. 

\smallskip

\noindent {$\bullet$} There exist 
an effective $\bR$--divisor $\Delta$ such that 
$$\eqalign{
\theta_{\wt L}+ \{\wt S+ \Delta\}= & \mu_0^\star \big(\tau \{T\}+ (1-\tau)\theta_L\big) \equiv \cr
\equiv & \mu_0^\star\big(\theta_L+ \tau m(c_1(K_X)+\theta_L)\big)\cr
}
$$
(the precise expression of $\Delta$ is not relevant for the moment, but one can easily get it from the 
relations (25), (26) and (27) above).

\noindent Therefore, it is enough to produce an effective $\bR$--divisor 
in the cohomology class $c_1(K_{\wt X}+ \wt S)+ \theta_{\wt L}$ in order to complete the proof of the theorem 0.1.
\smallskip

\noindent {$\bullet$} The adjoint class $c_1(K_{\wt X}+ \wt S)+ \theta_{\wt L}$ is pseudoeffective~; moreover, it contains the closed positive current 
$$(1+ m\tau)\wt \Lambda_{K_X+ L}+ \wt G$$
whose Lelong number at the generic point of $\wt S$ is equal to zero. 

\smallskip

\noindent {$\bullet$} By using a sequence of blow-up maps, we can even assume that the components $\displaystyle (Y_j)_{j\in J_p}$
of the representative (26) of 
$\theta_{\wt L}$ have empty mutual intersections. Indeed, this is a simple --but nevertheless crucial!-- result due to Hacon-McKernan, which we recall next. 

We denote by 
$B$ an effective $\bQ$-divisor, whose support do not contain $\wt S$, such that 
$\Supp B\cup \wt S$ has normal crossings and such that its coefficients are strictly smaller than 1.

\claim Lemma ([20])|There exist a birational map $\mu_1: \wh X\to \wt X$ such that 
$$\mu_1^\star (K_{\wt X}+ \wt S+ B)+ E_{\wh X}= K_{\wh X}+ S+ \Gamma$$
where $E_{\wh X}$ and $\Gamma$ are effective with no common components, 
$E_{\wh X}$ is 
exceptional and $S$ is the proper transform of $\wt S$~; moreover, the support of the divisor $\Gamma$ has normal crossings, 
its coefficients are strictly smaller than 1 and the intersection of any two components is empty.
\hfill\qed
\endclaim
The proof of the above lemma is by induction on the number of the components of 
$B$ having non-empty intersection~; the sets which we blow up to obtain $\mu_1$ are precisely the said intersections. Since from the start the components of $B$ have normal crossings, the restriction of the map $\mu_1$ to the proper transform $S$ of $\wt S$ will be isomorphic at the generic point of $S$. \hfill\qed
\smallskip
In our setting the $\bR$--divisor $B$ above is defined as
$$B:= \sum_{j\in J_p} (\tau a^{j}_T+ (1-\tau)a^{j}_L- (1+ m\tau)a^{j}_{K_X+L}
-a^j_{\wt X/X})Y_{j}$$
and we have
$$\mu_1^\star \big(c_1(K_{\wt X}+ \wt S)+ \big\{[B]+ (1-t)\wt \Lambda_L+ t\wt \Lambda_T\big\}\big)+ \{E_{\wh X}\}= c_1(K_{\wh X}+ S)+ \theta$$
where the cohomology class $\theta$ above contain a representative
which can be written as follows 
$$\sum_{j\in I}\rho^{j}W_j+ \wh \Lambda$$
where $0< \rho^{j}< 1$ for any $j$, the hypersurfaces $W_j$ are non-singular and 
$W_j\cap W_k= \emptyset$ if $j\neq k$, and $\wh \Lambda$ is smooth, semipositive, whose restriction to $S$ is positively defined at the generic point.

In addition, we remark that the new modification $\mu_1$ does not affect the bullets above. We summarize the discussion in this paragraph in the next statement (in which we equally adjust the notations).

\claim 1.D.1 Proposition|There exist a birational map $\mu: \wh X\to X$ and 
a class $\theta_{\wh L}\in \NS_\bR(\wh X)$ which contain the current 
$$\sum_{j\in J}\nu^{j}Y_j+ \wh\Lambda_L$$
where
$0< \nu^{j}< 1$, the hypersurfaces $Y_j$ above are smooth, they have empty mutual intersection and moreover the following hold :

{\itemindent 5mm

\item {i)} There exist a closed positive $(1, 1)$--current $\Theta$ on $\wh X$ with the property that 
$$\Theta\in c_1(K_{\wh X}+ S)+ \theta_{\wh L}$$
where $S\subset \wh X$ is a non-singular hypersurface which has transversal intersections with $(Y_j)$~;
\smallskip
\item {ii)} The support of the divisorial part of the current $\Theta$ is disjoint from the 
set $(S, Y_j)$~; 
\smallskip
\item {iii)} There exist a map 
$\mu_1: \wh X\to \wt X$ such that $S$ is not $\mu_1$--exceptional, and such that 
$\wh\Lambda_L$ is greater than the inverse image of a K\"ahler metric on $\wt X$
via $\mu_1$. Therefore, the form $\wh\Lambda_L$ is positive defined at the generic point of $\wh X$, and so is its restriction to the generic point of $S$~;
\smallskip
\item {iv)} There exist an effective $\bR$-divisor $\Delta$ on $\wh X$
such that 
$$\theta_{\wh L}+ \{S+ \Delta\}=  \mu^\star \Big(\theta_L+ \tau m\big(c_1(K_X)+\theta_L\big)\Big)+ \{ E\}$$
where $E$ is $\mu$--exceptional.
 
 }
\endclaim 
\noindent We remark that the relation $ii)$ imply that the Lelong number of $\Theta$ at the generic point of $S$ is equal to zero but nevertheless, the local potentials of  
$\Theta$ may be identically $-\infty$ when restricted to $S$. Therefore in order to 
be able to use the induction hypothesis, we have to regularize it. Certainly this creates some difficulties, which we overcome along the next paragraphs.


\medskip 

\subsection {\S 1.E Regularization and induction}

\smallskip
\noindent We consider now the family of approximations $(\Theta_\varepsilon)_{\varepsilon> 0}$
of the current $\Theta$ given by the theorem 1.A.1. For each $\varepsilon> 0$, the 
current $\Theta_\varepsilon$ has log poles, and we equally have
$$\nu(\Theta_\varepsilon, x)\leq \nu(\Theta, x)$$
for any $x\in \wt X$. In particular, we have $\nu(\Theta_\varepsilon, x)= 0$
for the generic point $x\in S$~; in other words, the restriction of the current $\Theta_\varepsilon$ to $S$ is well defined (i.e. its potential is not identically equal to $-\infty$
at each point of $S$).
Therefore, since the regularization process {\sl does not} change the 
cohomology class, the relation $i)$ of the proposition 1.D.1 above imply 
$$\Theta_{\varepsilon|S}\in c_1(K_{S})+ \theta_{\wh L{|S}}\leqno (28)$$
Next, we have the following decomposition
$$\Theta_{\varepsilon|S}= \sum_{j\in J}\rho^{\varepsilon, j}Y_{j|S}+ R_{\varepsilon}\leqno (29)$$
where the coefficients $(\rho^{\varepsilon, j})$ are positive real numbers, 
the generic Lelong number of $R_{\varepsilon}$ along $Y_j\cap S$ is zero,
and moreover we have 
$$R_\varepsilon\geq -\varepsilon \omega_{|S}$$
We remark that this current 
may be singular along some other {\sl hypersurfaces} of $S$. 

For each index $j\in J$ we will assume that the next limit
$$\rho^{\infty, j}:= 
\lim_{\varepsilon\to 0}\rho^{\varepsilon, j}$$
exist, and we introduce the following notation
$$I:= \{j\in J : \rho^{\infty, j}\geq \nu^{j}\}.\leqno (30)$$
The numerical identity (28) coupled with (29) show that we have
$$\sum_{j\in I}(\rho^{\infty, j}- \nu^{j})[Y_{j|S}]+
R_{\varepsilon}+ \sum_{j\in J}(\rho^{\varepsilon, j}-\rho^{\infty, j})[Y_{j|S}]
\in c_1(K_S)+ \theta_{L_S}\leqno (31)$$
where $\displaystyle \theta_{L_S}$ is the cohomology class of the current 
$$\sum_{j\in J\setminus I}(\nu^{j}-\rho^{\infty, j})[Y_{j|S}]+\wh \Lambda_{L|S}.$$
We infer that $\displaystyle \theta_{L_S}$ contains a K\"ahler current, by the property 1.D.1, $iii)$ and the fact that the above coefficients are positive. Its critical exponent is greater than 1, by the definition of the set $I$ in (30), the proposition 1.D.1 and the fact that the 
restrictions $Y_{j|S}$ are 
smooth and they do not intersect each other. Moreover, the adjoint class
$\displaystyle c_1(K_S)+ \theta_{L_S}$ is pseudoeffective, by the compactness property of closed positive currents.
Indeed, any weak limit of the currents $R_\varepsilon$ above will be positive, and the last term in the numerical identity (31) will converge to zero as $\varepsilon\to 0$.
\smallskip
Now we apply the induction hypothesis : the class 
$\displaystyle c_1(K_S)+ \theta_{L_S}$ contains a non-zero, effective $\bR$-divisor , which can be written as
$$T_S:= \sum_{i\in K}\lambda^{i}[W_i]$$
where $W_i\subset S$ are hypersurfaces. We consider now the current
$$\wh T_S:= \sum_{i\in K}\lambda^{i}[W_i]+ 
\sum_{j\in J\setminus I}\rho^{\infty, j}[Y_{j|S}]+  
\sum_{j\in I}\nu^{j}[Y_{j|S}]~;\leqno (32)$$
from the relation (31) we get 
$$\wh T_S\in c_1(K_{\wh X}+ S)+ \theta_{\wh L{|S}}.\leqno (33)$$

During our discussion of the case where $\displaystyle \Theta_{K_X+ L}$
has log poles and rational singularities, we have used the extension theorem [19]
in order to lift the section whose zero set is (in the actual context) $\wh T_S$. 
The obstructions in order to do the same thing now are the following :

{\itemindent 6mm

\item {$(\bullet)$} $\theta_{\wh L}$ does not necessarily correspond to a 
$\bQ$-bundle. In fact, we remark that 
this could happen even if at the very beginning it is the case for $\theta_L$ --the reason is
that the Lelong numbers of $\displaystyle \Theta_{K_X+ L}$ 
are not rational numbers in general~;
\smallskip
\item {$(\bullet)$} The current $\wh T_S$ do not (necessarily) correspond 
to the zero set of a $\bQ$-section~;
\smallskip
\item {$(\bullet)$} The current $\Theta_\varepsilon$ is greater than $-\varepsilon\omega$, but not 
positive in general.

}
\medskip
\noindent The way to deal with the first difficulty is by diophantine approximation, and will be explained in the next subsection. Then the second one is not a serious problem, since the method of Shokurov can be adapted 
to this setting. Finally, the slight negativity of the current $\Theta_\varepsilon$ has the following consequences : when we apply the ``invariance of plurigenera" iteration, we can perform {\sl 
only a finite number of steps}. However, we will show that as soon as the approximation is accurate enough,
the number of steps is sufficiently large to allow us to conclude (here we use in an essential manner the
strong positivity of $\wh \Lambda_L$).


\bigskip

\subsection {\S 1.F Approximation}

\medskip
\noindent In this paragraph we would like to prove the following diophantine approximation lemma (see [5], as well as [39] for similar considerations).

\claim  1.F.1 Lemma|For each $\eta> 0$, there exist a positive integer $q_\eta$, a $\bQ$--line bundle 
$\wh L _\eta$ on $\wh X$ and a closed positive current
$$\wh T_{S,\eta}:= \sum_{i\in K}\lambda^{i}_\eta[W_i]+ \sum_{j\in J\setminus I}
\rho^{\infty, j}_\eta[Y_{j|S}]+  
\sum_{j\in I}\nu^{j}_{\eta}[Y_{j|S}]~;\leqno (34)$$
such that :

\item {A.1} The multiple $q_\eta\wh L _\eta$ is a genuine line bundle, and the numbers
$$(q_\eta\lambda^{i}_\eta)_{i\in K}, (q_\eta\nu^{j}_\eta)_{j\in J}, 
(q_\eta\rho^{\infty, j}_\eta)_{j\in J}$$ 
are integers~;
\smallskip
\item {A.2}  We have $\displaystyle \wh T_{S,\eta}\in \{K_{\wh X}+ S+ \wh L _{\eta |S}\}$ ;
\smallskip
\item {A.3} We have $\Vert q_\eta\big(\wh L- \wh L _\eta\big)\Vert< \eta$,
$|q_\eta\big(\lambda^{i}_\eta- \lambda^{i}\big)| <\eta$ and the analog relation 
for the $(\rho^{\infty, j}, \nu^j)_{j\in J}$ 
(here $\Vert \cdot\Vert$ denotes any norm on 
the real Neron-Severi space of $\wh X$)~;
\smallskip
\item {A.4} For each $\eta_0> 0$, there exist a finite family $(\eta_j)$ such that 
$\{K_{\wh X}+ S\}+ \theta_{\wh L}$ belong to the convex hull of  
$\displaystyle \{K_{\wh X}+ S+ \wh L _{\eta_j}\}$
where $0< \eta_j< \eta_0$.\hfill\qed

\endclaim

\medskip
\proof \hskip 0.7mm (of the approximation lemma). 
We choose first an appropriate basis of the 
Neron-Severi space of $S$. 

Let $\alpha_1,..., \alpha_\rho$ be a basis of the $\NS(X)$.
The restrictions $(\alpha_{j|S})$ generate a subspace of $\NS(S)$~; let us assume that
a basis of this subspace is given by $\displaystyle (\alpha_{j|S})_{1\leq j\leq r}$. We complete this free family with $\beta_1,..., \beta_{\rho^\prime}$ to a basis of 
$\NS(S)$.

We compute the coordinates of the cohomology class of the current 
$$\wh T_S= \sum_{i\in K}\lambda^{i}[W_i]+ \sum_{j\in J\setminus I}
\rho^{\infty, j}[Y_{j|S}]+  
\sum_{j\in I}\nu^{j}[Y_{j|S}] 
\in \{K_{\wh X}+ S\}+ \theta_{\wh L|S}\leqno (35)$$
with respect to the basis 
$(\alpha, \beta)$ above. Since $\{K_{\wh X}+ S\}+ \theta_{\wh L}$ is defined on $\wh X$, there exist a set of real numbers $(b^{j})$ such that
$$\{K_{\wh X}+ S\}+ \theta_{\wh L}= \sum _{j= 1}^\rho b^{j}\alpha_j.$$
We equally have the {\sl rational numbers }
$(q_p^j)$ such that 
$$\alpha_{r+j|S}= \sum _{k=1}^rq^k_j\alpha_{k|S}$$
by our assumption, and therefore we obtain
$$\{K_{\wh X}+ S\}+ \theta_{\wh L|S}= \sum _{p= 1}^r \big(b^{p}+ \sum_{j=1}^{\rho-r}b^{r+j}q^p_j\big)\alpha_{p|S}.\leqno (36)$$
In order to simplify the writing, 
we denote by $a^p:= b^{p}+ \sum_{j=1}^{\rho-r}b^{r+j}q^p_j$.

Next, we are going to express the coordinates of
$c_1(W_j)$ and $c_1(Y_j)$ with respect to the basis $(\alpha, \beta)$.
There exist the rational numbers $(x^p_j), (y^q_j), (z^k_j)$
such that 

$$W_j\equiv \sum_{p=1}^{r}x^p_j\alpha_{p|S}+ 
\sum_{p=1}^{\rho^\prime}y^p_j\beta_p.$$
as well as
$$\eqalign {
Y_{j|S}\equiv & \sum_{k=1}^{\rho}z^k_j\alpha_{k|S} \cr
\equiv &  \sum_{k=1}^{r}\big(z^k_j+ \sum_{p=1}^{\rho- r}z^{p+r}_jq^k_p\big)\alpha_{k|S}
\cr
}
$$

By the relations (32), (33) and (36) we get the next equalities
$$\sum_{j\in K}\lambda^{j}y_j^p= 0\leqno (37_p)$$
$$\sum_{i\in K}\lambda^{i}x_i^k
+\sum_{j\in J\setminus I}\rho^{\infty, j}\big(z^k_j+ \sum_{p= 1}^{\rho- r}z^{p+r}_jq^k_p\big)
+ \sum_{j\in I}\nu^{j}\big(z^k_j+ \sum_{p=1}^{\rho- r}z^{p+r}_jq^k_p\big)= a^{k}
\leqno (38_k)$$
for $p= 1,..., \rho^\prime$ and $k= 1,...,r$.

\noindent Let $Z$ be the matrix whose coefficients are the rational numbers 
$(x_j^p), (y_j^q)$ as well as  $(z^k_j+ \sum_{p=1}^{\rho- r}z^{p+r}_jq^k_p)$ such that the above equations can be written as 
$$Z\Lambda= V\leqno (39)$$
where $V$ is the vector whose first $\rho^\prime$ coefficients are
zero, and the next ones are just the $(a^{k})$. Conversely, if we have a solution 
$(\lambda, \rho, \nu)$ of the linear system (39) corresponding to a vector $V$
whose first  $\rho^\prime$  components are zero, then the cohomology class of the corresponding $\wh T_S$ belong to the space generated by $(\alpha_j)$.

We denote by $H\subset \NS(S)$ the vector space
obtained by intersecting the image of $Z$ with the space defined by the vanishing 
of the first $\rho^\prime$ coordinates~; the equality (39) show that the vector $V$ belong to $H$. 
\smallskip

Next, the main fact is that the subspace $H$ is defined over the rational numbers i.e. there exist a basis $h_1,... h_l$ of $H$ which can be expressed in $(\alpha, \beta)$
with rational coefficients (this is a consequence of the fact that the matrix $Z$ has rational coefficients). Then we write our vector $V$ with respect to this basis of $H$,
as follows 
$$V= \sum_{p=1}^lc^ph_p.$$
We will use the next fact from the ``diophantine approximation" 
theory.

\claim 1.F.2 Fact|Let $x^1, ..., x^s$ be a finite set of real numbers~;
there exist a constant $C> 0$ such that for any $\eta> 0$, there exist 
$q_\eta\in \bZ_+$ such that :
\smallskip

{\itemindent 5mm
\item {(1)} For each $j= 1,..., s$, there exist integers $p^j$ so that we have 
$|q_\eta x^j- p^j|< C\eta$ we denote by $X$ the vector in $\bR^s$ whose coefficients are the $(x^j)$, and by $X_\eta$ the vector given by the rational approximations 
$\displaystyle {{p^j}\over {q_\eta}}$~;
\smallskip
\item {(2)} Given a positive real $\eta_0$, there exist a family of approximations 
$\displaystyle X_{\eta_j}$ as above, with $\eta_j< \eta_0$ such that the vector $X$ belongs to the convex hull of $\displaystyle X_{\eta_j}$.

}

\endclaim

\proof. Indeed, we show along the next lines that the above statement is a consequence of the Kronecker theorem (see e.g. Hardy-Wright, [21]).

We consider the maximal family, 
say $x^1,..., x^p$ such that 
the numbers $1, x^1,..., x^p$ are independent over $\bQ$. The integer $p$ 
above can be assumed to be greater than 1, because if not all $(x^j)$ are rational, and in this case the above result is obvious. 
Thus, we have
$$x^{p+j}= \sum_{k= 1}^p r^j_kx^k+ r^j_0$$
for all $j= 1,..., s-p$,  
where $(r^j_k)$ is a set of rational numbers. Let 
$A= (r^j_k)$ and $B= (r^j_0)$~; then our vector $(x^1,..., x^s)$ become an 
a point on the graph of the affine map
$$Y\to AY+ B$$
defined on $\bR^p\to \bR^{s-p}$. 

Let $\eta>0$ be a positive rational number. The theorem of Kronecker 
(see [21]) imply that the set 
$\displaystyle \big(\{\varepsilon_1qx^1\},..., \{\varepsilon_p qx^p\}\big)_{q\in \bZ_+}$ is dense in the cube $[0,1]^p$, for any choice of
the quantities $\varepsilon_j$ within the set $\{-1, 1\}$. Here we use the fact that if the family $1, x^1,..., x^p$ is independent over $\bQ$, then so is 
$1, \varepsilon_1x^1,..., \varepsilon_px^p$, for any choice of $(\varepsilon_j)$
as above.

We consider the vectors
$$X_{\varepsilon, q}:= \Big(x^1-{{\varepsilon_1}\over {q}}\{\varepsilon_1qx^1\},...,
x^p-{{\varepsilon_p}\over {q}}\{\varepsilon_pqx^p\}\Big)~;$$
our first remark is that they have rational coefficients, as one can easily see. 
Next, for any $(\varepsilon_j)$ as above, we will consider the vectors 
$X_{\varepsilon, q}$ such that $\{\varepsilon_jqx^j\}< \eta$, for any $j=1,..., p$.
Finally, we can select a finite family from the 
($X_{\varepsilon, q}$) above such that the vector $(x^1,..., x^p)$ is in the interior of the polygon they define. 

But then, the vector $(X, AX+B)$ is in the convex hull of the 
$$(X_{\varepsilon, q}, AX_{\varepsilon, q}+ B)$$
and this completes the proof of 1.F.1.
 \hfill\qed 

\bigskip 
We apply now the diophantine approximation statement above for the following set 
of real numbers : 
$$(c^j)_{1\leq j\leq l}, (b^{r+j})_{1\leq j\leq \rho-r}, (\rho^{\infty, j})_{j\in J}, 
(\nu^j)_{j\in J}, (\lambda^j)_{j \in K}.
$$

To complete the proof of the lemma, will trace our steps back in order to obtain the approximations 
of the $\bR$--sections/bundles we are interested in. 

We first define 
$$V_\eta:= \sum_jc^j_\eta h_j\in H$$
and its coefficients with respect to the basis $(\alpha, \beta)$
will define the rational approximations $(a^j_\eta)$ 
of the coefficients $(a^j)$ above~;
moreover, the coefficients corresponding to $(\beta_l)$ are zero.
Then we set
$$b^j_\eta:= a^j_\eta- \sum_{k=1}^{\rho-r}b^{r+k}_\eta q^j_k~; $$
and
$$F_\eta:\equiv \sum _{j= 1}^\rho b^{j}_\eta \alpha_j,$$
 that is to say we consider any $\bQ$--line bundle whose Chern class
 is given by the above expression.
 Finally, we can choose the vector $\Lambda_\eta$ such that
 $$Z\Lambda_\eta= V_\eta\leqno (40)$$
 with the property that $$\Vert q_\eta (\Lambda_\eta- \Lambda)\Vert< C\eta$$
 and such that $m_0q_\eta\Lambda_\eta$ has all the components integral numbers,
 where $m_0$ is a fixed positive integer. Such a vector can indeed be found, 
 by the usual theory of linear systems~;
the constant $m_0$ is due to the inverse of some square matrix extracted from $Z$.
It is at this point that we need to have at hand the approximations 
of some of the $(\lambda^j), (\rho^{\infty, j}), (\nu^j)$~; the others will be imposed by the linear system (40).

We recall that all the components of the vector $\Lambda$ are positive real numbers, therefore the approximation $\Lambda_\eta$ will have the same property, 
if $\eta\ll 1$.

 The approximation $\wh T_{S,\eta}$ of the the current $\wh T_S$ in our lemma 1.F.1 is obtained simply by plugging in 
 the coefficients of $\Lambda_\eta$.
 
 The equation (40) show that 
 $$\{\wh T_{S,\eta}\}= c_1\big(F_{\eta|S}\big)~;$$
moreover, the relation (1) in 1.F.2 
imply
$$\Vert q_\eta\big(\{K_{\wh X}+ S\}+ \theta_{\wh L}-\{F_\eta\}\big)\Vert< C\eta ;$$
moreover, we can assume that $q_\eta F_\eta$ is integral (here we are a bit sloppy, since some additional 
denominators can occur because of the coefficients $(x, y, q)$ above, but they are fixed independently of $\eta$, so we just ignore them). 

Finally, we remark that $\{K_{\wh X}+ S\}+ \theta_{\wh L}$ belong to the convex hull
of the $F_\eta$ thanks to the second part of 1.F.2. This is indeed the case, since
the vector $(b_1,..., b_\rho)$ is the image of $(c_1,..., c_l, b_{r+1},..., b_\rho)$ via a linear map, and we use precisely the same map (whose associated matrix with respect to the basis above is has rational coefficients) to define the approximations.

We define 
$\wh L_\eta$ such that 
$$K_{\wh X}+ S+ \wh L_\eta\equiv  F_\eta$$
and the proof of the lemma is finished.\hfill\qed


\bigskip
\subsection {\S 1.G The method of Shokurov}

\medskip 
Our concern in this paragraph will be to ``convert" the current $\wh T_{S,\eta}$ into
a genuine section $s_\eta$ of the bundle 
$\displaystyle q_\eta\big(K_{\wh X}+ S+ \wh L _\eta\big)$. To this end, we will apply a
classical trick of Shokurov, in the version revisited by Siu in his recent work [39].
A crucial point is that by a careful choice of the metrics we use, the $L^2$ estimates will allow us to have a very precise information concerning the vanishing of $s_\eta$.

\claim 1.G.1 Proposition|There exist a section 
$$s_\eta\in H^0\Big(S, q_\eta \big(K_{S}+ \wh L _{\eta|S}\big)\Big)$$
whose zero set contains the divisor 
$$q_\eta\Big(\sum_{j\in J\setminus I}\rho^{\infty, j}_\eta[Y_{j|S}]+  
\sum_{j\in I}\nu^{j}_{\eta}[Y_{j|S}]\Big)$$ 
for all $0<\eta\ll 1$.

\endclaim

\medskip
\proof \hskip 0.5mm (of the proposition). We express the bundle we are interested in as an
adjoint bundle as follows
$$q_\eta \big(K_{S}+ \wh L _{\eta|S}\big)= K_S+ (q_\eta-1) \big(K_{S}+ 
\wh L _{\eta|S}\big)+ \wh L_{\eta|S}.$$

\noindent In order to use the classical vanishing theorems, we have to endow the bundle
$$(q_\eta-1) \big(K_{S}+ 
\wh L _{\eta|S}\big)+ \wh L_{\eta|S}$$
with an appropriate metric.
We first consider first the $\bQ$--bundle $\wh L _{\eta}$~; we will construct a metric on
it from the decomposition 
$$c_1(\wh L _\eta)= \theta_{\wh L}+ \big(c_1(\wh L _\eta)-\theta_{\wh L}).$$
The second term above admits a smooth representative whose local weights are 
bounded by $\displaystyle {{\eta}\over {q_\eta}}$ in $\cC^\infty$ norm, by the approximation
relation $A.3$. As for the first one, we recall that the class $\theta_{\wh L}$
contains the representative
$$\sum_{j\in J}\nu^{j}Y_j+ \wh \Lambda_L~; \leqno (41)$$
where the (1,1)-form $\wh \Lambda_L$ has the positivity properties in $1.D.1, iii)$.
\smallskip 
\noindent Now, the first metric we consider on $\wh L_{\eta|S}$ is given by the 
perturbation of the current (41) as follows :
$$\sum_{j\in I}\max\big(\nu^{j}, \nu^{j}_{\eta}\big)Y_{j|S}+ 
\sum_{j\in J\setminus I }\nu^{j}Y_{j|S}+ \wh \Lambda_{L|S}+ \Theta(\eta)_{|S}\leqno (42)$$
where $\Theta(\eta)$ is a non-singular $(1,1)$--form on $\wh X$ in the class of
the current 
$$\sum_{j\in I}\Big(\nu^{(j)}-\max\big(\nu^{(j)}, \nu^{(j)}_{\eta}\big)\Big)[Y_j]$$
plus $c_1(\wh L_\eta)-\theta_{\wh L}$ ; 
we can assume that it is greater than $\displaystyle -C{{\eta}\over {q_\eta}}$.

The smooth term $\wh \Lambda_{L|S}$ is semi-positive on $S$ and strictly positive at the generic point~; in order to gain the strict positivity needed in the vanishing theorems, we recall the following standard result (see e.g. [15], [28]).

\claim Fact|There exist 
a K\"ahler metric $\Omega$ on $S$, a positive constant $C$ and 
a family of currents $\wh \Lambda_{L,S}^\varepsilon\in \{\wh \Lambda_{L|S}\}$  
such that~:
\smallskip
{\itemindent 5mm
\item {(1)} We have $\wh \Lambda_{L,S}^\varepsilon\geq \varepsilon \Omega$~;
\smallskip
\item {(2)} The Lelong number of $\wh \Lambda_{L,S}^\varepsilon$ at each point of the manifold is smaller than $\varepsilon C$.\hfill\qed

}
\endclaim

\smallskip
\noindent Thus, for any $\eta> 0$, we can define an $\varepsilon$ such that the corresponding current $\wh \Lambda_{L,S}^\varepsilon$ verify the next properties.

\item {$M_1)$} The current $\wh \Lambda_{L,S}^\varepsilon+ \Theta(\eta)$ 
dominates a small multiple of $\Omega$, which depends on $\eta$, but fortunately
this does not matter for the purposes of this paragraph~;

\item {$M_2)$} The Lelong number of $\wh \Lambda_{L,S}^\varepsilon$ at each point of $S$ is smaller than $\displaystyle C{{\eta}\over {q_\eta}}$.

\medskip

\noindent In conclusion, we can define a metric on $\wh L_{\eta|S}$ with the following curvature current 
$$\sum_{j\in I}
\max\big(\nu^{j}, \nu^{j}_{\eta}\big)Y_{j|S}+ 
\sum_{j\in J\setminus I}\nu^{j}Y_{j|S}+ 
\wh \Lambda_{L,S}^\varepsilon + \Theta(\eta)_{|S} ;
\leqno (43)$$
we remark that it is a K\"ahler current, and its critical exponent is 
strictly greater than 1, provided that $\eta\ll 1$.
\smallskip
Next, we define a singular metric 
on the bundle $(q_\eta-1) \big(K_{S}+ \wh L _{\eta|S}\big)$
whose curvature form is equal to $(q_\eta-1)\wh T_{S,\eta}$ and
we denote by $h_\eta$ the resulting metric on the bundle
$$(q_\eta-1) \big(K_{S}+ \wh L _{\eta|S}\big)+ L_{\eta|S}.$$
\smallskip
The current $q_\eta\wh T_{S,\eta}$ corresponds to the current of integration along the zero set of the section $u_\eta$ of the bundle 
$$q_\eta \big(K_{S}+ \wh L _{\eta|S}\big)+ \rho$$
where $\rho$ is a topologically trivial line bundle on $S$. 

By the Kawamata-Viehweg-Nadel vanishing theorem we have
$$H^j(S, q_\eta \big(K_{S}+ \wh L _\eta\big)\otimes \cI\big(h_\eta\big)\Big)= 0$$
for all $j\geq 1$, 
and the same is true for the bundle $q_\eta \big(K_{S}+ \wh L _\eta\big)+ \rho$, since
$\rho $ carries a metric with zero curvature.
Moreover, the section $u_\eta$ belong to the multiplier ideal of the metric 
$h_\eta$ above, as soon as $\eta$ is small enough, because the multiplier ideal of the metric
(43) on the bundle $\wh L_{\eta|S}$ will be trivial. Since 
the Euler characteristic of the two bundles is the same, we infer that
$$H^0\Big(S, q_\eta \big(K_{S}+ \wh L _\eta\big)\otimes \cI\big(h_\eta\big)\Big)\neq 0$$
We denote by $s_\eta$ any non-zero element in the group above~; we show now that its
zero set satisfy the requirements in the lemma. Indeed, locally at any point of $x\in S$
we have
$$\int_{(S, x)}{{|f_s|^2}\over {\prod _{j\in J\setminus I}|f_j|^{2\rho^{\infty, j}_\eta(q_\eta-1)+ 2\wt\nu^{(j)}_{\eta}}
\prod_{j\in I}|f_j|^{2\nu^{j}_{\eta}(q_\eta-1)+ 2\wt\nu^{j}_{\eta}}}}
d\lambda<\infty $$
where $\wt\nu^{j}_{\eta}:= \nu^{j}$ if $j\in J\setminus I$ and $\wt\nu^{j}_{\eta}:= 
\max\{\nu^{j}_{\eta}, \nu^{j}\}$ if $j\in I$ ; we denote by $f_s$ the local expression of the section $s$, and we denote by $f_j$ the local equation of $Y_j\cap S$.

But the we have 
$$\int_{(S, x)}{{|f_s|^2}\over {\prod _{j\in J\setminus I}|f_j|^{2\rho^{\infty, j}_\eta q_\eta}
\prod_{j\in I}|f_j|^{2\nu^{j}_{\eta}q_\eta}}}d\lambda<\infty$$
for all $\eta\ll 1$ (by the definition of the set $I$ and the construction of the metric on 
$\wh L_{\eta|S}$). Therefore, the lemma is proved.\hfill\qed

\bigskip

\subsection {\S 1.H The method of Siu and Hacon-McKernan}
\medskip
\noindent We have arrived now at the 
last step in our proof : {\sl for all $0< \eta\ll 1$, the section $s_\eta$ admit 
an extension on $\wh X$}. Once this is done, we just use the point A.4 of the approximation lemma 1.F.1, in order to infer the  existence of a $\bR$--section of the bundle 
$K_{\wh X}+ S+\wh L$, and then the relation 1.D.1, $iv)$ to conclude.

\medskip
\noindent In order to explain our approach in the simplest possible way, 
we consider next the {\sl approximation of the usual setup} of the extension of twisted pluricanonical sections ;
afterwards we will compare it with our current situation. We will use different notations, to avoid any confusion that may occur.

Let $\ol X$ be a projective, non-singular manifold ; let 
$$\cS:= \big\{\ol S, (\ol Y_j)_{j\in J}, \ol L, \ol A, T\big\}$$
be a set of objects on $\ol X$, with the following properties.
{\itemindent 8mm
\smallskip
\item {$\bf (U_0)$} $\ol S$ and $(\ol Y_j)_{j\in J}$ are non-singular hypersurfaces of 
$\ol X$, with normal crossings and such that $\forall j\neq k$ we have $\ol Y_j\cap \ol Y_k= \emptyset$ ; 
\smallskip
\item {$\bf (U_1)$} $\ol L$ is a $\bQ$-line bundle, which admits the following numerical decomposition
$$\ol L\equiv
\ol\Delta_1+ \ol\Delta_2$$
where 
$$\ol\Delta_1= \sum_{j\in J\setminus I}\alpha^j[\ol Y_j]$$
and 
$$\ol\Delta_2= \sum_{j\in I}\alpha^j[\ol Y_j]+ \ol \Lambda_2.$$
The numbers $\alpha^j$ above are rational and moreover $\alpha^j\in ]0, 1[$ ; also, 
$\ol \Lambda_2$ is a K\"ahler current, with logarithmic poles and 
rational cohomology class, whose generic Lelong number along $\ol S$  is equal to zero, i.e. 
$$\nu_{\ol S}(\ol \Lambda_2)= 0.$$
In this case the restriction $\ol \Lambda_{2|\ol S}$ is well defined and we denote by 
$$\displaystyle \nu := \sup_{x\in \ol S} \nu (\ol \Lambda_{2|\ol S}, x).$$
\smallskip
\item {$\bf (U_2)$} $T$ is a closed current with logarithmic poles 
in the cohomology class of the bundle 
$K_{\ol X}+ \ol S+ \ol L$ and $\ol A$ is an ample line bundle 
such that :
\smallskip
\item\item {$\bf (U_{2.1})$} We have 
$$T\geq -C_T\Theta(A)\leqno (\cR_T)$$
where $C_T$ is a positive constant, and $\Theta(A)$ is a positive and non-singular 
curvature form corresponding to a metric on $A$.
\smallskip
\item\item {$\bf (U_{2.2})$} The restriction of $T$ to $\ol S$ is well-defined 
and we have
$$T_{|\ol S}:= \sum_{j\in J}\ol \theta^j[\ol Y_{j|\ol S}]+ R$$
where $\ol\theta^j$ are positive real numbers, and $R$ is a closed current on $\ol S$,
whose generic Lelong numbers along $\ol Y_{j|\ol S}$ is equal to zero ;
\smallskip
\item {$\bf (U_3)$} Let $q\in \bZ_+$ be a positive integer, such that $q \ol L$ and
$q\alpha^j$ are integral, for any $j\in J$, and such that 
$$q\nu< 1\leqno (\cR_\nu)$$
that is to say, the multiple $q$ is not allowed to be greater than
the inverse of the maximum of the singularities of $\ol \Lambda_{2|S}$.
We assume that there exist a section 
$$s\in H^0\big(\ol S, q(K_{\ol S}+ \ol L)\big)$$
whose zero set contains the divisor
$$q\big(\sum_{j\in J\setminus I}\rho^j[ \ol Y_{j|\ol S}]+ \sum_{j\in I}\alpha^j[ \ol Y_{j|\ol S}]\big)$$
such that 
$$\rho^j\geq \theta^j- C_T^\prime\leqno (\cR_s)$$
for each $j\in J\setminus I$ ; in the above relation, $C_T^\prime$ is a positive constant.
\hfill\qed
\medskip
\noindent Since the coefficients of $q\ol \Delta_1$ are positive integers strictly less than $q$, we 
have the decomposition (see [4], [18])
$$q\ol \Delta_1= L_1+...+L_{q-1}$$
such that for each $m= 1,..., q-1$, we have 
$$L_m:= \sum_{j\in I_m\subset J\setminus I}\ol Y_j.$$
We denote by $\displaystyle L_{q}:= q \ol \Delta_2$ and
$$L^{(p)}:= p(K_{\ol X}+ \ol S)+L_1+...+ L_p$$
where $p= 1,..., q$. By convention, $L^{(0)}$ is the trivial bundle.

\smallskip
\item {$\bf (U_4)$} The ample bundle $A$ is assumed to be positive enough such that the next conditions hold.
{
\itemindent 6mm
\noindent
\item\item {$(\dagger)$} For each $0 \leq p\leq q-1$, the bundle
$L^{(p)}+ q A$ is generated by its global sections, which we denote by
$(s^{(p)}_j)$.
\smallskip
\item\item {$(\dagger^2)$} Any section of the bundle $L^{(q)}+
qA_{\vert S}$ admits an extension to 
$\ol X$.
\smallskip
\item\item {$(\dagger^3)$} We endow the bundle corresponding to $(\ol Y_j)_{j\in J}$ and $K_{\ol X}+ \ol S$
with non-singular metrics, and we denote by $\wt \varphi_m$ the induced 
metric on $L_m$. 
Then for each $m= 1,..., q$, the functions
$$\wt \varphi_{L_m}+ 1/3\varphi_A \hbox{ and } \wt \varphi_{K}+ 1/3\varphi_A$$
are strictly psh, where $\wt \varphi_{K}$ is the non-singular metric on the bundle $K_{\ol X}+ \ol S$.
\hfill\qed

}}

\medskip
\noindent Under the numerous assumptions/normalizations above, we formulate the next statement.

\claim Claim|The section 
$$s^{\otimes k}\otimes s^{(p)}_j\in H^0\bigl(\ol S, L^{(p)}+ kL^{(q)}+ q A_{\vert \ol S}\bigr)$$
extend to $\ol X$, for each $p= 0,..., q-1$, $j= 1,..., N_p$ and $k\in \bZ_+$ such that 
$$k\max (C_T, C_T^\prime)\leq 1/4.$$
\hfill\qed
\endclaim

\noindent As one can see, the main differences 
between the present situation and the usual version of the invariance of plurigenera
(see e.g. [9], [18], [19], [24], [26], [33], [37], [38], [41], [42], [43], [44])
is visible in the relations $\cR_T, \cR_\nu$ and $\cR_s$. Indeed, the ``standard" assumptions in the articles quoted above are as follows.
\smallskip
\item {$(\bullet)$} The current $T$ in $\bf (U_2)$ is positive, that it to say $C_T= 0$ ;
\smallskip
\item {$(\bullet)$} The current $\ol \Lambda_2$ in $\bf (U_1)$ is non-singular ;
\smallskip
\item {$(\bullet)$} The restriction $T_{|\ol S}$ and the support of $\ol \Delta_1$
have no common components and therefore we can take $C_T^\prime= 0$.
\medskip 
\noindent Under this perspective, the statement above 
can be seen as a natural generalization of the usual
setting ; in substance, we are about to say that the more general hypothesis 
we are forced to consider induce an {\sl effective limitation} of the number of iterations
we are allowed to perform.

The rest of the present paragraph is organized as follows. 
We start with the proof of the claim above in 1.H.A.
In the paragraph 1.H.B we show that the family of bundles $(\wh L_\eta)$ is ``uniform", in the sense that we have a precise control of the corresponding constants 
$C_T, C^\prime_T$ as $\eta\to 0$. We complete the proof of the non-vanishing theorem in
1.H.C, by showing that the claim imply the extension of $(s_\eta)$.\hfill\qed


\bigskip
\subsection {\S 1.H.A Proof of the claim}

\noindent \noindent We will review here the main steps of the arguments in the usual invariance of plurigenera ; to start with, we recall the following very useful integrability criteria
(see e.g. [4]).

\claim 1.H.A.1 Lemma{\rm ([4])}|Let $\Theta$ be a (1,1)-current with logarithmic poles on a manifold $\ol S$, such that 
$\Theta\geq -C\omega$, where $C$ is a positive real, and $\omega$ is a metric 
on $\ol S$. We consider equally the non-singular hypersurfaces $\ol Y_j\subset \ol S$ for 
$j=1,...,N$ such that $\ol Y_j\cap \ol Y_i= \emptyset$ if $i\neq j$, and such that 
the generic Lelong number of $\Theta $ along each of the $\ol Y_j$ is zero.
Then there exist a constant $\varepsilon_0:= \varepsilon_0(\{\Theta\}, C)$ depending only on 
the cohomology class of the current $\Theta$ and its lower bound
such that 
for all positive real numbers $\delta\in ]0, 1]$ and
$\varepsilon \leq \varepsilon_0$ we have
$$\int_{\ol S}\exp \big(-(1-\delta)\sum_j\varphi_{\ol Y_j}- \varepsilon \varphi_\Theta\big)d\lambda< \infty. $$ \hfill\qed
\endclaim
\noindent In the statement above, we denote by $\varphi_{\ol Y_j}$ the potential of the current 
$[\ol Y_j]$. We remark that strictly speaking the quantity under the sum above is not global, but we are only interested in the {\sl singularities} of the objects above.

\medskip

\noindent We will equally need the following version of the
Ohsawa-Takegoshi theorem, obtained by McNeal-Varolin in [29] (see also [2], [13], [32], [38])~; it will be our main technical tool in the proof of the claim.
\medskip

\claim 1.H.A.2 Theorem ([29])|Let $X$ be a projective
$n$-dimensional 
manifold, and let $S\subset X$ be a non-singular hypersurface.
Let $F$ be a line bundle,
endowed with a metric $h_F$. We assume the existence 
of some non-singular metric $h_S$ on the bundle $\cal O(S)$ such that :  

{\itemindent 6mm
\smallskip
\item {(1)} $\displaystyle {{\sqrt {-1}}\over {2\pi}}\Theta_F
\geq 0$ on $X$~;
\smallskip
\item {(2)} $\displaystyle {{\sqrt {-1}}\over {2\pi}}\Theta_F -  
\alpha
 {{\sqrt {-1}}\over {2\pi}}\Theta_S\geq 0$ 
for some $\alpha> 0$~;
\smallskip
\item {(3)} The restriction
of the metric $h_F$ on $S$ is well defined.

}
\smallskip

\noindent Then every section $u\in H^0\bigl(S, (K_X + S+
F_{\vert S})\otimes {\cal I}(h_{F|S})\bigr)$ admits an extension
$U$ to $X$ such that 
$$c_n\int_XU\wedge \ol U\exp \Big(-\varphi_F-\varphi_S-\log\big(|s|^2\log^2(\vert s\vert)\big)\Big)< \infty$$
where $s$ is a section whose zero set is precisely the hypersurface $S$ and its norm in the integral above is measured with respect to $h_S$. \hfill\qed
\endclaim
\medskip

\noindent We will use inductively the extension theorem 1.H.A.2, in order to derive a lower bound for the power $k$ we can afford in the invariance of plurigenera algorithm, under the conditions 
$\bf (U_j)_{\rm 1\leq j\leq 4}$~; the first steps are as follows.

\noindent (1)  For each $j= 1,...,N_0$, the section 
$s\otimes s^{(0)}_j\in H^0\bigl({\ol S}, L^{(q)}+ q A_{\vert {\ol S}}\bigr)$
admits an extension $U^{(q)}_j\in H^0\bigl(\ol X,  L^{(q)}+ q A\bigr)$,
by the property $\bf (U_4)$, $\dagger\dagger$.

\noindent (2) We use the sections $(U^{(q)}_j)$ to construct a metric $\varphi^{(q)}$
on the bundle $L^{(q)}+ q A$.

\noindent (3) Let us consider the section 
$s \otimes s^{(1)}_j\in H^0\bigl({\ol S}, L^{(1)}+ L^{(q)}+ q A_{\vert {\ol S}}\bigr)$. 
We remark that the bundle
$$L^{(1)}+ L^{(q)}+ q A = K_{\ol X}+ {\ol S}+ L_1+ L^{(q)}+ q A$$
can be written as $K_{\ol X}+ {\ol S}+ F$ where
$$F:=  L_1+ L^{(q)}+ q A $$
thus we have to construct a metric on $F$ which satisfy the 
curvature and integrability assumptions in the Ohsawa-Takegoshi-type theorem
above.

Let $\delta, \varepsilon$ be positive real numbers~; we endow the bundle 
$F$ with the metric given by
$$\varphi^{(q)}_{\delta, \varepsilon}:= (1-\delta)\varphi_{L_1}+ \delta\wt \varphi_{L_1}+ 
(1-\varepsilon)\varphi^{(q)}+ \varepsilon q(\varphi_A+ \varphi_{T})\leqno (44)$$
where the metric $\wt \varphi_{L_1}$ is smooth (no curvature requirements) and  
$\varphi_{L_1}$ is the singular metric induced by the divisors 
$\displaystyle (Y_j)_{j\in I_1}$, see $\bf (U_4)$, $(\dagger ^3)$.

We remark that the curvature conditions in the extension theorem will be fulfilled if
$$\delta< \varepsilon q$$
since we are interested in the integers $k$ such that 
$$k\max (C_T, C_T^\prime)\leq 1/4,$$ 
thus implicitly $C_T< 1/3$ and the negativity of the curvature induced by the term $\delta\wt \varphi_{L_1}$ will be absorbed by $A$.
We use here the relations in $\bf (U_4)$, $(\dagger^3)$ and $(\dagger^4)$.

\noindent Next we claim that the sections $s\otimes s^{(1)}_j$ are integrable with respect to the metric defined in (44), provided that the parameters $\varepsilon, \delta$ are chosen in an appropriate manner.
Indeed, we have to prove that 
$$\int_{\ol S}{{\vert s\otimes s^{(1)}_j\vert^2}\over {(\sum_r \vert s 
\otimes s^{(0)}_r\vert^2)^{1-\varepsilon}}}
\exp \big(-(1-\delta)\varphi_{L_1}- \varepsilon q\varphi_{T}\big)dV< \infty ;$$
since the sections $(s^{(0)}_r)$ have no common zeroes, it is enough to show that 
$$\int_{\ol S}\vert s\vert ^{2\varepsilon}
\exp \big(-(1-\delta)\varphi_{L_1}- \varepsilon q\varphi_{T}\big)dV< \infty$$
(we have abusively removed the smooth weights in the above expressions, to simplify the writing). 

Now the property $(\bf U_3)$ concerning the zero set of $s$ is used : the above integral is convergent, provided that we have
$$\int_{\ol S}\exp \big(-(1-\delta)\varphi_{L_1}- 
\varepsilon q(\varphi_{T}- 
\sum_{j\in J\setminus I}\rho^{ j} \varphi_{Y_j}-  
\sum_{j\in I}\alpha^{j}\varphi_{Y_j})\big)dV< \infty.$$
We remark that we have
$$T_{|{\ol S}}- \sum_{j\in J\setminus I}\rho^{j} [Y_{j|{\ol S}}]-  
\sum_{j\in I}\alpha^{j}[Y_{j|{\ol S}}])\leq  \sum _{j\in I}(\ol\theta^{j}-\alpha^{j})_+[Y_{j|{\ol S}}]+ C_T^\prime\sum _{j\in J\setminus I}[Y_{j|{\ol S}}]+ R\leqno (\star)$$
(by the property $\bf (U_3)$)
and therefore, we have an explicit measure of the size of the common part 
of the 1-codimensional components of the difference current above 
and those of $L_1$. We use the notation
$$x_+:= \max(x, 0)$$
in the previous expression.
\smallskip
We check next the 
conditions on $\varepsilon, \delta$ in order to insure the hypothesis of 1.H.A.1.
The cohomology class of the current 
$$\sum _{j\in I}(\ol\theta^{j}-\alpha^{j})_+[Y_{j|{\ol S}}]+ R$$
together with 1.H.C.1 will provide us with a quantity which we denote by $\varepsilon_0$ 
(in the context we are interested in, it can even be assumed to be independent on $\eta$, but this 
does not matter).
The singular part corresponding to $j\in J\setminus I$ in the expression 
$(\star)$ will be 
incorporated into the $\displaystyle (1-\delta)\varphi_{L_1}$ therefore we impose
the relation
$$1-\delta+ q\varepsilon C_T^\prime < 1.$$

In conclusion, the positivity and integrability conditions will be satisfied provided that
$$q\varepsilon C_T^\prime <\delta < \varepsilon q\leq 
\varepsilon_0\leqno (45)$$
We can clearly choose the parameters $\delta, \varepsilon$ in order to satisfy (45), again by the 
implicit assumption in the claim.
\smallskip

\noindent (4) We apply the extension theorem and we get
$U^{(q+ 1)}_j$, whose
restriction on ${\ol S}$ is precisely
$s\otimes s^{(1)}_j$.\hfill \qed
\medskip

\vskip 5pt The claim will be obtained by iterating
the
procedure (1)-(4) several times, and estimating carefully the influence of the negativity of
$T$ on this process. Indeed, assume that we already have the set of global sections 
$$U^{(kq+p)}_j\in H^0\bigl(\ol X, L^{(p)}+ kL^{(q)}+ q A\bigr)$$
which extend $s^{\otimes k}\otimes s^{(p)}_j$. They induce a metric on the above bundle, denoted by $\varphi^{(kq+p)}$. 

\noindent If $p< q-1$, then we define the family of 
sections 
$$s^{\otimes k}\otimes s^{(p+1)}_j\in H^0({\ol S}, L^{(p+1)}+ kL^{(q)}+ 
q A_{|{\ol S}})$$
on ${\ol S}$. 
As in the step (3) above we remark that we have
$$L^{(p+1)}= K_{\ol X}+ {\ol S}+ L_{p+1}+ L^{(p)}$$
thus according to the extension result 1.H.A.2, we have to exhibit a metric
on the bundle 
$$F:=  L_{p+1}+ L^{(p)}+kL^{(q)}+ q A$$
for which the curvature conditions are satisfied, and such that the family of sections above are $L^2$ with respect to it.
We define
$$\varphi^{(kq+p+1)}_{\delta, \varepsilon}:= (1-\delta)\varphi_{L_{p+1}}+ \delta\wt \varphi_{L_{p+1}}+ 
(1-\varepsilon)\varphi^{(kq+ p)}+ \varepsilon q\big(k\varphi_{T}+  \varphi_A+ {{1}\over {q}}\wt \varphi_{L^{(p)}}\big)\leqno (46)$$
and we check now the conditions that the parameters $\varepsilon, \delta$ have 
to satisfy.
\medskip

We have to absorb the negativity in the smooth curvature terms in (46), 
and the one from $T$. The Hessian of the term
$$1/3\varphi_A+ {{1}\over {q}}\wt \varphi_{L^{(p)}}$$
is assumed to be positive by $\bf (U_4)$, $\dagger^3$, but we also have a negative contribution 
$$-kC_T\Theta_A$$
induced by the current $T$. However we remark that we have 
$$kC_T< 1/3\leqno (47)$$
by the hypothesis of the claim,
and then the curvature of the metric defined in (44) will be positive, provided that 
$$\delta< \varepsilon q$$
again by $(\dagger ^3)$. 

\medskip 
\noindent 
Let us check next the $L^2$ condition~; we have to show that 
the integral below in convergent 

$$\int_{\ol S}{{|s^{\otimes k}\otimes s^{(p+1)}_j|^2}\over 
{(\sum_r \vert s^{\otimes k}\otimes s^{(p)}_r\vert^2)^{1-\varepsilon}}}
\exp \big(-(1-\delta)\varphi_{L_{p+1}}- kq 
\varepsilon\varphi_{T}\big)dV.$$
This is equivalent with 
$$\int_{\ol S}|s|^{2\varepsilon k}
\exp \big(-(1-\delta)\varphi_{L_{p+1}}- kq 
\varepsilon\varphi_{T}\big)dV< \infty.$$
In order to show the above inequality, we use the same trick as before : the vanishing set of the section 
$s$ as in $\bf (U_3)$ will allow us to apply the integrability lemma--the computations 
are strictly identical with those discussed in the point 3) above, but we give here some details. 

By the vanishing properties of the section $s$, the finiteness of the previous integral will be implied by the inequality

$$\int_{\ol S}\exp \big(-(1-\delta)\varphi_{L_{p+1}}- 
k\varepsilon q(\varphi_{T}- 
\sum_{j\in J\setminus I}\rho^{ j} \varphi_{Y_j}-  
\sum_{j\in I}\alpha^{j}\varphi_{Y_j})\big)dV< \infty.$$

In the first place, we have to keep the poles of $k\varepsilon qT$ ``small" in the expression of the metric (46), thus we impose 
$$k\varepsilon q\leq \varepsilon_0.$$
 \smallskip
The hypothesis in the integrability lemma also require 
$$1-\delta+ \varepsilon kq C_T^\prime < 1$$
because of the contribution of the common part of $\Supp L_{p+1}$ and $T$.
Combined with the previous relations, the conditions for the parameters become
$$ \varepsilon kqC_T^\prime< \delta< \varepsilon q< \varepsilon_0/k.$$

In conclusion,  we can choose the parameters
$\varepsilon, \delta$ so that the integrability/positivity conditions in the extension theorem are verified ; for example, we can take
\smallskip
\noindent $\bullet$ $\displaystyle \varepsilon:= {{\varepsilon_0}\over {2kq}}$

\noindent and 
\smallskip
\noindent $\bullet$ $\displaystyle \delta:= (1+ kC_T^\prime){{\varepsilon_0}\over
{4 k }}$.

\medskip
\noindent Finally, let us indicate how to perform the induction step if $p= q-1$ :
we consider the family of 
sections 
$$s^{k+1}\otimes s^{(0)}_j\in H^0({\ol S}, (k+1)L^{(q)}+ q A_{|{\ol S}}),$$ 
In the case under consideration, we have to exhibit a metric
on the bundle 
$$L_{q}+ L^{(q-1)}+kL^{(q)}+  q A~;$$
however, this is easier than before, since we can simply take 
$$\varphi^{q(k+1)}:= 
 q \varphi_{\ol \Delta_2}+ 
\varphi^{(kq+q-1)}\leqno (48)$$
where the metric on $\ol \Delta_2$ is induced by the decomposition in 
$(\bf U_1)$.
With this choice, the curvature conditions are satisfied~; as for the $L^2$ ones, we remark that we have
$$\eqalign{
 & \int_{\ol S}{{|s^{k+1}\otimes s^{(0)}_j|^2}\over {(\sum_r \vert s^k\otimes s^{(q-1)}_r\vert^2)}}
\exp \big(-q\varphi_{\Delta_2}\big)dV< \cr
 < & C\int_{\ol S}{{|s\otimes s^{(0)}_j|^2}\over {(\sum_r \vert s^{(q-1)}_r\vert^2)}}
\exp \big(-q\varphi_{\Delta_2}\big)dV\leq C\int_{\ol S} 
\exp \big(-q\varphi_{\wh \Lambda_{2}}\big)dV< \infty \cr
}$$
where the last relation holds because the vanishing properties of $S$ and of the fact that
$q\nu< 1$ (cf. $\bf (U_3)$.
The proof of the extension claim is therefore finished.\hfill \qed

\bigskip


\medskip 
\subsection {\S 1.H.B Uniformity properties of $(K_{\wh X}+ S+ \wh L_\eta)_{\eta > 0}$}

\medskip
\noindent We come back now to  
the family of approximations 
$(\wh L_\eta)$ we have produced in the paragraph 1.F ; our main concern in this 
paragraph will be to derive its uniformity properties, in order to apply the previous considerations.
We list them below ; the constant $C$ which appear in the next 
statement is independent of $\eta$.
{\itemindent 7.5mm 
\smallskip
\item {$\bf (P_1)$} There exist a closed (1,1)--current $\Theta_\eta\in \{K_{\wh X}+ S+\wh L_\eta\}$ such that :
\smallskip
\item\item {$\bf (P_{1.1})$} It has logarithmic poles~;
\smallskip
\item\item {$\bf (P_{1.2})$} It is greater than $\displaystyle -C{{\eta}\over {q_\eta}}\omega$~;
\smallskip
\item\item {$\bf (P_{1.3})$} Its restriction to $S$ is well defined, and we have
$$\Theta_{\eta|S}= \sum _{j\in J}\theta^j_\eta[Y_{j|S}]+ R_{\eta, S}.$$
Moreover, the generic Lelong number of $R_{\eta, S}$ along $Y_{j|S}$ is zero, for any $j\in J$ and $\theta^j_\eta\leq \rho^{\infty, j}_\eta+ C{{\eta}\over {q_\eta}}$.
\smallskip

\item {($\bf P_2$)} The bundle $\wh L_{\eta}$ can be endowed with a metric whose curvature current is given by 
$$\sum_{j\in J}\nu^{j}_\eta[Y_j]+ 
\wh \Lambda_{L}+ \Xi(\eta)$$
where $\Xi(\eta) $ is non-singular and greater than
$-C\eta/q_\eta$~; we equally have $Y_j\cap Y_i=\emptyset$, if $i\neq j$.
As for $\wh \Lambda_{L}$, we have the next important property : 
there exist a smooth, projective manifold $\wt X$ and a birational map
$\mu_1:\wh X\to \wt X$ such that  
 $$\wh \Lambda_{L}\geq \mu_1^\star\wt \omega$$
 where $\wt \omega$
 is a K\"ahler metric on $\wt X$ and $S$ {\sl is not $\mu_1$-exceptional}
 (see the proposition 1.D.1.).

\smallskip

\item {($\bf P_3$)} The section $s_\eta\in H^0\big(S, q_\eta(K_S+ \wh L_\eta)\big)$
vanishes along the divisor
$$q_\eta\Big(\sum_{j\in J\setminus I}\rho^{\infty, j}_\eta[Y_{j|S}]+  
\sum_{j\in I}\nu^{j}_{\eta}[Y_{j|S}]\Big)$$ 
for all $0<\eta\ll 1$ \hfill\qed

}
\bigskip

\noindent The property $\bf (P_3)$ is a simple recapitulation of facts which were completely proved during the previous paragraphs ; let us give few details concerning 
$\bf (P_1)$ and $\bf (P_2)$.
\smallskip
By using the regularization theorem of Demailly, we have obtained in the section 1.E a family of currents with log poles 
$\Theta_\varepsilon\in \{K_{\wh X}+ S+\wh L\}$ such that $\Theta_\varepsilon\geq -\varepsilon \omega$, and whose restriction to $S$ equals
$$\Theta_{\varepsilon|S}= \sum_{j\in J}\rho^{\varepsilon, j}Y_{j|S}+ R_{\varepsilon}.$$
We first ``move" the current $\Theta_\varepsilon$ in the class $\{K_{\wh X}+ S+\wh L_\eta\}$
by a smooth form ; we denote by $\Theta_{\varepsilon, \eta}$ the result, and we remark that we have 
$$\Theta_{\varepsilon, \eta}\geq -\big(\varepsilon+ C{{\eta}\over {q_\eta}}\big)\omega$$
by the approximation lemma 1.F.1. 
Next,  we take $\varepsilon $ small enough such that :
\smallskip
\noindent $\bullet$ $\displaystyle \varepsilon\leq C{{\eta}\over {q_\eta}}$ ;
\smallskip
\noindent $\bullet$  $|\rho^{\infty, j}_\eta- \rho^{\varepsilon, j}|\leq C{{\eta}\over {q_\eta}}$ for each $j\in J$.
\medskip 

\noindent With this choice, the corresponding current will be our $\Theta_\eta$ ;
it satisfy all the requirements in the $\bf (P_1)$ above. \hfill\qed

\smallskip
The curvature current in the first Chern class of $\wh L_\eta$ can be obtained as in the section 1.G ; nevertheless, we will give here the full details of the construction.
We consider the decomposition
$$\wh L_\eta= L+ (\wh L_\eta-L)$$
and we remark that by the approximation lemma 1.F.1, the second term of the previous decomposition can be endowed with a non-singular metric whose curvature form is greater than $-C\eta/q_\eta$. We replace the coefficients $\nu^j$ of $L$ (see the property 1.D.1) with their rational approximations $\nu^j_\eta$ (cf. 1.F.1), and the negativity of the error
term is again dominated by $-C\eta/q_\eta$ ; thus, the claim  $\bf (P_2)$ is verified. \hfill\qed

\bigskip
\noindent We will establish now the relations between $\bf (U_0)-(U_4)$ and $\bf (P_1)-(P_3)$.
\smallskip
\noindent $\bullet$ $\ol X:= \wh X$, $\ol S= S$ and $\ol L:= \wh L_\eta$ ; the corresponding  coefficients/forms are $\alpha^j:= \nu^j_\eta$ and $\{\ol \Lambda_2\}:= 
\{\wh \Lambda_L+ \Xi(\eta)\}$ (we will choose later the representative corresponding to 
$\ol \Lambda_2$).
\smallskip
\noindent $\bullet$ The hypersurfaces $(Y_j)_{j\in J}$ are the same for all $\eta> 0$. 
\smallskip
\noindent $\bullet$ $T:= \Theta_\eta$, therefore
 $\ol \theta^j:= \theta^{j}_\eta$.
 \smallskip
\noindent $\bullet$ $q:= q_\eta$ and $s:= s_\eta$ ; thus we have $C_T^\prime:= C{{\eta}\over {q_\eta}}$, where we insist on the fact that $C$ {\sl does not depend on $\eta$}.

\medskip
\noindent In order to apply the claim we still have to clarify the following points : we have to 
identify the decomposition of $\wh L_\eta$ as in $\bf (U_1)$ and 
``trade" $\wh \Lambda_L+ \Xi(\eta)$ for a K\"ahler current, with an estimate for $\nu$ as it is required 
by $\bf (U_1)$ and $\bf (U_3)$. Also, we have to choose the bundle $A$ with the property $\bf (U_4)$.
This will be discussed along the next lines.
\smallskip

We first remark that we have the next decomposition  
$$\wh L_\eta\equiv  \sum_{j\in J\setminus I}\nu^{j}_\eta[Y_{j}]+  
\sum_{j\in I}\nu^{j}_\eta[Y_{j}]+
\wh \Lambda_{L}+ \Xi(\eta)~;$$
the coefficients corresponding to the indexes $j\in J$ are positive, as soon as $\eta$ is small enough. 

\noindent Let us introduce the next notations :
\smallskip
{\itemindent 2mm
\item {\noindent $\bullet$} $\Delta_1:= \sum_{j\in J\setminus I}\nu^{j}_\eta[Y_{j}]$. It 
is an effective $\bQ$--bundle whose critical exponent is greater than 1, and such that 
the multiple $q_\eta \nu^{j}_\eta$ is a positive integer strictly smaller than $q_\eta$, for each $j\in J\setminus I$~;
\smallskip
\item {\noindent $\bullet$} $\Delta_2:= 
\sum_{j\in I}\nu^{j}_\eta[Y_{j}] + \wh \Lambda_{L}+ \Xi(\eta) $. It is equally a $\bQ$-bundle whose critical exponent 
is greater than 1 and such that $q_\eta \Delta_2$ is integral. One of the facts which will be relevant in what follows is that the section $s_\eta$
vanishes along the singular part of $q_\eta\Delta_2$. 

}
\medskip
\noindent By the property $\bf (P_2)$, we can find a representative of the class
 $\{ \wh \Lambda_{L}\}$ which dominates a K\"ahler metric (see also the previous section) ; in general we cannot avoid that this representative acquire some singularities.
 However, in the present context 
 we will show that there exist a K\"ahler current in the above class which is ``restrictable" to $S$. It is at this point that 
 we use the full force of the property $\bf (P_2)$. 
 
Indeed, we consider the exceptional divisors $(E_j)$ of the map
$\mu_1$ (see the paragraph 1.D)~; the hypersurface $S$ do not belong to this set, and 
then the class 
$$\wh \Lambda_{L}-\sum_{j}\varepsilon^jE_j$$
is ample on $\wh X$, for some positive reals $\varepsilon^j$. Once a set of such parameters is chosen, we fix a K\"ahler form
$$\Omega\in \{\wh \Lambda_{L}-\sum_{j}\varepsilon^jE_j\}$$
and for each $\delta\in [0, 1]$ we define 
$$\wh \Lambda_{L, {\delta}}:= (1-\delta)\wh \Lambda_{L}+ \delta \big(\Omega+
\sum_j\varepsilon^jE_j\big)\in \{\wh \Lambda_{L}\}.\leqno (49)$$
For each $\eta > 0$, there exist $\delta> 0$ such that 
\smallskip
\item {$M_3)$} The current $\wh \Lambda_{L, {\delta}}+ \Xi(\eta)$ dominates
a K\"ahler form on  $\wh X$~;
\smallskip
\item {$M_4)$} The Lelong number of the restriction $\wh \Lambda_{L, {\delta}|S}$ at each point of $S$
does not exceed $\displaystyle C{{\eta}\over {q_\eta}}$, where $C$ is a constant independent of $\eta$.

\smallskip
\noindent  One can take $\delta:= \varepsilon_0^{-1}\displaystyle {{\eta}\over {q_\eta}}$ where $\varepsilon_0$ is small (but fixed)
and the properties $M_3)$ and $M_4)$ are clearly satisfied,
since the negative part of $\Xi(\eta)$ is given 
by the property $\bf (P_2)$ above. We denote by
$\wh \Lambda_{L, {\eta}}$ the corresponding current, and we observe that it satisfy
the next property as well.
\smallskip
\item {$M_5)$} The restriction of $\wh \Lambda_{L, {\eta}}$ to $S$ is well defined.

\medskip
\noindent Thus, the $\bQ$--divisor $\Delta_2$ is linearly equivalent to

$$\Delta_2\equiv 
\sum_{j\in I}\nu^{j}_\eta[Y_{j}] + 
 \wh \Lambda_{L, \eta}+ \Xi(\eta) ;\leqno (50)$$ 
it is a K\"ahler current and we have $q_\eta\max \nu(\wh \Lambda_{L, \eta|S}, s)\ll1$, for $\eta$ small enough ;
therefore, the first inequality in $\bf (U_3)$ is satisfied.

\medskip
\noindent Precisely as in [4], [15], [18] there exist a decomposition
$$q_\eta\Delta_1= L_1+...+L_{q_\eta-1}$$
such that for each $m= 1,..., q_\eta-1$, we have 
$$L_m:= \sum_{j\in I_m\subset J\setminus I}Y_j.$$

We denote by $\displaystyle L_{q_\eta}:= q_\eta\Delta_2$ and
$$L^{(p)}:= p(K_X+ S)+L_1+...+ L_p\leqno (51)$$
where $p= 1,..., q_\eta$. By convention, $L^{(0)}$ is the trivial bundle.
\medskip
\noindent Finally, it is possible to find an ample bundle $(A, h_A)$ {\sl independent of $\eta$} whose curvature form is positive enough such that the next relations hold.

{
\itemindent 6mm
\noindent
\item {$(\dagger)$} For each $0 \leq p\leq q_\eta-1$, the bundle
$L^{(p)}+ q_\eta A$ is generated by its global sections, which we denote by
$(s^{(p)}_j)$.
\smallskip
\item {$(\dagger^2)$} Any section of the bundle $L^{(q_\eta)}+
q_\eta A_{\vert S}$ admits an extension to 
$\wt X$.
\smallskip
\item {$(\dagger^3)$} We endow the bundle corresponding to $(Y_j)_{j\in J}$
with a non-singular metric, and we denote by $\wt \varphi_m$ the induced 
metric on $L_m$. 
Then for each $m= 1,..., q_\eta$, the functions
$$\wt \varphi_{L_m}+ 1/3\varphi_A$$
are strictly psh.
\smallskip
\item {$(\dagger^4)$} For any $\eta > 0$ we have
$$\Theta_\eta\geq -{{\eta}\over {q_\eta}}\Theta_A.  $$
thus $C_T:= {{\eta}\over {q_\eta}}$.

}
\medskip
\noindent In other words, the bundle $``A"$ in ${\bf U_4}$ will be $A$ in the present context.\hfill\qed

\bigskip
\claim Remark|{\rm Concerning the construction and  the properties 
of $\wh \Lambda_{L, \delta}$, we recall the very nice result in [17], stating that 
{\sl if $D$ is an $\bR$-divisor which is nef and big, then its associated augmented base locus can be determined numerically.} 
}\endclaim
\bigskip

\subsection {\S 1.H.C End of the proof}

\smallskip
\noindent We show next that the sections $s_\eta$ can be lifted to $\wh X$
as soon as $\eta$ is small enough, by using the claim proved in 1.H.A.    

Indeed, we consider the extensions $U^{(kq_\eta)}_j$ of the sections 
$s_\eta ^{\otimes k}\otimes s^{(0)}_j$~; they can be used to define a metric on the bundle
$$kq_\eta(K_{\wh X}+S+ \wh L_\eta)+ q_\eta A$$
whose $kq_\eta^{\rm th}$ root it is defined to be 
$h^{(\eta)}_k$. 

As usual, we write the bundle we are interested in i.e. 
$q_\eta(K_{\wh X}+S+ \wh L_\eta)$ as an adjoint bundle ; 
we have
$$\eqalign{ q_\eta(K_{\wh X}+S+ \wh L_\eta)= & K_{\wh X}+ S+ (q_\eta- 1)(K_{\wh X}+S+ \wh L_\eta)+ \wh L_\eta= \cr
= & K_{\wh X}+ S+ (q_\eta- 1)\big(K_{\wh X}+S+ \wh L_\eta+ 1/kA\big)+ \wh L_\eta- {{q_\eta-1}\over 
{k}}A \cr
}
$$

\medskip 
\noindent Given the extension theorem 1.H.A.2, we need to construct a metric on the bundle
$$(q_\eta- 1)\big(K_{\wh X}+S+ \wh L_\eta+ 1/kA\big)+ \wh L_\eta- {{q_\eta-1}\over 
{k}}A. $$
On the first factor of the above expression we will use 
$(q_\eta- 1)\varphi_k^{(\eta)}$ (that is to say, the $(q_\eta-1)^{\rm th}$ power of the metric given by
$h_k^{(\eta)}$). 

We endow the bundle $\wh L_\eta$ with a metric whose curvature is given by the expression 
$$\sum_{j\in J\setminus I}\nu^{j}_\eta[Y_{j}]+  
\sum_{j\in I}\nu^{j}_\eta[Y_{j}]+
\wh \Lambda_{L, \delta}+ \Xi(\eta)~;$$
here 
we take $\delta$ {\sl independent of $\eta$}, but small enough such that the critical exponent of the 
resulting metric on $\wh L_{\eta|S}$ is still greater than 1.
Finally, we multiply with the ${{q_\eta-1}\over 
{k}}$ times $h_A^{-1}$. 

\smallskip The corresponding constants $C_T$, respectively $C^\prime_T$ in 1.H.A are in the present context $\displaystyle {{\eta}\over {q_\eta}}$, respectively $\displaystyle C{{\eta}\over {q_\eta}}$ ; thus by the claim, we 
are free to choose $k$ e.g. such that $k= q_\eta\big[\eta^{-1/2}\big]$
(where $[x]$ denotes the integer part of the real $x$). Then the metric above is not identically $\infty$ when restricted to $S$, and its 
curvature will be 
strongly positive as soon as $\eta\ll 1$. Indeed, the curvature of $\wh L_\eta$ is greater 
than a  K\"ahler metric on $\wh X$ {\sl which is independent of $\eta$}
because of the factor $\wh \Lambda_{L, \delta}$. 

Moreover, the $L^2$ conditions in the theorem 1.H.A.2 are satisfied, since 
the norm of the section $s_\eta$ with respect to the metric $q_\eta\varphi_k^{(\eta)}$
is {\sl pointwise bounded}, and since the critical exponent of the metric on 
$\wh L_{\eta|S}$ is greater than 1. 
\medskip
\noindent In conclusion, we obtain an extension of the section $s_\eta$, and the theorem 0.1 is completely proved.\hfill\qed

\medskip
\claim Remark|{\rm The {\sl exact} vanishing properties of 
$s_\eta$ given by the proposition 1.G.1 are crucial. Indeed, assume that instead of 
the divisor
$$q_\eta\Big(\sum_{j\in J\setminus I}\rho^{\infty, j}_\eta[Y_{j|S}]+  
\sum_{j\in I}\nu^{j}_{\eta}[Y_{j|S}]\Big)$$ 
the section $s_\eta$ only vanishes along
$$(q_\eta-1)\Big(\sum_{j\in J\setminus I}\rho^{\infty, j}_\eta[Y_{j|S}]+  
\sum_{j\in I}\nu^{j}_{\eta}[Y_{j|S}]\Big).$$ 
This may look {\sl innocent}, since $q_\eta\to \infty$ anyway, but 
we remark that under these circumstances we have
$$C_T^\prime= C{{\eta}\over {q_\eta}}+ \max_j\Big\{{{\rho^{\infty, j}}\over {q_\eta}}\Big\}$$
and the whole extension process collapse, since 
we cannot insure 
$${{q_\eta}\over {k}}\to 0$$
anymore. \hfill\qed

}

\endclaim

\bigskip

\section{\S 2. Metrics with minimal singularities and holomorphic sections}

\medskip
\noindent In this section we will consider the following geometric context. Let $X$ be a non-singular, projective $n$-dimensional manifold, and let 
$L\to X$ be a $\bQ$-line bundle with the following metric property :
\smallskip
\item {($\bullet$)} The (1,1)--class $c_1(L)$ contains a K\"ahler current 
$\Theta_L$ whose critical exponent is strictly greater than 1. 

\medskip
\noindent We consider the following graded algebra
$$\cR(X, L):= \bigoplus_{k\in q\bZ_+}H^0\big(X, k(K_X+ L)\big)\leqno (52)$$
where $q\in \bZ_+$ is a positive integer such that 
$qL$ is a line bundle. In order to study its properties,
we assume that a non-singular metric $\wt h$ on $K_X+ L$ is given, and 
for any 
$$u\in H^0\big(X, k(K_X+ L)\big)$$
we will denote by 
$|u|_k$ the norm of $u$ with respect to the metric $\wt h^{\otimes k}$. \hfill\qed
\medskip

In the paragraph 2.A below we 
define {\sl a metric with 
minimal singularities} which is adapted to the ring $\cR(X, L)$. For some technical reasons (which will only appear in the second paragraph 2.B), we are forced to take into account the sections of the multiples of $K_X+ L$
{\sl twisted with a topologically trivial line bundle},
even if our ultimate goal would be to understand the structure of the ring $\cR(X, L)$.

Next, we will consider the relative threshold of the minimal metric with respect to a metric given by a finite number of sections of the multiples of $K_X+ L$ ; 
if the singularities of these metrics do not coincide, then the non-vanishing theorem will provide us with an
$\bR$--section of $K_X+L$ 
which has precisely the same vanishing order along {\sl some} divisor, say $S$, of a modification of $X$ as the metric with minimal singularities (compare with [39]). 
The results we obtain in the paragraph 2.A will show in particular that the said vanishing order is
a rational number. \hfill\qed
\subsection {\S 2.A  Metrics with minimal singularities}
\smallskip
\noindent Along the following lines, we will only consider the case of adjoint bundles, since it is in this setting that the main properties we will establish afterwards hold ; 
however, one could define the objects below in a more general context. A general reference for the notions discussed in this subsection is the article [14].

Let $\rho\to X$ be a topologically trivial line bundle, endowed with a metric
$h_\rho$ whose curvature form is equal to zero. Given 
$$u\in H^0\big(X, k(K_X+ L)+\rho\big)$$
we denote by $|u|_{k,\rho}^2$ the poinwise norm of $u$, measured with the 
$\wt h^{k}$ twisted with $h_\rho$.

\noindent We introduce the following class of functions on $X$
$$\cF:= \Big\{f= {1\over k}\log|u|_{k,\rho}^2 :  
 k\in q\bZ_+, u\in H^0\big(X, k(K_X+ L)+\rho\big), \hbox{ s.t. }\sup_X|u|_{k,\rho}^2= 1  \Big\},$$
and we remark that we have 
$$\Theta_{\wt h}(K_X+ L)+ {{\sqrt {-1}}\over {2\pi}}\ddbar f\geq 0$$
for any $f\in \cF$ (since the curvature of $\rho$ with respect to $h_\rho$ is equal to zero).
Thus the curvature current associated to the metric 
$$\exp(-f)\wt h$$
on the bundle $K_X+ L$ is positive.

We will denote by $h_{\rm min}$ the metric on $K_X+ L$ given 
by the smallest upper semicontinuous
majorant of the family $\cF$ above (see [14]). We denote by 
$$\Theta_{\rm min}:= \Theta_{h_{\rm min}}(K_X+ L)$$
the curvature current associated to this metric ; by definition we have
$$|v|^2\exp(-k\varphi_{\rm min})\leq O(1)$$
for any $v\in H^0\big(X, k(K_X+ L)+\rho\big)$. 
\medskip

\noindent Next {\sl we fix} a topologically trivial line bundle $\rho_0$ on $X$ and an integer 
$m_0\in q\bZ_+$ ; we will construct another metric on the bundle $K_X+ L$
as follows. We consider the set of potentials
$$\cF_0:= \Big\{f= {{1}\over {km_0}}\log|u|_{km_0,k\rho_0}^2 : u\in H^0\big(X, km_0(K_X+ L)+k\rho_0\big), \sup_X|u|_{km_0,k\rho_0}^2= 1  \Big\}$$
and we denote by $h_{\rm min}^{\rho_0}$ the metric on $K_X+ L$ given 
by the smallest upper semicontinuous
majorant of the family $\cF^{\rho_0}$ ; let $\Theta_{\rm min}^{\rho_0}$ be the associated curvature current.
\medskip
\noindent We state now the main result of this subsection.

\claim 2.A.1 Theorem|We have 
$$h_{\rm min}= h_{\rm min}^{\rho_0}.$$
\hfill\qed 
\endclaim

\claim 2.A.2 Remark|{\rm If the adjoint bundle $K_X+ L$ is big, then the above result is a consequence of the regularization theorem [11]. Also, the results of Campana-Peternell 
(see [8]) suggest that the strict positivity of $L$ in the theorem above may be superfluous. 
}
\endclaim
\proof \hskip 1.5mm (of the theorem 2.A.1). The relation 
$$\varphi_{\rm min}\geq\varphi_{\rm min}^{\rho_0}$$
is implied by the definition ; in order to obtain an inequality in the opposite sense, 
we consider a section
$$u_1\in H^0\big(X, m_1(K_X+ L)+\rho_1\big)\leqno (53)$$ 
whose norm is smaller than 1.
\smallskip 
\noindent As a consequence of an argument due to Shokurov 
(already employed in the section 1.G), we have the next statement.
\claim 2.A.2 Lemma|For any $k\in \bZ_+$,  there exist a section
$$u_k\in H^0\big(X, km_1m_0(K_X+ L)+km_1\rho_0\big)\leqno (54)$$
such that the next integral condition is satisfied
$$\int_X{{c_ku_k\wedge \ol {u_k}\exp({-km_1\varphi_{\rho_0}})}\over 
{{(c_1u_1\wedge \ol{u_1}})^{{{km_1m_0-1}\over{m_1}}}\exp\big(-
{{km_1m_0-1}\over {m_1}}\varphi_{\rho_1}\big)}}\exp(-\varphi_L)d\lambda
= 1.\leqno (55)$$
\hfill\qed
\endclaim
\noindent  In the lemma above we denote by $\varphi_L$ the local weight of any metric on 
$L$ which satisfy the properties at the beginning of 2.A : the corresponding curvature 
is a K\"ahler current whose critical exponent is greater than 1. The quantities $c_k$ and $c_1$ are the usual ones, such that the wedge at the
denominator and numerator in the previous lemma are reals. We remark that the quantity under the integral sign in the above formula is a globally defined measure.
\smallskip
\noindent 
\proof. For any positive integer $k$, we consider the bundles
$$km_1m_0(K_X+ L)+ km_1\rho_0\leqno (56)$$
and
$$km_1m_0(K_X+ L)+ km_0\rho_1\leqno (57)$$
as well as the multiplier ideal
$$I_k:= \cI\big((km_1m_0-1)\log|u_1|^{{2}\over {m_1}}_{m_1, \rho_1}+ \varphi_L\big)\leqno (58)$$
We have
$$\chi\Big(\big(km_1m_0(K_X+ L)+ km_1\rho_0\big)\otimes I_k\Big)=
\chi \Big(X, \big(km_1m_0(K_X+ L)+ km_0\rho_1\big)\otimes I_k\Big)$$
by the usual arguments (i.e. the existence of a finite and free resolution of the multiplier sheaf $I_k$, and the 
fact that $\rho_j$ are topologically trivial).

Thanks to the Kawamata-Viehweg-Nadel vanishing theorem we equally know that the respective higher cohomology groups are equal to zero, so in conclusion
$$ H^0\Big(X, \big(km_1m_0(K_X+ L)+ km_1\rho_0\big)\otimes I_k\Big)=
H^0\Big(X, \big(km_1m_0(K_X+ L)+ km_0\rho_1\big)\otimes I_k\Big).\leqno (59)$$

\noindent We claim now that the section $u_1^{\otimes km_0}$ belong to the right hand side cohomology group in (54). To verify this claim, we have to show that the following integral converge
$$\int_X|u_1|^{{{2}\over {m_1}}}\exp (-\varphi_L)d\lambda,$$
and indeed this is the case, since the critical exponent of the current $\Theta_L$ is greater than 1.
\medskip 
\noindent Therefore we infer the existence of a non-identically zero section 
$$u_k\in H^0\Big(X, \big(km_1m_0(K_X+ L)+ km_1\rho_0\big)\otimes I_k\Big).\leqno (60)$$
By the construction of the ideal $I_k$, we see that we can normalize the section $u_k$ such that the integral condition (54) is satisfied ; the lemma 2.A.2 is therefore proved.
\hfill\qed
 \medskip
The finiteness of (54) show the existence of a section 
$$v_k\in H^0\big(X, m_1(K_X+ L)+ km_1\rho_0-km_0\rho_1\big)$$
such that 
$$u_k= u_1^{\otimes (km_0-1)}\otimes v_k ; \leqno (61)$$
the integral relation (54) become
$$\int_X{{c_{1}v_k\wedge \ol v_k}\over {{\big(c_{1}u_1\wedge\ol u_1\big)}^{{{m_1-1}\over{m_1}}}}}
\exp\big(-\varphi_L-km_1\varphi_{\rho_0}+ (km_0-1+ 1/m_0)\varphi_{\rho_1}\big)= 1.
\leqno (62)$$

\medskip
We will use the family of sections $(u_k)_{k\in \bZ_+}$ in order to compare 
$\varphi_{\rm min}$ and $\varphi_{\rm min}^{\rho_0}$. A specific normalization was chosen for the sections defining the potentials in $\cF$ ; thus, we
have to estimate the sup norm of $u_k$ along the next lines. Our main technical tools
will be the standard convexity properties of the psh functions ; we use the same notation
$``C" $ for all the constants which will occur during the following computations, even if they are not the same inside the same line, as long as they do not depend on $k$.

We define
$$\exp(f_k):= |v_k|^2_{m_1, km_1\varphi_{\rho_0}- (km_0-1)\varphi_{\rho_1}} ;\leqno (63)$$
the next step in our proof is to show the existence of a positive constant 
$C= C(m_1)$ large enough, so that we have 
$$\sqrt {-1} \ddbar f_k\geq -C\omega\leqno (64)$$
and moreover
$$-\log C\leq \max_Xf_k\leq \log C.\leqno (65)$$
The inequality (64) is a consequence of the fact that the curvature form of 
the metrics on the bundles $\rho_j$ is equal to zero ; 
let us give some explanations about (65). By using the 
notations introduced in (63), the equality (62) become
$$\int_X{{\exp\big(f_k-f_L)}\over {|u_1|_{m_1, \rho_1}^{2{{m_1-1}\over{m_1}}}}}dV_\omega
= 1.
\leqno (66)$$
where $f_L$ is the (global) distortion function between the
metric $\varphi_L$ and the non-singular metric on $L$ induced by $\wt h$ on $K_X+ L$ and 
$\det (\omega)$ on $-K_X$.
The section $u_1$ is normalized such that
$$\max_X(|u_1|_{m_1, \rho_1})= 1$$
and then the right hand side part of (65) is a consequence of the mean inequality
for the psh functions. 

In order to obtain the first inequality, 
we consider a log-resolution 
$\mu: \wh X\to X$ of the function
$$\displaystyle \psi:= f_L+ {{m_1-1}\over{m_1}}\log |u_1|_{m_1, \rho_1}^2$$ 
and we have
$$\psi\circ\mu:= \wh \psi+ \sum_{j\in J}a^j_\psi\log |s_j|^2$$
as well as
$$K_{\wh X/X}:= \sum_{j\in J}a^j_{\wh X/X}[W_j]$$
where $\wh \psi$ is a smooth function globally defined on
$\wh X$, and the hypersurfaces 
$$W_j:= (s_j= 0)$$
have normal crossings. We stress at this point on the fact that $\mu$ {\sl does not} depend
on $k$.

\noindent For each positive integer $k$ we decompose the inverse 
image of $f_k$ as follows
$$f_k\circ\psi:= \wh f_k+ \sum_{j\in J}a^j_k\log |s_j|^2$$
where $\wh f_k$ is non singular along any of $(W_j)$,
and we observe that there exist a metric $\wh \omega$ on $\wh X$ such that we have
$$\sqrt {-1}\ddbar \wh f_k\geq -\wh \omega\leqno (67)$$
for any $k$--this is a direct consequence of the relation (64).
With this notations, the equality (66) become
$$\int_X\exp\big(\wh f_k- \wh \psi+ \sum_{j\in J}(a^j_k+ a^j_{\wh X/ X}- a^j_\psi)\log |s_j|^2\big)dV_{\wh \omega}
= C
\leqno (68)$$
for some constant $C> 0$, uniform with respect to $k$
(which appear instead of 1 because of (67)).
The finiteness of the above integral show that we have
$$a^j_k+ a^j_{\wh X/ X}- a^j_\psi> -1$$
for any $j\in J$. We remark that the positive reals $a^j_k$ defined above are in fact {\sl integers}, since they correspond to the vanishing order of the inverse image of $v_k$ along $W_j$. Therefore we obtain
$$a^j_k\geq [a^j_\psi- a^j_{\wh X/ X}]$$
for all $j\in J$ and we re-write the formula (68) as follows
$$\int_X\exp\big(\wt f_k- \wt \psi\big)dV_{\wh \omega}
= C,
\leqno (69)$$
where we use the following notations
$$\wt f_k:= \wh f_k+ \sum_{j\in J}(a^j_k- [a^j_\psi- a^j_{\wh X/ X}])\log |s_j|^2$$
and
$$\displaystyle \wt \psi:= \wh\psi+ \sum_{j\in J}\{a^j_\psi- a^j_{\wh X/ X}\}\log |s_j|^2.$$
We remark that the equality (69) is similar to (66), but in addition multiplier ideal of the function $\wt \psi$ is trivial.

\noindent Let 
$$\wt C_k:= \max _{\wh X}(\wt f_k) ;$$
as a consequence of the relations (67) and (69), there exist a positive constant $\wt C$ such that
$$\wt C_k\leq C$$
for any $k\in \bZ_+$. The equality (69) imply
$$C\exp (-\wt C_k)= \int_X\exp\big(\wt f_k- \wt C_k-\wt \psi\big)dV_{\wh \omega}$$
and thus 
$$\exp (-\wt C_k) \leq C\int_X\exp\big(-\wt \psi\big)dV_{\wh \omega}:= C< \infty$$
In conclusion, the sequence $(\wt C_k)$ is bounded from below as well.
\smallskip
By the usual properties of the quasi-psh functions (see e.g. [11], [12]), there exist a function
$\wt f_\infty\in L^1(\wh X)$ such that 
$$\wt f_k\to  \wt f_\infty$$
as $k\to\infty$.
This show in particular the validity of the inequality (65), since
$$f_k\circ \mu= \wt f_k+ \sum_{j\in J}[a^j_\psi- a^j_{\wh X/ X}]\log |s_j|^2$$
and thus the sequence $(f_k)$ cannot tend to $-\infty$.

\hfill\qed
\smallskip
The important consequence of the previous considerations is the existence of a limit for the sequence $(f_k)$.
The normalization of $u_1$ and the relation (56) show that 
$$\max_X|u_k|_{km_1m_0, km_1\rho_0}\leq e^{C(m_1)}$$
and finally we get
$$\eqalign{
\varphi_{\rm min}^{\rho_0}\geq & -{{C(m_1)}\over {km_1m_0}}+ {{1}\over {km_1m_0}}\log |u_k|^2\geq\cr
\geq & -{{C(m_1)}\over {km_1m_0}}+ {{km_0-1}\over {km_1m_0}}\log |u_1|^2+ 
{{1}\over {km_1m_0}}f_k.\cr
}$$
We let $k\to \infty$ ; the first and the third term in the last inequality above tend to zero, and thus we get
$$
\varphi_{\rm min}^{\rho_0}\geq {{1}\over {m_1}}\log |u_1|^2.$$
The section $u_1$ above is arbitrary, thus the theorem is proved.
\hfill\qed

\bigskip

\subsection {\S 2.B Constructing sections with minimal vanishing order}
\medskip
\noindent In this paragraph we would like to point out an important property
of the zeroes of the $\bR$--sections produced by 0.1, in connection with Siu's 
proof of the finite generation problem (see [39], [40]). The same hypothesis/conventions
as in the beginning of the section are in force ; in addition, 
given an integer $\alpha$ large enough, we consider the following {\sl truncation metric}
$$\varphi_{\alpha}:= \log \Big(\sum_{k=1}^\alpha\varepsilon_k\!\sum_{j\in J_k}|f^k_j|^{2\over k}\Big)\leqno (70)$$
where $\varepsilon_k$ are positive real numbers, and
$(f^k_j)$ are local expressions of a family of sections of $k(K_X+ L)$~;
let $\Theta_{\alpha}$ be the corresponding current.

Precisely as in the paragraph 1.C, we will consider $\mu: \wh X\to X$
a log-resolution of the currents $\Theta_L$ and $\Theta_{\alpha}$~; we have

$$\mu^\star(\Delta)= \sum_{j\in J}a^{j}_\Delta[Y_{j}]\leqno (71)$$
and
$$\mu^\star(\Theta_L)= \sum_{j\in J}a^{j}_L[Y_{j}]+ \wh \Lambda_L\leqno (72)$$
as well as 
$$\mu^\star(\Theta_{\alpha})= \sum_{j\in J}a^{j}_\alpha [Y_{j}]+ \wh \Lambda_\alpha\leqno (73)$$
where $\wh \Lambda_L$, respectively $\wh \Lambda_\alpha$ are non-singular and semi-positive $(1,1)$--forms on $\wh X$ which are positively defined at the generic point of this manifold.
We consider next the inverse image of the minimal current via $\mu$ :
$$\mu^\star(\Theta_{\rm min})= \sum_{j\in J}a^{j}_{\rm min} [Y_{j}]+ 
\wh \Lambda_{\rm min}\leqno (74)$$
where $\wh \Lambda_{\rm min}$ is a closed positive current, whose generic Lelong numbers along the hypersurfaces $Y_j$ above is equal to zero.
Moreover, by the definition of the minimal metric we have 
$$a^j_\alpha\geq  
a^{j}_{\rm min}\leqno (75)$$
for all $j\in J$. We equally have the pointwise
inequality 
$$\varphi_{\wh \Lambda_{\rm min}}\geq \sum_{j\in J}\big(a^j_\alpha- 
a^{j}_{\rm min})\log |f_j|^2\leqno (76)$$
modulo an irrelevant constant. 

\medskip
\noindent The main result of the current subsection is the following.

\claim 2.B.1 Theorem|If at least one of the inequalities (75) is strict, then
there exist a topologically trivial line bundle $\rho\to X$ and a section $u\in H^0\big(X, m(K_X+ L)+\rho\big)$ 
such that the vanishing order of $\mu^\star(u)$ along 
$\displaystyle Y_{j_0}$ is precisely $ma^{j_0}_{\rm min}$ for some index $j_0\in J$. In particular,
we have $a^{j_0}_{\rm min}\in \bQ$.\hfill\qed
\endclaim

\medskip

\noindent We remark that if all the inequalities (75) are 
{\sl equalities}, then (76) show that all the local potentials of the current 
$\wh \Lambda_{\rm min}$ are bounded. In other words, the metric with minimal singularities is equivalent with its truncation $\varphi_\alpha$, and this 
imply the finite generation of the ring associated to $K_X+ L$, according to [39].
\medskip
\proof.  We consider the relative threshold associated to the following objects :
$$\tau:= \sup\{t\in \bR_+:  \int_{\wh X}\exp\big(t(\varphi_{\wh D}-\varphi_\alpha
)+ \varphi_{\wh D}+ \varphi_{\wt X/X}- \varphi_L\circ \mu\big)d\lambda<\infty \},\leqno (77)$$
where we use the notation
$$\wh D:= \sum_{j\in J}a^{j}_{\rm min} [Y_{j}].$$
We observe that $\tau$ verify the next relations
$$0< \tau< \infty\leqno (78)$$
by the same arguments as in the proof of 0.1--we remark that the latter inequality is a consequence of our assumption above.

The perturbation argument we have used in 1.C still apply in the present setting ; there exist a unique $S\subset \{Y_j\}$ such that we have the next relation
$$\leqno (79)
\eqalign{
\mu^\star\big( & K_X+ \tau(\Theta_\alpha-\Theta_{\rm min})+ 
\Theta_L- \Theta_{\rm min}\big)+ (1+\tau) \wh \Lambda_{\rm min}
\equiv \cr 
\equiv &K_{\wh X}+ S+ \sum_{j\in J}
\big(\tau (a^j_\alpha-a^j_{\rm min})+ a^j_L-a^j_{\rm min}-a^j_{\wh X/X}\big)[Y_j]+
\tau\wh \Lambda_{\alpha}+ \wh \Lambda_{L}\cr
}$$
where the coefficients of $Y_j$ are strictly smaller than 1, and the form 
$\wh \Lambda_{L}$ is positively defined. 
The relation (83) is equivalent with
$$\leqno (80)
\eqalign{
(1+\tau) \wh \Lambda_{\rm min}+ &
\sum_{j\in J_n}
\big(a^j_{\rm min}+a^j_{\wh X/X}-\tau (a^j_\alpha-a^j_{\rm min})- a^j_L\big)[Y_j]
\equiv \cr
\equiv &  K_{\wh X}+ S+ \sum_{j\in J_p}
\big(\tau (a^j_\alpha-a^j_{\rm min})+ a^j_L-a^j_{\rm min}-a^j_{\wh X/X}\big)[Y_j]+
\tau\wh \Lambda_{\alpha}+ \wh \Lambda_{L}\cr
}$$
We use the notation
$$\wh L:= \sum_{j\in J_p}
\big(\tau (a^j_\alpha-a^j_{\rm min})+ a^j_L-a^j_{\rm min}-a^j_{\wh X/X}\big)[Y_j]+
\tau\wh \Lambda_{\alpha}+ \wh \Lambda_{L}$$
and we remark that $\wh L$ is a big $\bR$--line bundle on $\wh X$, whose critical exponent is greater than 1 ; moreover, the restriction $\wh L_{|S}$ has the same properties. 

Next we invoke the non-vanishing theorem 0.1 : as a by-product of its proof, we get the family of approximations $\wh L_\eta$ together with a corresponding family of effective $\bQ$--sections $\wh U_\eta$ of the bundle 
$$K_{\wh X}+ S+ \wh L_\eta\leqno (81)$$
{\sl whose restriction to $S$ is non-zero}, as they were obtained as extensions of non-zero sections defined on $S$.
The bundle $K_{\wh X}+ S+ \wh L$ is obtained as a convex combination of the bundles of type (85), therefore we can assume the existence of a
$\bR$-divisor
$$T:= \sum _{i\in I}\lambda^i [Z_i]\in \{K_{\wh X}+ S+ \wh L\}\leqno (82)$$
such that $\card (I)< \infty$ and such that $T$ is non-singular along $S$.
\smallskip
Then the current
$$\leqno (83)
\eqalign{\wh T:= & T+ \big((1+\tau)a^0_{\rm min}+ a^0_{\wh X/X}\big)[S]+ \cr
+ & \sum_{j\in J_p}\big((1+\tau)a^j_{\rm min}+ a^j_{\wh X/X}\big)[Y_j]+ 
\sum_{j\in J_n}(\tau a^j_\alpha+ a^j_L)[Y_j]\cr
}$$
belong to the class of the bundle 
$$K_{\wh X/X}+  (1+\tau)\mu^\star(K_X+ L),\leqno (84)$$
where the index $j= 0$ in (83) corresponds to $S$. By the Hartogs 
principle, we obtain
$$\wh T= (1+\tau)\mu^\star \Theta+ \sum_{j\in J}a^j_{\wh X/X}[Y_j]\leqno (85)$$
where 
$\Theta$ is an effective $\bR$--divisor in the Chern class of $K_X+ L$. 

\noindent The current $\Theta$ is the one we seek ; it is obvious that 
the Lelong number of $\mu^\star (\Theta)$ along $S$ above is precisely the same as the Lelong number of the inverse image of the current $\Theta_{\rm min}$. The rationality statement in 2.B.1 can be obtained as in [39] : 
we write 
$$\mu^\star(\Theta)= a^0_{\rm min}[S]+ \sum_ja^j_\Theta [W_j]\in \mu^\star \big(c_1(K_X+ L)\big)$$
where $W_j\subset \wh X$ are hypersurfaces and $a^j_\Theta$ are positive real numbers.

\noindent If $a^0_{\rm min}\not \in \bQ$, then we use the rationality of $L$ and infer the existence of an effective 
$\bQ$--divisor
$$b^0S+ \sum_j b^jW_j\in \mu^\star \big(c_1(K_X+ L)\big)$$
such that $b^0< a^0_{\rm min}$ (see [39] for a more complete discussion). The multiplication of the above $\bQ$--section with a 
divisible enough positive integer transform it into a section of
$$m_0(K_X+ L)+ \rho_0$$
where $\rho_0\to X$ is a topologically trivial line bundle. By the theorem 2.A.1, the minimal 
metric on $K_X+ L$ constructed with the sections of multiples of $K_X+ L$ coincide with the minimal metric associated to the family of sections of the multiples of $m_0(K_X+ L)+ \rho_0$, thus we get a contradiction, and the theorem 2.B.1 is completely proved.

 \hfill\qed

\vfill
\eject

\section{References}

\bigskip

{\eightpoint

\bibitem [1]&Bedford, E., Taylor.:& The Dirichlet problem for a complex Monge-Amp\`ere equation~;& Invent. Math.  37  (1976), no. 1, 1--44&

\bibitem [2]&Berndtsson, B.:& On the Ohsawa-Takegoshi extension theorem~;& Ann.\ Inst.\ Fourier (1996)&

\bibitem [3]&Berndtsson, B., P\u aun, M.:& Bergman kernels and the pseudo-effectivity of the relative canonical bundles~;& arXiv:math/0703344, to appear in Duke Math. Journal&

\bibitem [4]&Berndtsson, B., P\u aun, M.:& A Bergman kernel proof of the Kawamata subadjunction theorem, I and II~;& in preparation&

\bibitem [5]&Birkar, C., Cascini, P., Hacon, C., McKernan, J.:&\ Existence of minimal models for varieties of log general type~;& on the web&

\bibitem [6]&Bonavero, L.:& In\'egalit\'es de morse holomorphes singuli\`eres. (French) [Singular holomorphic Morse inequalities]~;&  J.\ Geom. Anal.  8  (1998),  no. 3, 409--425&

\bibitem [7]&Boucksom, S.& Divisorial Zariski decompositions on compact complex manifolds~;&  Ann.\ Sci. Ecole Norm. Sup. (4)  37  (2004),  no. 1, 45--76&

\bibitem [8]&Campana, F., Peternell, Th.& Geometric stability of the cotangent bundle and the universal cover of a projective manifold~;& arXiv: math/0405093&

\bibitem [9]&Claudon, B.:& Invariance for multiples of the twisted canonical bundle~;& math.AG/0511736, Ann. Inst. Fourier (Grenoble)  57  (2007),  no. 1, 289--300&

\bibitem [10]&Corti, A.:& Flips for 3-folds and 4-folds~;& Oxford Lecture Ser. 
Math. Appl. 35(2007), 189 pp, Oxford Univ. Press&

\bibitem [11]&Demailly, J.-P.:& Regularization of closed positive currents and Intersection Theory~; &J. Alg. Geom. 1 (1992) 361-409&

\bibitem [12]&Demailly, J.-P.:& A numerical criterion for very ample line bundles~; &
J. Differential Geom.  37  (1993),  no. 2, 323--374&

\bibitem [13]&Demailly, J.-P.:&  On the Ohsawa-Takegoshi-Manivel  
extension theorem~;& Proceedings of the Conference in honour of the 85th birthday of Pierre Lelong, 
Paris, September 1997, Progress in Mathematics, Birkauser, 1999&

\bibitem [14]&Demailly, J.-P., Peternell, Th., Schneider, M. :& Pseudoeffective line bundles on compact K\"ahler manifolds~; &
Internat. J. Math.  12  (2001),  no. 6&

\bibitem [15]&Demailly, J.-P.:& K\"ahler manifolds and transcendental techniques in algebraic geometry~;&  Plenary talk and Proceedings of the Internat. Congress of Math., Madrid (2006), 34p, volume I&

\bibitem [16]&Druel, S.:& Existence de mod\`eles minimaux pour les vari\'et\'es de type g\'en\'eral~;& Expos\'e 982, S\'eminaire Bourbaki, 2007/08&

\bibitem [17]&Ein, L., Lazarsfeld, R., Musta\c t\u a, M., Nakamaye, M., Popa, M. :&Asymptotic invariants of base loci ~;&  Ann.\ Inst.\ Fourier (Grenoble)  {\bf 56}  (2006),  no. 6, 1701--1734&

\bibitem [18]&Ein, L., Popa, M.:& Adjoint ideals and extension theorems~;& preprint in preparation~; june 2007&

\bibitem [19]&Hacon, C.,  McKernan, J.:& Boundedness of pluricanonical maps of varieties of general type~;&
 Invent.\ Math.\ Volume {\bf 166}, Number 1 / October, 2006, 1-25&
 
 \bibitem [20]&Hacon, C., McKernan, J.:& On the existence of flips~;
 &\ math.AG/0507597&
 
 \bibitem [21]&Hardy, G.H., Wright, E.M.:& An introduction to the theory of numbers~;&\ Oxford University Press, 1938&
 
 \bibitem [22]&Kawamata, Y.:&\ A generalization of Kodaira-RamanujamÕs vanishing 
theorem~;&\ Math. Ann. 261 (1982), 43Ð46&

\bibitem [23]&Kawamata, Y.:& Pluricanonical systems on minimal algebraic varieties~;&\ Invent. Math.  79  (1985),  no. 3&

\bibitem [24]&Kawamata, Y.:& On the extension problem of pluricanonical forms~;& Contemporary. Math.  241  (1999),  no. 3&

\bibitem [25]&Kawamata, Y.:& Finite generation of a canonical ring~;&\ arXiv:0804.3151&

\bibitem [26]&Kim, Dano.:&Ph.D. Thesis~;& Princeton, 2006&

\bibitem [33]&Koll\'ar, J.,  Mori, S.:&Birational geometry of algebraic varieties.~;& Cambridge University Press, Cambridge, 1998&

\bibitem [28]&Lazarsfeld, R.:& Positivity in Algebraic Geometry~;& Springer, Ergebnisse der Mathematik und ihrer Grenzgebiete&

\bibitem [29]&McNeal, J., Varolin, D.:&Analytic inversion of adjunction: $L\sp 2$ extension theorems with gain~;&  Ann.\ Inst.\ Fourier (Grenoble)  {\bf 57}  (2007),  no. 3, 703--718&

\bibitem [30]&Nadel, A. M.:&\ Multiplier ideal sheaves and Kahler-Einstein metrics of positive scalar curvature~;&\ Ann. of Math. (2) 132 (1990), no. 3, 
549Ð596&

\bibitem [31]&Nakayama, N.:&Zariski decomposition and abundance~;&  MSJ Memoirs  {\bf 14}, Tokyo (2004)&
  
 \bibitem [32]&Ohsawa, T., Takegoshi, K.\ :& On the extension of $L^2$
holomorphic functions~;& Math.\ Z.,
{\bf 195} (1987), 197--204&

\bibitem [33]&P\u aun, M.:&\ Siu's Invariance of Plurigenera: a One-Tower Proof~;&\ preprint IECN (2005), J. Differential Geom. 76 (2007), no. 3, 485Ð493&

\bibitem [34]&Shokurov, V.:& A non-vanishing theorem~;&\ Izv. Akad. Nauk SSSR
(49) 1985&

\bibitem [35]&Skoda, H.:& Sous-ensembles analytiques d'ordre fini ou infini dans $  C\sp{n}$ (French)~;& Bull. Soc. Math. France  100  (1972), 353--408&
  
\bibitem [36]&Siu, Y.-T.:& Analyticity of sets associated to Lelong numbers and the extension of closed positive currents& Invent. Math.  27  (1974), 53--156&
 
\bibitem [37]&Siu, Y.-T.:& Invariance of Plurigenera~;& Inv.\ Math.,
{\bf 134} (1998), 661-673&

\bibitem [38]&Siu, Y.-T.:& Extension of twisted pluricanonical sections with plurisubharmonic weight and invariance of semipositively twisted plurigenera for manifolds not necessarily of general type~;& Complex geometry (G\"ottingen, {\bf 2000}),  223--277, Springer, Berlin, 2002&

\bibitem [39]&Siu, Y.-T.:& A General Non-Vanishing Theorem and an Analytic Proof of the Finite Generation of the Canonical Ring~;& arXiv:math/0610740&

\bibitem [40]&Siu, Y.-T.:&\ Finite Generation of Canonical Ring by Analytic Method~;&\ arXiv:0803.2454&  

\bibitem [41]&Takayama, S:& On the Invariance and Lower Semi--Continuity
of Plurigenera of Algebraic Varieties~;& J. Algebraic Geom.  {\bf 16 } (2007), no. 1, 1--18&

\bibitem [42]&Takayama, S:& Pluricanonical systems on algebraic varieties of general type~;&Invent.\ Math.\
Volume {\bf 165}, Number 3 / September, 2005, 551-587&

\bibitem [43]&Tsuji, H.:& Extension of log pluricanonical forms from subvarieties~;& math.CV/0511342 &
 
\bibitem [44]&Varolin, D.:&  A Takayama-type extension theorem~;&  math.CV/0607323, to appear in Comp.\ Math&

\bibitem [45]&Viehweg, E.:&  Vanishing theorems~;& J. Reine Angew. Math. 335 (1982)&

\bibitem [46]&Yau, S-T.:& On the Ricci curvature of a compact K \"ahler manifold and the complex Monge-Amp\`ere 
equation. I.;& Comm. Pure Appl. Math. 31 (1978), no. 3, 339Ð411&

}

\bigskip
\noindent
(version of July 19, 2008, printed on \today)
\bigskip\bigskip
{\parindent=0cm
Mihai P\u aun,
paun@iecn.u-nancy.fr

\end

\bigskip

\subsection {\S 2.C Induction and extension}
\smallskip
\noindent Very roughly, our goal in this section is to show the existence of a finite family of
global sections of $\cA_k$ such that the restriction to $S$ of the inverse image of the metric induced by these sections minus the {\sl minimal vanishing order} along 
$S$ is equivalent with the restriction to $S$ of the inverse image of the minimal metric
minus the {\sl minimal vanishing order} along 
$S$. In other words, the restriction to $S$ of the (modified) minimal metric has log poles.
\smallskip 
\noindent We use the section $\wh T$ we have produced during the 
preceding subsection to show the following :
\smallskip
\noindent $\bullet$  The threshold $\tau$ is rational ;
\smallskip
\noindent $\bullet$ The restriction of the current 
$$\mu^\star\Theta_{\rm min}-a^0_{\rm min}[S]$$
to $S$ is well defined.
\medskip 
\noindent Indeed, the first bullet is clear since $a^0_{\rm min}$ is a positive rational number,
and 
we have an explicit expression for $\tau$, see ?. The second assertion above is equally clear,
by the definition of the minimal metric as in ? ; the upshot is that we can reduce ourselves to the case of $\bQ$-line bundles, as we will show along the next lines.

By the relations (63), (64) and (65) we infer that :
$$
\leqno (66)
\eqalign{
(1+ \tau)\mu^\star\Theta_{\rm min}+ & \sum_{j\in J}a^j_{\wh X/X}[Y_j]\equiv
K_{\wh X}+ S+ \wh L+ 
\big((1+\tau)a^0_{\rm min}+ a^0_{\wh X/X}\big)[S]+ \cr
+ & \sum_{j\in J_p}\big((1+\tau)a^j_{\rm min}+ a^j_{\wh X/X}\big)[Y_j]+ 
\sum_{j\in J_n}(\tau a^j_\alpha+ a^j_L)[Y_j]\cr
}
$$
and therefore we obtain
$$\leqno (?)
\eqalign{
(1+ \tau)\big(\mu^\star\Theta_{\rm min}-a^0_{\rm min}[S] & -\sum_{j\in J_p}a^j_{\rm min}[Y_j]- {{1}\over {1+\tau}}\sum_{j\in J_n}(\tau a^j_\alpha+ a^j_L-a^j_{\wh X/X})[Y_j]\big)\cr 
 & \equiv K_{\wh X}+ S+ \wh L\cr
 }
 $$
In order to illustrate what will follow, we
assume that for any $j\in J_p$ the coefficient $a^j_{\rm min}$ is a rational, 
and also that the components of the divisor part of $\wh L$ are mutually disjoint.
We denote 
by $\Delta_S$ the effective $\bR$--divisor which is the common part of the restriction of the current in the left hand side of (?) and the components of $\wh L$.

We consider a section $u\in \cA_k$, and we denote by $Z_u$ the divisor associated to 
$\mu^\star u$. Then the $\bQ$-divisor
$$1/kZ_u- a^0_{\rm min}[S] -\sum_{j\in J_p}a^j_{\rm min}[Y_j]- {{1}\over {1+\tau}}\sum_{j\in J_n}(\tau a^j_\alpha+ a^j_L-a^j_{\wh X/X})[Y_j]\leqno (?)$$
is effective, since by the definition (?) we have
$$\tau a^j_\alpha+ a^j_L-a^j_{\wh X/X}\leq (1+\tau)a^j_{\rm min}.$$
On the other hand, the $\bQ$--divisor (?) is more singular than the left hand side of 
(?), so in particular the restriction of the corresponding section to 
$S$ is either identically zero (if the vanishing along $S$ is not minimal), or it ``vanish" on 
$\Delta_S$. By induction, the ring associated to the sections of multiples of $K_S+ \wh L$
which vanish on $\Delta_S$ (see ? for the precise condition) is of finite type.

By the invariance of plurigenera we infer that any sections which belong to this ring admit an extension to $\wh X$, which moreover can be injected into to $\cR_\Delta(X, L)$ by multiplication with the corresponding multiple of the divisor 
$$a^0_{\rm min}[S] +\sum_{j\in J_p}a^j_{\rm min}[Y_j]+ {{1}\over {1+\tau}}\sum_{j\in J_n}(\tau a^j_\alpha+ a^j_L-a^j_{\wh X/X})[Y_j]$$
(see the inequalities (?) and (?)). In conclusion, the following restriction 
$$\mu^\star\Theta_{\rm min}-a^0_{\rm min}[S]_{|S}$$
has the same singularities as the algebraic metric induced by the generators of the ring
$\displaystyle \cR_{\Delta_S}(S, \wh L)$, and this is what we want to establish in this subsection. \hfill\qed
\medskip 
\noindent We will follow this program in the general setting--as it is no obvious for us that the additional assumptions above should hold {\sl a-priori}. 
\smallskip
We recall the expression of $\wh L$ in (?) ; it clearly show that 
$$(1+\tau)\sum_{j\in J_p}a^j_{\rm min}[Y_j]+ \wh L$$
is a $\bQ$--bundle ; therefore, for any positive sequence of real numbers 
$(\eta^j)_{j\in J_p}$ such that 
$$(1+\tau)a^j_{\rm min}-\eta^j\in \bQ$$
the Chern class of the bundle 
$$\wh L(\eta):= \sum_{j\in J_p}\eta^j[Y_j]+ \wh L\leqno (?) $$
is rational. The critical exponent of the current (?) is greater than 1, and we will choose the 
parameters $(\eta^j)$ such that the $\bQ$--bundle (?) has the same property ; we will impose more conditions on the set of approximations later, as they will appear along the following lines.

The equivalence (?) become

$$\leqno (?)
\eqalign{
(1+ \tau)\big(\mu^\star\Theta_{\rm min}-a^0_{\rm min}[S] &
-\sum_{j\in J_p}\big(a^j_{\rm min}-{{1}\over {1+\tau}} \eta^j\big)[Y_j]- \cr 
& -{{1}\over {1+\tau}}\sum_{j\in J_n}(\tau a^j_\alpha+ a^j_L-a^j_{\wh X/X})[Y_j]\big) 
 \equiv K_{\wh X}+ S+ \wh L_\eta.\cr
 }
 $$
At this point we use again the lemma of Hacon-McKernan [20], to infer the existence of
a composition of blow-up maps $\mu_1: X_1\to \wh X$ such that 
$$\mu_1^\star\big(K_{\wh X}+ S+ \wh L(\eta)\big)+ E\equiv K_{X_1}+ S_1+ \wt L\leqno (?)$$
where the notations in the previous formula are as follows.
\smallskip
\noindent $\bullet$ $E$ is an exceptional $\bQ$-divisor.
\smallskip
\noindent $\bullet$ $S_1$ is the proper transform of $S$ ; in particular, $S_1$ is not 
$\mu_1$--exceptional.
\smallskip
\noindent $\bullet$ $\wt L$ is a $\bQ$--line bundle, such that 
$$\wt L= \sum_{j\in J_1}\nu^j[W_j]+ \Lambda_1\leqno (?)$$
where $0<\nu^j< 1$ are rational numbers (which depend on the perturbation $(\eta^j)$ we use)
and such that $W_j\cap W_k= \emptyset$ for any $j\neq k$.
We have 
$$\Lambda_1:= \mu_1^\star(\tau\wh \Lambda_{\alpha}+ \wh \Lambda_{L})\leqno (?)$$
and we remark that this differential form is independent of $(\eta^j)$ (see the equality (?))~;  one can see that so is the map $\mu_1$.

We denote by $T_1$ the inverse image of the current in the left hand side of (?) plus $E$, so that we have
$$T_1\equiv K_{X_1}+ S_1+ \wt L.\leqno (?)$$
The main properties of $T_1$ which will be useful for us are listed below.
\smallskip
{\itemindent 9.5 mm
\item {$\bf 2.C.1$} For each $j\in J_p$, we denote by $\ol Y_j$ the proper transform of $Y_j$ under the map $\mu_1$. Then we have
$$T_1\geq \sum _{j\in J_p}\eta^j[\ol Y_j]\leqno (?)$$
and moreover we remark that the singular part of $\wt L$ has similar properties, i.e.
$$\sum_{j\in J_1}\nu^j[W_j]\geq \sum _{j\in J_p}\eta^j[\ol Y_j]\leqno (?)$$
\smallskip
\item {$\bf 2.C.2$} The current $T_1$ admits a well-defined restriction to $S_1$ ; we write
$$T_{1|S_1}:= \sum _{j\in J_1}\rho^j[W_{j|S_1}]+ R_1\leqno (?)$$
\hfill\qed
}

\medskip
\noindent Indeed, the relation (?) is a consequence of the definition of $T_1$, and of the expression (?) ; as for 
(?), it holds because the support of the difference between $\wt L$ and $\mu^\star \big(\wh L(\eta)\big)$ is $\mu_1$-exceptional. The assertion 2.C.2 is a consequence of the second bullet
at the beginning of 2.C.

\smallskip
We introduce the following notation
$$\Delta_{S_1}:= \sum_{j\in J_1}\min (\rho^j, \nu^j)[W_{j|S_1}] ;\leqno (?)$$
it is an effective $\bR$-divisor, and it represents the common part of the singularities of
$\displaystyle T_{1|S_1}$ and $\displaystyle \wt L_{|S_1}$ in codimension 1. The motivation of considering precisely this divisor is clear : it appear as ``obstruction" to the lifting of the sections
of $m(K_S+ \wt L_{|S})$ to $X_1$.

We consider the ... the rest tonight!

Therefore by induction, the ring
$$R(S, \wh L):= \bigoplus _{m\in m_0\bZ}H^0\big(S, m(K_S+  \wh L)\big)\leqno (66)$$
is of finite type. Moreover we have the next statement.

\claim 1.I.2 Lemma|Any holomorphic section of the bundle $m(K_S+ \wh L_{|S})$ extends to $\wh X$.
\endclaim

\proof.

\vskip 3cm

By the method explained 
in [?], page ?, we can assume that all the coefficients of $\Theta$ are rational ; here is where we use the fact that $L$ is a $\bQ$-line bundle ({\sl 
some additional explanations at this point, saying that if if the coefficient corresponding to $S$ is not rational, we can decrease it, and get a contradiction}). \hfill\qed
\smallskip



\section{\S 2. Metrics with minimal singularities and holomorphic sections}

\medskip
\noindent In this section we will consider the following geometric context. Let $X$ be a non-singular, projective $n$-dimensional manifold, and let 
$L$, respectively $\Delta$ be a $\bQ$-line bundle, respectively an effective 
$\bR$-divisor such that :
\smallskip
\item {($\bullet$)} The (1,1)--class $c_1(L)$ contains a K\"ahler current 
$\Theta_L$ ; 
\smallskip
\item {($\bullet$)} We have $\Theta_L\geq [\Delta]$ and the critical exponent of the difference $\Theta_L- [\Delta]$ is strictly greater than 1.
\medskip
\noindent We consider the following graded algebra
$$\cR_\Delta(X, L):= \bigoplus_{k\in q\bZ_+}\cA_k\leqno (47)$$
where we denote by 
$$\cA_k:= \big\{u\in H^0\big(X, k(K_X+ L)\big) : |u|_k^2\leq C\exp(k\varphi_\Delta)\big\} ;$$
the weight $\varphi_\Delta$ is the potential of the current $[\Delta]$ and
the inequality above is supposed to hold pointwise on $X$. 
We assume that a non-singular metric $\wt h$ on $K_X+ L$ is given, and 
$|u|_k$ above is the norm of $u$ with respect to the metric induced
on $k(K_X+ L)$.
The positive constant 
$C$ is dependent on $u$, and the integer $q$ in the definition of the algebra above is such that $qL$ become a genuine line bundle.

Strictly speaking, the inequality in the definition of $\cA_k$ makes no sense, i.e. it is not intrinsic, but this can be easily arranged by considering a function on $X$ which 
have precisely the same singularities as $\varphi_\Delta$. We prefer to state it like this, for the simplicity of notation ; 
its meaning is to impose a restriction on the zeroes of $u$ according to the size of $\Delta$. \hfill\qed
\medskip

In the paragraph 2.A below we 
define {\sl a metric with 
minimal singularities} which is adapted to the ring $\cR_\Delta(X, L)$. For some technical reasons (that will only appear in the second paragraph 2.B), we are forced to take into account the sections of the multiples of $K_X+ L$
{\sl twisted with a topologically trivial line bundle},
even if our ultimate goal would be to understand the structure of the ring $\cR_\Delta(X, L)$.

Next, we will consider the relative threshold of the minimal metric with respect to a metric given by a finite number of sections of the above ring ; the non-vanishing theorem will provide us with an
$\bR$--section of $K_X+L$ 
which has precisely the same vanishing order along {\sl some} divisor, say $S$, of a modification of $X$ as the metric with minimal singularities (compare with [39], [40]). 
The results we obtain in the paragraph 2.A will show that the said vanishing order is
a rational number. \hfill\qed

\subsection {\S 2.A  Metrics with minimal singularities}
\smallskip
\noindent Along the following lines, we will only consider the case of adjoint bundles, since it is in this setting that the main properties we will establish afterwards hold ; 
however, one could define the objects below in a more general context.

Let $\rho\to X$ be a topologically trivial line bundle, endowed with a metric
$h_\rho$ whose curvature form is equal to zero. Given 
$$u\in H^0\big(X, k(K_X+ L)+\rho\big)$$
we denote by $|u|_{k,\rho}^2$ the poinwise norm of $u$, measured with the 
$k^{\rm th}$ power of $\wt h$ twisted with $h_\rho$.

We denote by 
$$\cA_{k, \rho}:= \big\{u\in H^0\big(X, k(K_X+ L)+\rho\big) : |u|_{k,\rho}^2\leq C\exp(k\varphi_\Delta)\big\} ;$$
it is a $\rho$-twisted version of the set $\cA_{k}$ previously considered.

\noindent We introduce the following class of functions on $X$
$$\cF_\Delta:= \Big\{f= {1\over k}\log|u|_{k,\rho}^2,  
\hbox{ where } k\in \bZ_+, u\in \cA_{k, \rho}, \hbox{ such that }\sup_X|u|_{k,\rho}^2= 1  \Big\},$$
and we remark that we have 
$$\Theta_{\wt h}(K_X+ L)+ {{\sqrt {-1}}\over {2\pi}}\ddbar f\geq 0$$
for any $f\in \cF_\Delta$ (since the curvature of $\rho$ with respect to $h_\rho$ is equal to zero).
Thus the curvature current associated to the metric 
$$\exp(-f)\wt h$$
on the bundle $K_X+ L$ is positive.

We will denote by $h_{\rm min}$ the metric on $K_X+ L$ given 
by the smallest upper semicontinuous
majorant of the family $\cF_\Delta$ above. We denote by 
$$\Theta_{\rm min}:= \Theta_{h_{\rm min}}(K_X+ L)$$
the curvature current associated to this metric ; by definition we have
$$|v|^2\exp(-k\varphi_{\rm min})\leq O(1)$$
for any $v\in \cA_{k, \rho}$. 
\medskip

\noindent Next {\sl we fix} a topologically trivial line bundle $\rho_0$ on $X$ and an integer 
$m_0\in \bZ_+$, and we will construct another metric on the bundle $K_X+ L$
as follows. We consider the set of potentials
$$\cF_\Delta^{\rho_0}:= \Big\{f= {{1}\over {km_0}}\log|u|_{km_0,k\rho_0}^2,  \hbox{ where } k\in \bZ_+, u\in \cA_{km_0, k\rho_0}, \sup_X|u|_{km_0,k\rho_0}= 1 \Big\}$$
and we denote by $h_{\rm min}^{\rho_0}$ the metric on $K_X+ L$ given 
by the smallest upper semicontinuous
majorant of the family $\cF_\Delta^{\rho_0}$ ; let $\Theta_{\rm min}^{\rho_0}$ be the associated curvature current.
\medskip
\noindent We state now the main result of this subsection.

\claim 2.A.1 Theorem|We have 
$$h_{\rm min}= h_{\rm min}^{\rho_0}.$$
\endclaim
\proof \hskip 1.5mm (of the theorem 2.A.1). The relation 
$$\varphi_{\rm min}\geq\varphi_{\rm min}^{\rho_0}$$
is implied by the definition ; in order to obtain an inequality in the opposite sense, 
we consider a section
$$u_1\in \cA_{m_1, \rho_1}.\leqno (48)$$ 

\noindent As a consequence of an argument due to Shokurov 
already employed during the section 1.G, we have the next statement.
\claim 2.A.2 Lemma|For any $k\in \bZ_+$,  there exist a section
$$u_k\in \cA_{km_1m_0, km_1\rho_0}\leqno (49)$$
such that the next integral condition is satisfied
$$\int_X{{c_ku_k\wedge \ol {u_k}\exp({-km_1\varphi_{\rho_0}})}\over 
{{(c_1u_1\wedge \ol{u_1}})^{{{km_1m_0-1}\over{m_1}}}\exp\big(-
{{km_1m_0-1}\over {m_1}}\varphi_{\rho_1}\big)}}\exp(-\varphi_L)d\lambda
= 1.\leqno (50)$$
\hfill\qed
\endclaim
\noindent  In the statement above we denote by $\varphi_L$ the local weights of any metric on 
$L$ which satisfy the properties at the beginning of 2.A : the corresponding curvature is a K\"ahler current which dominates $[\Delta]$ and the critical exponent of the difference is greater than 1. The quantities $c_k$ and $c_1$ are the usual ones, such that the wedge at the
denominator and numerator in the previous lemma are reals. We remark that the quantity under the integral sign in the above formula is a globally defined measure.
\smallskip
\noindent 
\proof. For any positive integer $k$, we consider the bundles
$$km_1m_0(K_X+ L)+ km_1\rho_0\leqno (51)$$
and
$$km_1m_0(K_X+ L)+ km_0\rho_1\leqno (52)$$
as well as the multiplier ideal
$$I_k:= \cI\big((km_1m_0-1)\log|u_1|^{{2}\over {m_1}}_{m_1, \rho_1}+ \varphi_L\big)\leqno (53)$$
We have
$$\chi\Big(\big(km_1m_0(K_X+ L)+ km_1\rho_0\big)\otimes I_k\Big)=
\chi \Big(X, \big(km_1m_0(K_X+ L)+ km_0\rho_1\big)\otimes I_k\Big)$$
by the usual arguments (i.e. the existence of a finite and free resolution of the multiplier sheaf $I_k$, and the 
fact that $\rho_j$ are topologically trivial).
Thanks to the Kawamata-Viehweg-Nadel vanishing theorem we equally know that the respective higher cohomology groups are equal to zero, so in conclusion
$$ H^0\Big(X, \big(km_1m_0(K_X+ L)+ km_1\rho_0\big)\otimes I_k\Big)=
H^0\Big(X, \big(km_1m_0(K_X+ L)+ km_0\rho_1\big)\otimes I_k\Big).\leqno (54)$$

\noindent We claim now that the section $u_1^{\otimes km_0}$ belong to the right hand side cohomology group in (54). To verify this claim, we have to show that the following integral converge
$$\int_X|u_1|^{{{2}\over {m_1}}}\exp (-\varphi_L)d\lambda,$$
and indeed this is the case, since $u_1\in \cA_{m_1, \rho_1}$ and by definition we have
$$|u_1|^{{{2}\over {m_1}}}\leq C\exp(\varphi_\Delta).$$
Since the critical exponent of the current $\Theta_L-[\Delta]$ is greater than 1, we get
$$\int_X\exp(\varphi_\Delta- \varphi_L)d\lambda< \infty$$
and thus the dimension of the right hand side vector space in (54) is strictly positive.
\medskip 
\noindent Therefore we infer the existence of a non-identically zero section 
$$u_k\in H^0\Big(X, \big(km_1m_0(K_X+ L)+ km_1\rho_0\big)\otimes I_k\Big).\leqno (55)$$
By the construction of the ideal $I_k$, we see that we can normalize the section $u_k$ such that the integral condition (50) is satisfied ; we show now that $u_k$ has the vanishing properties along $\Delta$ as required.
 
The finiteness of (50) show the existence of a section 
$$v_k\in H^0\big(X, m_1(K_X+ L)+ km_1\rho_0-km_0\rho_1\big)$$
such that 
$$u_k= u_1^{\otimes (km_0-1)}\otimes v_k.\leqno (56)$$
The integral relation (50) become
$$\int_X{{v_k\wedge \ol v_k}\over {{\big(u_1\wedge\ol u_1\big)}^{{{m_1-1}\over{m_1}}}}}
\exp\big(-\varphi_L-km_1\varphi_{\rho_0}+ (km_0-1+ 1/m_0)\varphi_{\rho_1}\big)= 1.
\leqno (57)$$

\medskip
\noindent In order to obtain ``enough vanishing" of $v_k$ along $\Delta$, we will construct
next a slight perturbation of the metric on $L$. 

To this end, we decompose the divisor $\Delta$ as follows
$$\Delta= \sum_{j\in Q}\theta^j[Y_j]+ \sum_{j\in R}\theta^j[Y_j]$$
where the coefficients corresponding to the index set $Q$ are positive rational numbers, and the coefficients corresponding to $R$ belong to $\bR_+\!\!\setminus \bQ$.

Our next observation is that 
we can assume that $m_1$ is large and divisible enough. Indeed,
we can always consider an appropriate power of $u_1$, without affecting the normalization in the definition of the metric $\varphi_{\rm min}$.

We consider the next decomposition
$$L= \Delta+ (L-\Delta) ;$$
the Chern class of the $\bR$--bundle $L-\Delta$ contains the K\"ahler current
$$\Theta_L-[\Delta] 
,\leqno (58)$$ 
whose critical exponent is greater than 1, by hypothesis. We endow the component 
$\Delta$ in the previous decomposition with a metric whose curvature form is equal to
the next expression
$$\sum_{j\in Q}\theta^j[Y_j]+ \sum_{j\in R}{{1+ [m_1\theta^j]}\over {m_1}}[Y_j]+ \Theta (m_1)\leqno (59)$$
where $\Theta (m_1)$ is a non-singular (1,1)--form in the cohomology class of
the divisor
$$ \sum_{j\in R}\Big(\theta^j- {{1+ [m_1\theta^j]}\over {m_1}}\Big)[Y_j].$$
We remark that the cohomology 
class as the current in (58) is the same as the one of $[\Delta]$. We will assume that the real numbers
$$\theta^j- {{1+ [m_1\theta^j]}\over {m_1}}$$
are small enough, so that the curvature of the bundle $L$ associated to the metric 
induced by (58) and (59) is positive, and such that the critical exponent of the difference between this current and $[\Delta]$ is greater than 1. It is at this point that we may be forced to consider a 
large enough multiple of the initial $m_1$. 
\smallskip
We suppose that the metric $\varphi_L$ we have worked with during the construction of $u_k$  
was precisely the one obtained above. Then we infer the following relations.

\smallskip
\noindent $\bullet$ Let $j\in Q$ ; we denote by $f_j$ the local equation of the hypersurface $Y_j$ and we remark that we have
$$|u_1|^2\leq C|f_j|^{2m_1\theta^j}\leqno (60)$$
by hypothesis. But then we have 
$$|u_1|^{2(1-1/m_1)}|f_j|^{2\theta^j}\leq C|f_j|^{2m_1\theta^j}\leqno (61)$$
and we see that (57) coupled with (61) imply the next inequality
$$|v_k|^2\leq C|f_j|^{2m_1\theta^j}.\leqno (62)$$
because we can assume further that $m_1$ is divisible by all the 
denominators of the coefficients $\theta^j$ for $j\in Q$.
The relations (56) and (62) show that 
$$|u_k|^2\leq C|f_j|^{2km_1m_0\theta^j}$$
and as a consequence, the section $u_k$ has the desired vanishing properties along
the components of $\Delta$ with rational coefficients.

\smallskip
\noindent $\bullet$ Let $j\in R$ ; by hypothesis we have
$$|u_1|^2\leq C|f_j|^{2(1+ [m_1\theta^j])}\leqno (63)$$
and the additional singularities of $\varphi_L$ along $Y_j$ as in (59),
 together with the 
finiteness of the integral in (57)  
show that the zero set of $v_k$ contains the divisor
$$(1+ [m_1\theta^j])Y_j ;$$
so that we have
$$|u_k|^2\leq C|f_j|^{2km_0(1+ [m_1\theta^j])}$$
and therefore the second case is equally settled. 
\smallskip
\noindent In conclusion, the section $u_k$ belong to the set $\cA_{km_1m_0, km_1\rho_0}$ and the lemma 2.A.2 is proved. \hfill\qed
\medskip
We will use the family of sections $(u_k)_{k\in \bZ_+}$ in order to compare 
$\varphi_{\rm min}$ and $\varphi_{\rm min}^{\rho_0}$. A specific normalization was chosen for the sections defining the potentials in $\cF _\Delta$ ; thus, we
have to estimate the sup norm of $u_k$ along the next lines. Our main technical tools
will be the standard convexity properties of the psh functions.

We denote by 
$$\exp(f_k):= |v_k|^2_{m_1, km_1\varphi_{\rho_0}- (km_0-1)\varphi_{\rho_1}} ;\leqno (64)$$
the next step in our proof is to show the existence of a positive constant 
$C= C(m_1)$ large enough,  independent of $k$, so that we have 
$$\sqrt {-1} \ddbar f_k\geq -C\omega\leqno (65)$$
and moreover
$$-\log C\leq \max_Xf_k\leq \log C.\leqno (66)$$
The inequality (65) is clear ; let us give some explanations about (66). By using the 
notations introduced in (64), the equality (57) become
$$\int_X{{\exp\big(f_k-f_L)}\over {|u_1|_{m_1, \rho_1}^{2{{m_1-1}\over{m_1}}}}}dV_\omega
= 1.
\leqno (67)$$
where $f_L$ is the (global) distortion function between the
metric $\varphi_L$ and the non-singular metric on $L$ induced by $\wt h$ on $K_X+ L$ and 
$\det (\omega)$ on $-K_X$.
The section $u_1$ is normalized such that
$$\max_X(|u_1|_{m_1, \rho_1})= 1$$
and then the second inequality of (66) is a consequence of the mean inequality
for the psh functions. In order to obtain the first relation, we consider a log-resolution 
$\mu: \wh X\to X$ of $\displaystyle f_L+ {{m_1-1}\over{m_1}}\log |u_1|_{m_1, \rho_1}^2$. The finiteness of the integral (67) and the change of variable formula will lead us to the following
situation
$$\int_{\wh X}\exp \big(\wh f_k-\psi\big)dV= 1\leqno (68)$$ 
where $\psi$ is a quasi-psh function on $\wh X$ with trivial multiplier sheaf ; we equally have
$\wh f_k:= f_k\circ \mu+ \sum_ja^j\log|s_j|$
where $s_j$ are global sections, and $a^j$ may also be negative.
As is the previous context on $X$, the hessian of
$\wh f_k$ is bounded from below by a uniform constant ; we stress on the fact that
the objects $\psi, a^j, s_j$ {\sl do not depend on $k$}.
 
Let 
$$\wh C_k:= \max_{\wh X} \wh f_k ;$$
the relation (68) become
$$\int_{\wh X}\exp \big(\wh f_k-\wh C_k -\psi\big)dV
= \exp(-\wh C_k).
\leqno (69)$$
We already know that the sequence $\wh C_k$ is bounded from above ;
the relation (69) show that 
$$\exp(-\wh C_k)\leq \int_{\wh X}\exp \big(-\psi\big)dV$$
so our sequence is bounded from below as well. 

But then the functions $(f_k)$ have similar properties ; this can be seen as follows.
We can assume that $\wh C_k= 0$, by the preceding discussion. Let $x_k\in \wh X$
be a point such that $\wh f_k(x_k)= 0$. The mean value property of $f_k$ show that
$$\int_{B(x_k, 1)}\wh f_k dV\geq -C\leqno (70)$$
where $C$ above is independent of $k$ (but eventually depending on $m_1$). By the negativity of the integrant in (70), we have
$$\int_{\Omega}\wh f_k dV\geq -C\leqno (71)$$
for any open set $\Omega\subset B(x_k, 1)$. We choose $\Omega$ such that it does not intersect
the inverse image of the singularities of the function $f_L+ {{m_1-1}\over{m_1}}\log |u_1|^2$, and then the relation (66) is completely proved. Indeed, even if the point $x_k$ varies,
the volume of the set $B(x_k, 1)$ has uniform bounds, and the same thing is true for the
set $\Omega$.\hfill\qed
\smallskip
The important consequence of the previous considerations is the existence of a limit for the sequence $(f_k)$.
The normalization of $u_1$ and the relation (56) show that 
$$\max_X|u_k|_{km_1m_0, km_1\rho_0}\leq e^{C(m_1)}$$
and finally we get
$$\eqalign{
\varphi_{\rm min}^{\rho_0}\geq & -{{C(m_1)}\over {km_1m_0}}+ {{1}\over {km_1m_0}}\log |u_k|^2\geq\cr
\geq & -{{C(m_1)}\over {km_1m_0}}+ {{km_0-1}\over {km_1m_0}}\log |u_1|^2+ 
{{1}\over {km_1m_0}}f_k.\cr
}$$
We let $k\to \infty$ ; the first and the third term in the last inequality above tend to zero, and thus we get
$$
\varphi_{\rm min}^{\rho_0}\geq {{1}\over {m_1}}\log |u_1|^2.$$
The section $u_1$ above is arbitrary, thus the theorem is proved.
\hfill\qed

\bigskip

\subsection {\S 2.B Constructing sections with minimal vanishing order}
\medskip
\noindent In this paragraph we would like to point out an important property
of the zeroes of the $\bR$--sections produced by 0.1, in connection with Siu's 
approach of the finite generation problem (see [39], [40]). The same hypothesis/conventions
as in the beginning of the section are in force ; in addition, 
given an integer $\alpha$ large enough, we consider the following truncation of the minimal metric
$$\varphi_{\alpha}:= \log \Big(\sum_{k=1}^\alpha\varepsilon_k\!\sum_{j\in J_k}|f^k_j|^{2\over k}\Big)\leqno (72)$$
where $\varepsilon_k$ are positive real numbers, and
$(f^k_j)$ are local expressions of a family of sections of $k(K_X+ L)+ \rho_k$, where 
$\rho_k$ is a topologically trivial line bundle (which may vary with $k$)~;
let $\Theta_{\alpha}$ be the corresponding current.

Precisely as in the paragraph 1.C, we will consider $\mu: \wh X\to X$
a log-resolution of the currents $\Theta_L$ and $\Theta_{\alpha}$~; we have

$$\mu^\star(\Delta)= \sum_{j\in J}a^{j}_\Delta[Y_{j}]\leqno (73)$$
and
$$\mu^\star(\Theta_L)= \sum_{j\in J}a^{j}_L[Y_{j}]+ \wh \Lambda_L\leqno (74)$$
as well as 
$$\mu^\star(\Theta_{\alpha})= \sum_{j\in J}a^{j}_\alpha [Y_{j}]+ \wh \Lambda_\alpha\leqno (75)$$
where $\wh \Lambda_L$, respectively $\wh \Lambda_\alpha$ are non-singular and semi-positive $(1,1)$--forms on $\wh X$ which are positively defined at the generic point of this manifold.
We consider next the inverse image of the minimal current via $\mu$ :
$$\mu^\star(\Theta_{\rm min})= \sum_{j\in J}a^{j}_{\rm min} [Y_{j}]+ 
\wh \Lambda_{\rm min}\leqno (76)$$
where $\wh \Lambda_{\rm min}$ is a closed positive current, whose generic Lelong numbers along the hypersurfaces $Y_j$ above is equal to zero.
Moreover, by the definition of the minimal metric we have 
$$a^j_\alpha\geq  
a^{j}_{\rm min}\leqno (77)$$
for all $j\in J$. We equally have the pointwise
inequality 
$$\varphi_{\wh \Lambda_{\rm min}}\geq \sum_{j\in J}\big(a^j_\alpha- 
a^{j}_{\rm min})\log |f_j|^2\leqno (78)$$
modulo an irrelevant constant. 

We recall that by hypothesis we have $\Theta_L\geq [\Delta]$ ; also, the sections we use in the 
construction of the metrics $\varphi_{\rm min}, \varphi_\alpha$
vanish along $\Delta$ according to their respective definition. Therefore we have
$$a^j_\alpha\geq a^j_{\rm min}\geq a^j_\Delta\leqno (79)$$
and also
$$a^j_L\geq a^j_\Delta.\leqno(80)$$
\medskip
\noindent The main result of the current subsection is the following.

\claim 2.B.1 Theorem|If at least one of the inequalities (77) is strict, then
there exist a section $u\in \cA_{m, \rho}$ and an index $j_0\in J$
such that the vanishing order of $\mu^\star(u)$ along 
$\displaystyle Y_{j_0}$ is precisely $ma^{j_0}_{\rm min}$. \hfill\qed
\endclaim

\medskip

\noindent We remark that if all the inequalities (77) are 
{\sl equalities}, then (78) show that all the local potentials of the current 
$\wh \Lambda_{\rm min}$ are bounded. In other words, the metric with minimal singularities is equivalent with its truncation $\varphi_\alpha$, and this 
imply the finite generation of the ring associated to $K_X+ L$, according to [39].
\medskip
\proof.  We consider the relative threshold associated to the following objects :
$$\tau:= \sup\{t\in \bR_+:  \int_{\wh X}\exp\big(t(\varphi_{\wh D}-\varphi_\alpha
)+ \varphi_{\wh D}+ \varphi_{\wt X/X}- \varphi_L\circ \mu\big)d\lambda<\infty \}.\leqno (81)$$We observe that $\tau$ verify the next relations
$$0< \tau< \infty\leqno (82)$$
by the same arguments as in the proof of 0.1--we remark that the latter inequality is a consequence of our assumption above.

The perturbation argument we have used in 1.C still apply in the present setting ; there exist a unique $S\subset (Y_j)$ such that we have the next relation
$$\leqno (83)
\eqalign{
\mu^\star\big( & K_X+ \tau(\Theta_\alpha-\Theta_{\rm min})+ 
\Theta_L- \Theta_{\rm min}\big)+ (1+\tau) \wh \Lambda_{\rm min}
\equiv \cr 
\equiv &K_{\wh X}+ S+ \sum_{j\in J}
\big(\tau (a^j_\alpha-a^j_{\rm min})+ a^j_L-a^j_{\rm min}-a^j_{\wh X/X}\big)[Y_j]+
\tau\wh \Lambda_{\alpha}+ \wh \Lambda_{L}\cr
}$$
where the coefficients of $Y_j$ are strictly smaller than 1, and the form 
$\wh \Lambda_{L}$ is positively defined. 
The relation (83) is equivalent with
$$\leqno (84)
\eqalign{
(1+\tau) \wh \Lambda_{\rm min}+ &
\sum_{j\in J_n}
\big(a^j_{\rm min}+a^j_{\wh X/X}-\tau (a^j_\alpha-a^j_{\rm min})- a^j_L\big)[Y_j]
\equiv \cr
\equiv &  K_{\wh X}+ S+ \sum_{j\in J_p}
\big(\tau (a^j_\alpha-a^j_{\rm min})+ a^j_L-a^j_{\rm min}-a^j_{\wh X/X}\big)[Y_j]+
\tau\wh \Lambda_{\alpha}+ \wh \Lambda_{L}\cr
}$$
We use the notation
$$\wh L:= \sum_{j\in J_p}
\big(\tau (a^j_\alpha-a^j_{\rm min})+ a^j_L-a^j_{\rm min}-a^j_{\wh X/X}\big)[Y_j]+
\tau\wh \Lambda_{\alpha}+ \wh \Lambda_{L}$$
and we remark that $\wh L$ is a big $\bR$--line bundle on $\wh X$, whose critical exponent is greater than 1 ; moreover, the restriction $\wh L_{|S}$ has the same properties. 

Next we invoke the non-vanishing theorem 0.1 : as a by-product of its proof, we get the family of approximations $\wh L_\eta$ together with a corresponding family of effective $\bQ$--sections $\wh U_\eta$ of the bundle 
$$K_{\wh X}+ S+ \wh L_\eta\leqno (85)$$
{\sl whose restriction to $S$ is non-zero}, as they were obtained as extensions of non-zero sections defined on $S$.
The bundle $K_{\wh X}+ S+ \wh L$ is obtained as a convex combination of the bundles of type (85), therefore we can assume the existence of a
$\bR$-divisor
$$T:= \sum _{i\in I}\lambda^i [Z_i]\in \{K_{\wh X}+ S+ \wh L\}\leqno (86)$$
such that $\card (I)< \infty$ and such that $T$ is non-singular along $S$.
\smallskip
Then the current
$$\leqno (87)
\eqalign{\wh T:= & T+ \big((1+\tau)a^0_{\rm min}+ a^0_{\wh X/X}\big)[S]+ \cr
+ & \sum_{j\in J_p}\big((1+\tau)a^j_{\rm min}+ a^j_{\wh X/X}\big)[Y_j]+ 
\sum_{j\in J_n}(\tau a^j_\alpha+ a^j_L)[Y_j]\cr
}$$
belong to the class of the bundle 
$$K_{\wh X/X}+  (1+\tau)\mu^\star(K_X+ L),\leqno (88)$$
where the index $j= 0$ in (86) corresponds to $S$. By the Hartogs 
principle, we get 
$$\wh T= (1+\tau)\mu^\star \Theta+ \sum_{j\in J}a^j_{\wh X/X}[Y_j]\leqno (89)$$
where in the formula (89) we use the notations (16) for the relative canonical bundle, and
$\Theta$ is an effective $\bR$--section of $K_X+ L$. 

\noindent The $\bR$--section $\Theta$ is the one we seek ; one can easily see that 
the coefficient of $\mu^\star (\Theta)$ along $S$ above is precisely the same as the one of the inverse image of the minimal current. The rationality statement in 2.B.1 can be obtained as in [39] : we remark that we use at this point the 
theorem 2.A.1. 

Moreover, the inequalities (79) and (80) show that 
$$\wh T\geq (1+ \tau)\mu^\star \Delta$$
and therefore the theorem is completely proved. \hfill\qed

\end